\documentclass{article}

\usepackage{Jared-MK-Thin}
\usepackage{graphicx} 
\usepackage{empheq}

\makeatletter
\renewcommand\subparagraph{%
\@startsection {subparagraph}{5}{\z@ }{3.25ex \@plus 1ex
\@minus .2ex }{-1em}{\normalfont \normalsize \bfseries }}%
\makeatother

\usepackage{thmtools}
\usepackage{thm-restate}
\usepackage{cleveref}

\newtheorem{Assumption}{Assumption}
\newtheorem{Remark}{Remark}

\newcommand{\RV}{\mathcal{V}}
\newcommand{\bnu}{\overline{\nu}}

\newcommand{\bg}{\overline{g}}
\newcommand{\bh}{\overline{h}}
\newcommand{\bT}{\overline{T}}

\newcommand{\bGamma}{\overline{\Gamma}}
\newcommand{\eps}{\epsilon}
\newcommand{\p}{\partial}

\newcommand{\FPz}{\underset{z=0}{FP} \; }
\newcommand{\FPeps}{\underset{\epsilon\to 0}{FP} \;}

\newcommand{\Fa}{F_{\alpha}}

\newcommand{\va}{v_{\alpha}}
\newcommand{\vb}{v_{\beta}}
\newcommand{\vg}{v_{\gamma}}
\newcommand{\vd}{v_{\delta}}
\newcommand{\vx}{v_x}
\newcommand{\vs}{v_{\sigma}}
\newcommand{\vt}{v_{\tau}}
\newcommand{\vo}{v_{\omega}}
\newcommand{\F}{\mathcal{P}}

\newcommand{\bW}{\overline{W}}
\newtheorem*{theorem*}{Theorem}
\newtheorem*{proposition*}{Proposition}
\newcommand{\bA}{\overline{A}}
\newcommand{\bnabla}{\overline{\nabla}}
\newcommand{\tbGamma}{\tilde{\overline{\Gamma}}}
\newcommand{\bN}{\overline{N}}

\title{Variations of Renormalized Volume for Minimal Submanifolds of Poincare-Einstein Manifolds}
\author{Jared Marx-Kuo}
\date{}

\begin{document}

\maketitle

\begin{abstract}
\noindent We investigate the asymptotic expansion and the renormalized volume of minimal submanifolds, $Y^m$ of arbitrary codimension in Poincare-Einstein manifolds, $M^{n+1}$. In particular, we derive formulae for the first and second variations of renormalized volume for $Y^m \subseteq M^{n+1}$ when $m < n + 1$. We apply our formulae to the codimension $1$ and the $M = \H^{n+1}$ case, exhibiting a small correction to \cite{alexakis2010renormalized} when $n = 2$. Furthermore, we prove the existence of an asymptotic description of our minimal submanifold, $Y$, over the boundary cylinder $\partial Y \times \R^+$, and we further derive an $L^2$-inner-product relationship between $u_2$ and $u_{m+1}$ when $M = \H^{n+1}$. Our results apply to a slightly more general class of manifolds, which are conformally compact with a metric that has an even expansion up to high order near the boundary.
\end{abstract}
\section{Introduction}
We consider the half-space model of $\H^{n + 1} = \{(y, x) \; | \; y \in \R^n, \;\; x \in \R^+\} $ equipped with the complete metric
\[
g = \frac{dy_1^2 + \dots + dy_{n}^2 + dx^2}{x^2}
\]
Renormalized volume arises by trying to make sense of the $m$-dimensional volume of noncompact $Y^m \subseteq \H^{n+1}$ which intersect $\partial \H^{n+1}$ in a compact $(m-1)$-submanifold, $\gamma$, $C^{m+1,\alpha}$-embedded in $\partial \H^{n+1} = \R^n$. The hyperbolic metric is singular along the boundary $\partial \H^{n+1} = \{x = 0\}$, and the $m$-dimensional volume is a priori infinite. But because $\partial Y = \gamma \subseteq \R^n$ is prescribed and $Y$ is minimal, we know the precise manner in which the volume of appropriate cutoffs diverge. The original definition of renormalized volume comes from an asymptotic expansion of the $m$-dimensional volume of $Y \cap \{x > \epsilon\}$ as $\eps \to 0$
\[
\int_{x > \epsilon} dA_{Y} = a_0\epsilon^{-m + 1} + \dots + a_{m-1} \epsilon^{-1} + a_m + O(\epsilon)
\]
and then defining the \textbf{renormalized volume}
\[
\mathcal{V}(Y):= a_m
\]
The process of expanding in $\eps$ is known as \textit{Hadamard regularization}, and it can be used to compute renormalized volume in more general contexts, including Poincar\'e-Einstein (hereon labeled as ``PE") spaces. Though $\mathcal{V}(Y)$ no longer represents the ``volume" of $Y$, it is a Riemannian invariant that reflects the topology and conformal geometry of $Y$ when $m$ is even (cf \cite{alexakis2010renormalized}, Proposition 3.1). When $m$ is odd, the definition depends on the choice of representative of the conformal infinity of $g$, but the ``conformal anomaly" is computable and of physical interest. \nl \nl
%
Our goal is to compute formulae for the first and second variations of renormalized volume for minimal submanifolds of PE spaces. This requires us to prove regularity of minimal submanifolds in PE spaces, which is needed to formally expand the volume form as $\eps \to 0$.  Renormalized volume is typically defined using Hadamard regularization (notable exceptions \cite{o2021uniqueness} \cite{albin2009renormalizing}). We find it more convenient to use Riesz regularization \ref{Riesz}, an equivalent way of defining renormalized volume. Formulae for variations of renormalized volume appear for $Y^2 \subseteq \H^3$ in \cite{alexakis2010renormalized}, and this paper was the primary motivation for our work. We prove results for $Y^m \subseteq M^{n+1}$ of arbitrary dimension and codimension with $M$ PE. When $m$ is odd, renormalized volume depends on the choice of representative in the conformal class of the metric. However, any two such choices lead to definitions of renormalized volume that differ by a boundary integral, depending only on the curvature of $\gamma = \partial Y$, and not the ``global" data of $Y$ in the interior. The first and second variations of renormalized volume for $m$ odd are similarly well defined up to a ``local" boundary integral.
\subsection{Background}
Renormalized volume was originally studied in high energy physics and string theory. We state its physical significance here for historical record: for a $k$-brane in string theory, one can associate a $k$-dimensional submanifold, $Y$, of an ambient manifold, $\overline{X}^{n+1}$. The expected value of the Wilson line operator of the boundary, $W(\partial Y)$, is then given by $\exp(- T \RV(Y))$ where $T$ is the string tension and $\RV(Y)$ is the renormalized volume \cite{graham1999conformal}. Henningson and Skenderis \cite{henningson1998holographic} were the first to compute renormalized volume (in the literature, ``Weyl Anomaly") for low dimension odd examples, and Graham and Witten developed the mathematical theory shortly after. \nl \nl
We are interested in the renormalized volume of minimal submanifolds $Y^m \subseteq M^{n+1}$ where $M$ is a conformally compact, asymptotically hyperbolic, and has an even metric to high order in terms of a ``boundary defining function" $x$. While $M = \H^{n+1}$ is the primary example, we are generally motivated by PE spaces and their deep history. Graham and Lee \cite{graham1991einstein} first discuss existence of PE metrics on $B^{n+1}$ with $\partial M = S^n$. Graham and Witten \cite{graham1999conformal} (and similarly, Graham and Reichert, \cite{graham2017higher}) is the most relevant work for us. They show that renormalized volume is mathematically defined for even dimensional submanifolds of PE spaces, and that a graphical expansion for $Y$ minimal is even in its bdf to high order (\textit{assuming} the expansion exists). One of the main results of this paper is to show existence of such an expansion in aribtrary codimension. There is a long history of showing regularity in codimension one, including Lin \cite{lin1989dirichlet}, Guan, Spruck, Szapiel \cite{guan2009hypersurfaces}, Tonegawa \cite{tonegawa1996existence}, Han, Sehn, Wang \cite{han2016optimal}, and Jiang \cite{jiang2016boundary}. More recently, Mazzeo and Alexakis \cite{alexakis2010renormalized} derive a formula for the first and second variation of renormalized area for $Y^2 \subseteq \H^3$. They also show that these variations record a Dirichlet-to-Neumann type operator, and we generalize the variation formulas. Nguyen and Fine investigate renormalized area through their work on weighted monotonicity theorems with applications to the renormalized area of minimal surfaces in \cite{nguyen2021weighted}. They have further work on minimal $Y^2 \subseteq \H^4$ in preparation.
\subsection{Statement of Results}
\begin{figure}[h!]
\centering
\includegraphics[scale=0.25]{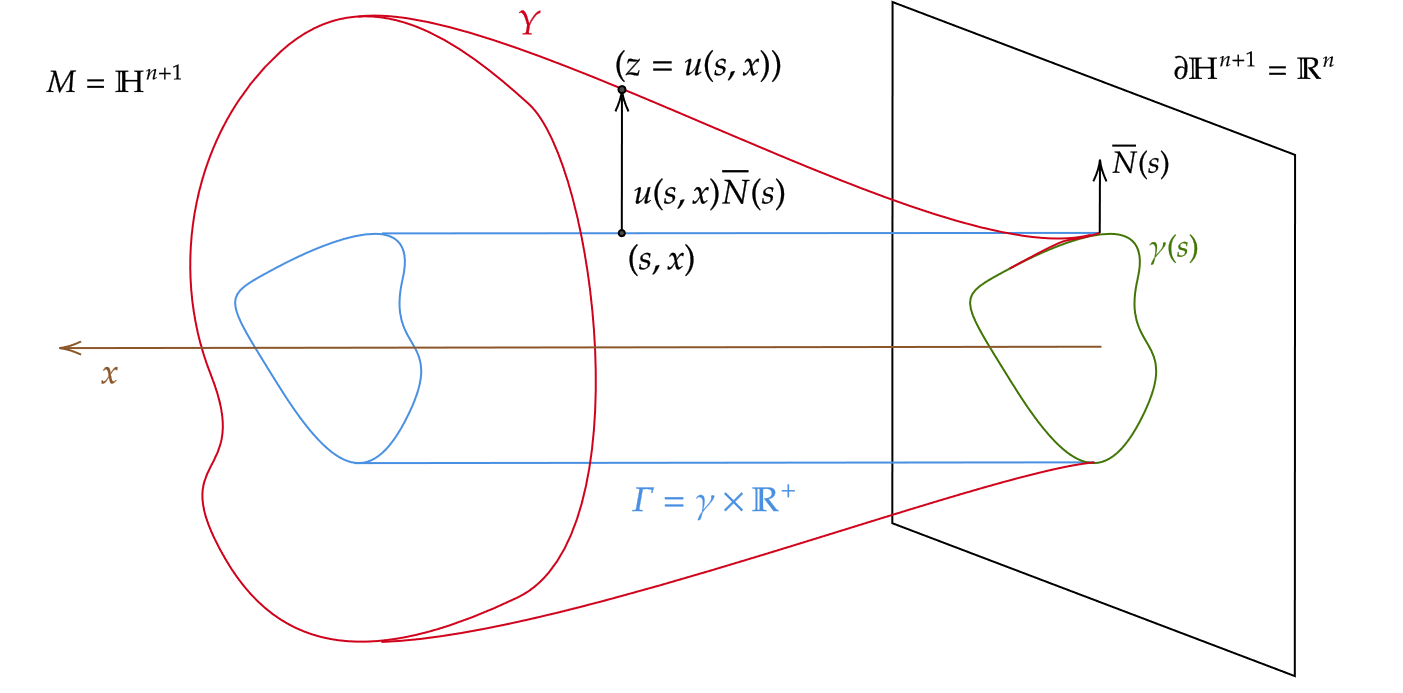}
\caption{}
\label{rvpicturegraphical}
\end{figure}

\begin{figure}[h!]
\centering
\includegraphics[scale=0.25]{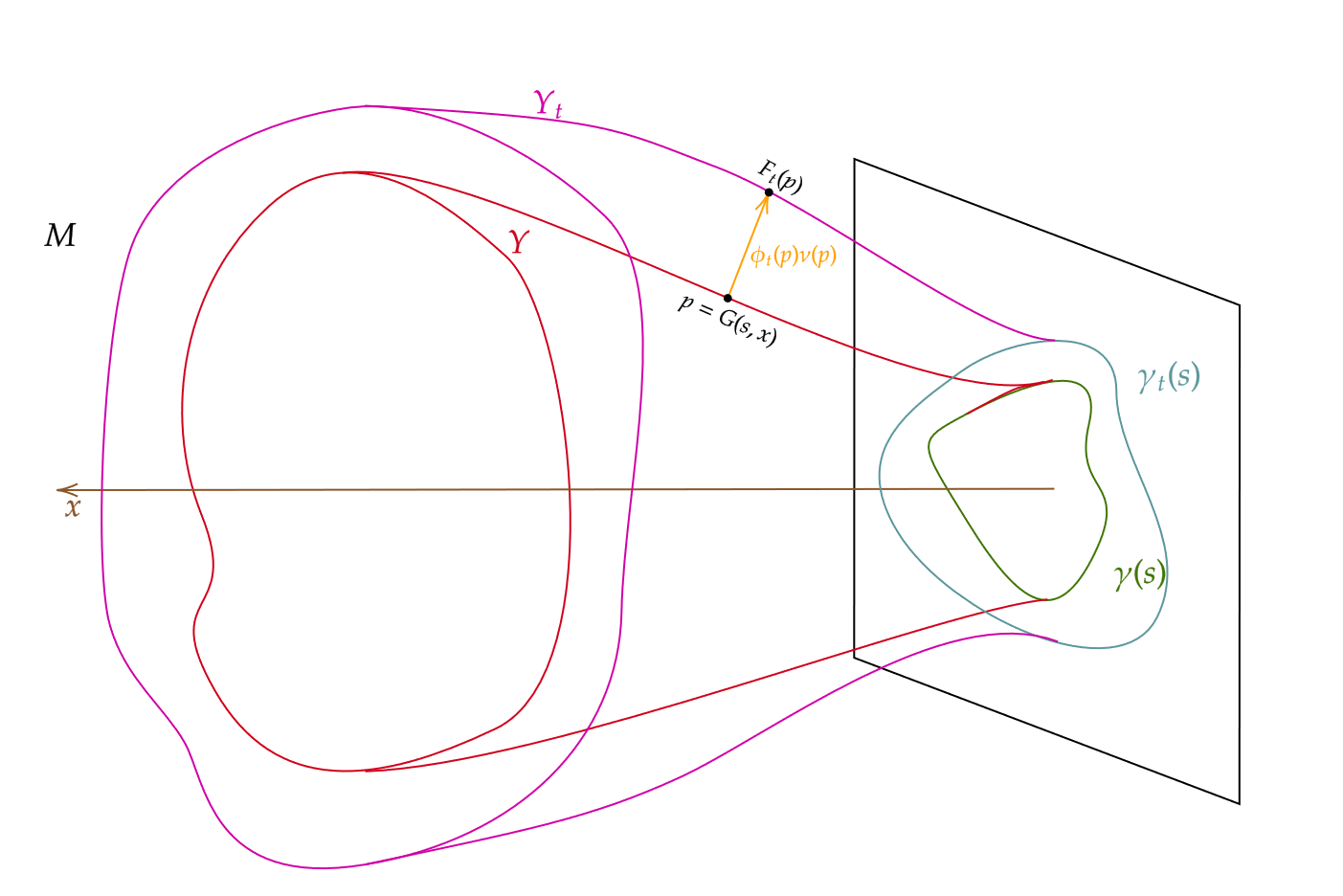}
\caption{}
\label{YtOverY}
\end{figure}

We work with $(M^{n+1}, g)$ PE and $Y^m \subseteq M^{n+1}$ minimal ($m \geq 2$), conformally compact with boundary $\gamma = \partial Y = Y \cap \partial M$. We require that $Y$ be embedded in some neighborhood of its boundary, $\gamma = \partial Y$. WLOG we assume that $\gamma$ is connected and $C^{m+1,\alpha}$ embedded in $\partial M$. Let $x$ be a bdf for $M$ in a neighborhood of $\partial M$ and consider the cylinder over the boundary:
\[
\Gamma = \gamma \times [0, \eps) = \{(x, s) \; | \; s \in \gamma, \qquad 0 \leq x < \eps\}
\]
We assume $Y$ is graphical over $\Gamma$ in a neighborhood of the boundary (see figure \ref{rvpicturegraphical}) and describe $Y$ via the exponential map
\begin{equation} \label{GraphicalParam}
Y \cap \{x \leq \eps \} = \{\overline{\exp}_{\Gamma}(u(s,x))\}
\end{equation}
where $\overline{\exp}$ denotes the exponential map taken with respect to the compactified metric, $\bg = x^2 g$, restricted to elements of $N(\Gamma)$. Here $u = u^i \bN_i \in N \Gamma$ where $\{\bN_i(s,x)\}$ is a normal frame for $\Gamma$ and
%
\noindent $u$ satisfies a degenerate elliptic equation coming from $Y$ being minimal. In \S \ref{Graphical}, we establish regularity of $u$ and prove theorem \ref{AsymptoticExpansion}
\begin{theorem*}
For $Y^m \subseteq M^{n+1}$ minimal and $u = u^i \bN_i$ satisfying \eqref{GraphicalParam}, we have
\begin{empheq}[box=\fbox]{equation*}
u^i(s,x) = \begin{cases}
	u_2^i(s) x^2 + u_4^i(s) x^4 + \dots + u^i_m(s) x^m + u^i_{m+1}(s) x^{m+1} + \dots  & \text{m even}\\[2ex]
	u_2^i(s) x^2 + u_4^i(s) x^4 + \dots + u^i_{m+1}(s) x^{m+1} + U^i(s) x^{m+1} \log(x) + u_{m+2}(s) x^{m+2} + \dots & \text{m odd }
\end{cases}
\end{empheq}
for $C^{m+1, \alpha}$ coefficients $u_{k}(s)$ and $U(s)$.
\end{theorem*}
\noindent \S \ref{Graphical} contains the full details. To make similar statements to the above but more concisely, we recall notation from \cite{albin2009renormalizing}: let 
\begin{align*}
f &: \Gamma \to \R \\
f(s,x) &= f_0(s)  + f_1(s) x + \dots + f_m(s) x^m + O(x^{m+1})
\end{align*}
and define
\begin{align*}
\F&: C^{\infty}(\Gamma) \to \Z \\
\F(f) &= \begin{cases}
0 & \text{if $f$ is } O(x^{m+1}) \\
1 & \text{if $f$ is even below } x^m \text{ and not } O(x^{m+1}) \\
-1 & \text{if $f$ is odd below } x^m \text{ and not } O(x^{m+1}) \\
\text{undefined} & \text{else}
\end{cases}
\end{align*}
When $m$ is odd, we define the above but replacing $m \to m+1$ and allowing for $x^{m+1} \log(x)$ terms. We will often omit the case of $\F = 0$ and write $\F = 1$ or $\F = -1$ for our computations, i.e. any statement of $\F = \pm 1$ should be interpreted as $\F \in \{0, \pm 1\}$ (see \S \ref{ParityDefinition} for a full definition and convention). We note that theorem \ref{AsymptoticExpansion} becomes, $\F(u) = 1$ (or $\F(u) = 0$, implicitly). With this, we informally state theorem \ref{SFFParityMinimal} in codimension $1$
\begin{theorem*}
Suppose that $Y^m \subseteq M^{n+1}$ minimal with $\overline{h} = \overline{g} \Big|_{TY}$ even up to order $x^m$. Let $p \in Y$, $\bA: \text{Sym}^2(TY) \to N(Y)$ denote the second fundamental form, $\bnu$ be a normal to $Y$, both with respect to $\bg$. Then $\bA$ and its covariant derivatives are \textit{even} up to order $x^m$.
\end{theorem*}
\noindent See \S \ref{SFFParity} for the full theorem. We also consider variations of $Y$ among the space of minimal submanifolds. We can describe a smooth family of minimal submanifolds as 
\[
Y_t = \exp_Y (S_t)
\]
for $S_t \in N(Y)$ a smooth function of $t$ and $\exp_Y$ the exponential map with respect to $h = g \Big|_Y$. Let $\dot{S}:= F_*(\partial_t) \Big|_{t = 0}$ and $\ddot{S} = \nabla_{F_*(\partial_t)} F_*(\partial_t) \Big|_{t = 0}$. Both satisfy Jacobi equations when $\{Y_t\}$ is a family of minimal submanifolds, giving regularity and parity. in codimension $1$ we can write 
\begin{align*}
\dot{S} &= \dot{\phi}(s,x) \bnu(s,x) \\
\ddot{S} &= \ddot{\phi}(s,x) \bnu(s,x)
\end{align*}
for $\bnu$ a normal to $Y$ with respect to $\bg$ (see figure \ref{YtOverY}).
\noindent We informally state theorem \ref{Variation}
\begin{theorem*}
For $\{Y_t\}$ a family of minimal submanifolds, and $\dot{S} = \dot{\phi}(s,x) \bnu$, $\ddot{S} = \ddot{\phi} \bnu$:
\begin{empheq}[box=\fbox]{align*}
\forall i, \quad \F(\dot{\phi}^i) &= 1  \\
\forall i, \quad \F(\ddot{\phi}^i) &= 1
\end{empheq}
\end{theorem*}
\noindent i.e. $\dot{\phi}$ and $\ddot{\phi}$ are even in $x$ to high order - see section \S \ref{AsymptoticVariational} for full details. We remark that in order to compute an equation for $\ddot{S}$, we must compute the second variation of mean curvature (i.e. third variation of area). The author was unable to find this result in the literature, so it is stated in proposition \ref{SecondVariationOfMC}. In codimension $1$, we get proposition \ref{codim1SecondVarOfMC}
\begin{proposition*}
For $\{Y_t\}$ a family of minimal submanifolds, and $\dot{S} = \dot{\phi}(s,x) \bnu$, $\ddot{S} = \ddot{\phi} \bnu$, we have that 
\begin{align*}
\frac{d^2}{dt^2} H(t)\Big|_{t = 0} & = [J_Y(\ddot{\phi}) + G(\dot{\phi}, \nabla \dot{\phi}, D^2 \dot{\phi})] \nu = 0
\end{align*}
where $G$ is an error term and quadratic in its arguments. Furthermore $\F(G(\dot{\phi}, \nabla \dot{\phi}, D^2 \dot{\phi})) = 1$.
\end{proposition*}
\noindent We then compute the first and second variations of renormalized volume in theorem \ref{MainThm}. In codimension $1$, $n$ even, we get theorems \ref{codim1FirstVar} \ref{codim1SecondVar}
\begin{theorem*}
\begin{empheq}[box=\fbox]{align*}
\frac{d}{dt} \RV(Y_t)\Big|_{t = 0} &= -(n+1) \int_{\gamma} \dot{\phi}_0(s) u_{n+1}(s) \; dA_{\gamma}(s) \\
\frac{d^2}{dt^2} \RV(Y_t)\Big|_{t = 0} &= \int_{\gamma} \Big( -(n+1) \ddot{\phi}_0 u_{n+1} + (1 - n) \dot{\phi}_0(s) \dot{\phi}_{n+1}(s) \\
& \quad - 8n\dot{\phi}_0(s)^2  u_2 u_{n+1}(s)\Big) \; dA_\gamma(s)  
\end{empheq}
\end{theorem*}
\noindent The full theorem in arbitrary codimension and $m$ odd is stated in \S \ref{MainResult}. We note that while $M$ being PE is the most natural setting, our results hold for a slightly larger class of manifolds - namely those that are conformally compact with a metric, $g$, that splits as in \eqref{splits} with $k(x,s)$ satisfying \eqref{kEquation}. \nl \nl
\noindent \rmk\; Let $\mathcal{D}: \gamma \mapsto u_{n+1}$ be the Dirichlet-to-Neumman map for minimal surfaces with boundary. Then we see that this map is recovered from the first variation of renormalized volume, in the sense that $D \RV'(Y)(\dot{\phi}_0) = \langle \dot{\phi}_0, u_{n+1} \rangle$ and we can let $\dot{\phi}_0$ be arbitrary. This is because $\dot{\phi}_0$ is clearly independent from $u_{n+1}$. \nl \nl
As an application of the second variation formula and regularity of $u$, we prove the following projection relationship, proposition \ref{EvenL2Result}
\begin{proposition*}
For $n$ even, $Y^n \subseteq \H^{n+1}$ minimal with graphical expansion given by $u(s,x)$, we have
\[
\boxed{\langle u_2, u_{n+1} \rangle_{L^2(\gamma)} = 0} 
\]
\end{proposition*}
\noindent \rmk \; The author is unaware of an analogous result for variations of area in the euclidean setting but notes the computations for the first and second variation in \cite{carstea2023inverse} which may lead to a similar orthogonality relation in the asymptotically euclidean case.
\subsection{Outline of Proofs}
This paper has $3$ goals:
\begin{itemize}
\item In \S \ref{Graphical}, we show that $Y$ can be described in Fermi coordinates by a graphical function $u$. We prove that $u$ has an even  expansion as as we approach the boundary, and that $u$ is highly regular (formally ``polyhomogeneous") in this domain
\begin{itemize}
\item The proof relies on Allard's regularity theorem, as well as geometric microlocal techniques from \cite{rafe1991elliptic}, and standard PDE arguments. The author suspects that Allard's theorem can be avoided in establishing the regularity of $Y$, but have yet to find such a proof. 
\item Several authors have contributed to the existence of graphical expansions as in theorem \ref{AsymptoticExpansion} for minimal hypersurfaces in hyperbolic space, including Lin \cite{lin1989dirichlet}, Guan, Spruck, Szapiel \cite{guan2009hypersurfaces}, Tonegawa \cite{tonegawa1996existence}, Han, Sehn, Wang \cite{han2016optimal}, and Jiang \cite{jiang2016boundary}. These authors primarily use classical PDE techniques, and by contrast, we use methods from geometric microlocal analysis to establish regularity. We emphasize that the expansion in higher codimension and the variational formulae are novel.

\item This immediately shows that $h = g \Big|_Y$ and $\bA_Y$ have corresponding even expansions in $x$ as well
\end{itemize}
\item In \S \ref{AsymptoticVariational}, we consider a family of submanifolds close to $Y$, $\{Y_t\}$ with $Y_{t = 0} = Y$. Each $Y_t$ can be written as $Y_t = \exp_Y(S_t(p))$ for some $S_t \in N(Y)$. We show that $\dot{S}$ and $\ddot{S}$ are regular and admit even asymptotic expansions, by computing the first and second variations of mean curvature.
\item In \S \ref{MainResult} and \S \ref{ProofMainResult}, we prove a formula for the first and second variations of renormalized volume for families of minimal submanifolds $\{Y_t\} \subseteq M^{n+1}$. Such formulae appear for minimal surfaces in $\H^3$ in \cite{alexakis2010renormalized}, and we extend their results to $Y^m \subseteq M^{n+1}$ for $m$ and $n$ arbitrary, and $M$ a PE manifold. Past research (\cite{alexakis2010renormalized} \cite{albin2009renormalizing} \cite{graham1999conformal}) focuses on $m$ even, however we extend our results to $m$ odd as well.

\item In \S \ref{codim1}, we specialize our result to the codimension $1$ case, i.e. $m = n$, yielding a slight correction to the second variation formula in \cite{alexakis2010renormalized}. We use this to prove an $L^2$-orthogonality result in \S \ref{Killing}
\end{itemize}
\subsection{Acknowledgements}
The author wishes to thank Rafe Mazzeo for providing the inspiration for this problem, as well as his time spent across many meetings. The author also wishes to thank Otis Chodosh for suggesting the application in \S \ref{Killing}, as well as Brian White, Joel Spruck, and Robin Graham for their background work in this area and thoughtful insight. The author thanks the referee for their in-depth feedback. The author is supported by an NSF Graduate Research Fellowship and a Stanford Graduate Fellowship in Science and Engineering.
\section{Preliminaries}
\label{prelims}

\subsection{Defining Renormalized Volume}
Consider $M^m$ a Poincare-Einstein manifold. For $x: \overline{M} \to \R^{\geq 0}$ a special bdf, the metric splits in Graham-Lee Normal form as
\begin{equation} \label{splits}
g = \frac{dx^2 + k(s,x)}{x^2}
\end{equation}
with $(s,x)$ smooth coordinates on a neighborhood, $U$, of $\partial M$, with $U \cong \partial M \times [0, b)$ for $b > 0$. Here, $k(s,x)$ is a smooth tensor on $T \partial M$, i.e. $k(\partial_x, \cdot) \equiv 0$, and it has an even expansion in $x$ up to order $x^{m-2}$ ($x^{m-1}$) when $m$ is even (odd), i.e. 
\begin{align} \label{kEquation}
m \text{ even} \implies k(s,x) &= k_0(s) + x^2 k_2(s) + \dots + k_{m-2}(s) x^{m-2} + k_{m-1}(s) x^{m-1} + k_m(s) x^m + O(x^{m+1}) \\ \nonumber
m \text{ odd} \implies k(s,z,x) &= k_0 + x^2 k_2 + \dots + x^{m-1} k_{m-1} + x^{m-1} \log(x) K + x^{m+1} k_{m+1} + O(x^{m+2})
\end{align}
\noindent Throughout the paper, we will assume that the boundary metric, $k_0(s)$, is smooth.
\begin{Remark} \label{NoKn+1Term}
Such expansions can be found most readily in Fefferman--Graham \cite{fefferman2011ambient}, Chapters 3 and 4. In particular, we highlight that there is no $k_m x^m$ term for $m$ \textit{odd}, and this is seen through the same techniques used to derive the expansions in \cite{fefferman2011ambient} theorem 3.10.
\end{Remark}
\noindent Using the above expansions, we compute
\[
d\text{Vol} = \sqrt{\det g} dx \wedge ds = \frac{1}{x^m} \sqrt{\det k} \; dx \wedge ds
\]
In \cite{graham1999volume}, Graham showed that for $m$ even,
\begin{align*}
q(x,s) := \sqrt{\det k} & = q_0(s) + q_2(s) x^2 + \dots + q_m(s) x^m + q_{m+1}(s) x^{m+1} + \dots
\end{align*}
i.e. $\text{Tr}(k_{m-1}) = 0$ and $q$ is even up to order $m$. For $m$ even, define
\begin{align*}
\RV(M) &:= \FPeps \int_{x > \eps} d \text{Vol} \\
\int_{x > \epsilon} d \text{Vol} &= \left( \int_{x > b} + \int_{\epsilon}^b \right) \frac{q(x,s)}{x^m} dx ds \\
& = \int_{x > b} d\text{Vol} + \int_{x = \epsilon}^b [x^{-m} q_0(s) + x^{-m + 2} q_2(s) + \dots + q_m(s) + R(s,x)] ds dx \\
&= I(b) + c_0(s) \frac{b^{-m + 1} - \epsilon^{-m + 1}}{1-m} + \dots + c_m(s) (b - \epsilon) + F(b,\epsilon)
\end{align*}
where $R(s,x) = O(x)$ and
\begin{align*}
F(b,\epsilon) &:= \int_{x = \epsilon}^{b} R(s,x) ds dx \\
c_{2k}(s) &:= \int_{\partial M} q_{2k}(s) ds 
\end{align*}
Renormalized volume is then computable as
\begin{align*}
\RV(M) &= \FPeps I(b) + \sum_{k = 0}^{m/2} c_{2k}(s) \frac{b^{-m+1} - \epsilon^{-m + 1}}{1 - m} + F(b, \epsilon) \\
&= I(b) + F(b,0) + \sum_{k = 0}^{m/2} c_{2k(s)} \frac{b^{-m+1}}{1-m}
\end{align*}
A priori, using $x$ seems arbitrary, as there could be several functions like $x$ for which we have an asymptotic expansion in $\epsilon$. Formally, we require $x$ to be a ``special bdf" which we define in  \S \ref{PoincareEinstein}. One can show that renormalized volume is a geometrically natural quantity to consider as it is:
\begin{itemize}
\item Independent of the parameter $b$
\item Independent of the choice of special bdf, $x$, or equivalently independent of the representative in the conformal infinity, $k_0 = \bg \Big|_{\gamma}$
\end{itemize}
The former fact follows by keeping track of boundary terms when integrating applying the FTC. The latter is discussed in \cite{graham1999conformal} among other sources, and are also shown in \S \ref{ConformalInvariance}. These properties only hold for $m$ even. Renormalized volume is defined similarly for odd dimensional submanifolds and is done in \S \ref{equivOfDef}. However, the renormalized volume depends on the choice of $x$, and hence depends on the choice of representative of the conformal infinity. \nl \nl
Note that to have an expansion for $k(s,x)$ (and hence $q(s,x)$) in the first place, there needs to be \textit{some regularity} of the metric as we approach the boundary. When we handle the case of $Y^m \subseteq M^{n+1}$ with the metric induced by restriction, this amounts to \textit{regularity of $Y$ itself}. Thus, we prove regularity of $Y$, which implicitly proves that renormalized volume is mathematically defined for our class of $Y^m \subseteq M^{n+1}$ minimal. 
\subsection{Brief Review of Poincar\'e-Einstein Manifolds}
\label{PoincareEinstein}
The splitting of the metric in \eqref{splits} is motivated by Graham-Lee Normal Form \cite{graham1991einstein} \cite{fefferman2011ambient} for Poincar\'e-Einstein (PE) manifolds. A Riemannian manifold $(M,g)$ is Einstein if $g$ satisfies the Einstein equations. The manifold is Poincare if $g$ is \textit{conformally compact}, i.e. the boundary is compact and there exists a function 
\[
\rho: \overline{M} \to \R^{\geq 0} \qquad \st \qquad \{\rho = 0\} = \partial M, \qquad \overline{\nabla} \rho|_{\partial M} \neq 0
\]
and $\overline{g} = \rho^2 g$ is a nondegenerate metric on $\overline{M}$. Here, $\overline{\nabla}:= \nabla^{\overline{g}}$ and we call $\overline{g} := \rho^2 h$ the \textbf{compactified metric}. Furthermore, $\rho$ is a \textbf{boundary defining function} (bdf). We are interested in $\overline{g} |_{\partial M}$ and how it determines $\overline{g}$ on the interior. Note that if $\varphi: \overline{M} \to \R^+$ is a positive smooth function, then $\tilde{\rho} = \varphi \rho$ is also a boundary defining function. As a result, we can consider the conformal class $[\overline{g} |_{\partial M}]$, which we call the \textit{conformal infinity}. For PE manifolds with a chosen representative, $k_0$, in the conformal infinity, there exists a bdf $x$, for which $\overline{g}$ splits as in equation \ref{splits}. Moreover, $k(s,x)$ is regular up to order $x^m$ as shown in \cite{graham1991einstein}. The bdf $x$ is special if 
\[
||d \log(x)||_{g} = 1
\]
holds in a neighborhood of $\partial M$. Furthremore, by equation \eqref{splits} we have $k_0 = x^2 g \Big|_{\partial M}$. Given these conditions, $x$ is unique (see \cite{fefferman1985Conformal} for details). Renormalized volume is conformally invariant for $m$ even in the sense that it does not depend on the choice of $k_0 \in [\bg|_{\partial M}]$ and the corresponding special bdf used. Thus, we can define renormalized volume for $m$ even as long as we use a special bdf (see \cite{albin2009renormalizing} \cite{graham1999conformal}). 
\begin{example}
Consider the Poincare Ball model of hyperbolic space $M = \H^3$. The metric on $\H^3$ is 
\[
g = \frac{4}{(1 - r^2)^2} \left[ dr^2 + r^2 d\phi^2 + r^2 \sin^2\phi d \theta^2 \right]
\]
is Einstein. We want to find a special bdf, $\rho$, for $\H^3$. We assume that it is rotationally symmetric, i.e. $\rho_{\theta} = \rho_{\phi} = 0$. With this, we compute 
\[
1 = ||d \log(\rho)||_{g}^2 = \frac{\rho_r^2}{\rho^2} g^{rr} = \partial_r(\log(\rho))^2 \frac{(1 - r^2)^2}{4}
\]
we take the negative root and get 
\[
\partial_r (\log(\rho)) = \frac{-2}{1 - r^2}
\]
Integrating and exponentiating, we compute
\[
\rho = A \frac{1 - r}{1 + r}
\]
for $0 \leq r \leq 1$ and some constant $A$. Note that as long as $A \neq 0$, we have $\rho^{-1}(0) = \{ r = 1\} = S^2$, which is the boundary of $\overline{\H^3}$. Suppose that we want to prescribe the standard metric on this boundary. i.e. $k_0(\theta) = \sin^2 \phi d \theta^2 + d \phi^2$. Then we have that 
\[
k_0 = \rho^2 h\Big|_{r = 1} = \frac{4 A^2}{(1 + r)^4} [dr^2 + r^2 d \phi^2 + r^2 \sin^2 \phi d \theta^2] \Big|_{r = 1} = \frac{4 A^2}{16} [d \phi^2 + \sin^2 \phi d \theta^2]
\]
so we choose $A = 2$ so that $\rho$ is positive. We see that $\bg \Big|_{r = 1} = \bg \Big|_{\partial M} = k_0$. Note that
\begin{align*}
\overline{\nabla} \rho & = \overline{g}^{ij} (\partial_i \rho) \partial_j = \overline{g}^{rr} (\partial_r \rho) \partial_r = \frac{(1 + r)^4}{16} \cdot \frac{-4}{(1 + r)^2} \partial_r\\
\overline{\nabla} \rho \Big|_{r = 1} &=  - \partial_r
\end{align*}
which is non-zero.
\end{example}
\noindent We can also compute the renormalized volume of $Y = \H^2 \subseteq \H^3 = M$ in this model.
\begin{example}
Consider the Poincare Ball model of hyperbolic space with $\H^2 \subseteq \H^3$ represented as the geodesic disk (see figure \ref{fig:rvpictureballexample}). The restricted metric on $\H^2$ corresponds to when $\phi = \pi/2$ 
\[
h:= g \Big|_{\H^2} =  \frac{4}{(1 - r^2)^2} [ dr^2 + r^2 d \theta^2]
\]
Because $\rho = \frac{2(1 - r)}{1 + r}$ is rotationally symmetric, it is the special bdf for $\H^2$ by the same computation.
\begin{figure}[h!]
\centering
\includegraphics[scale=0.25]{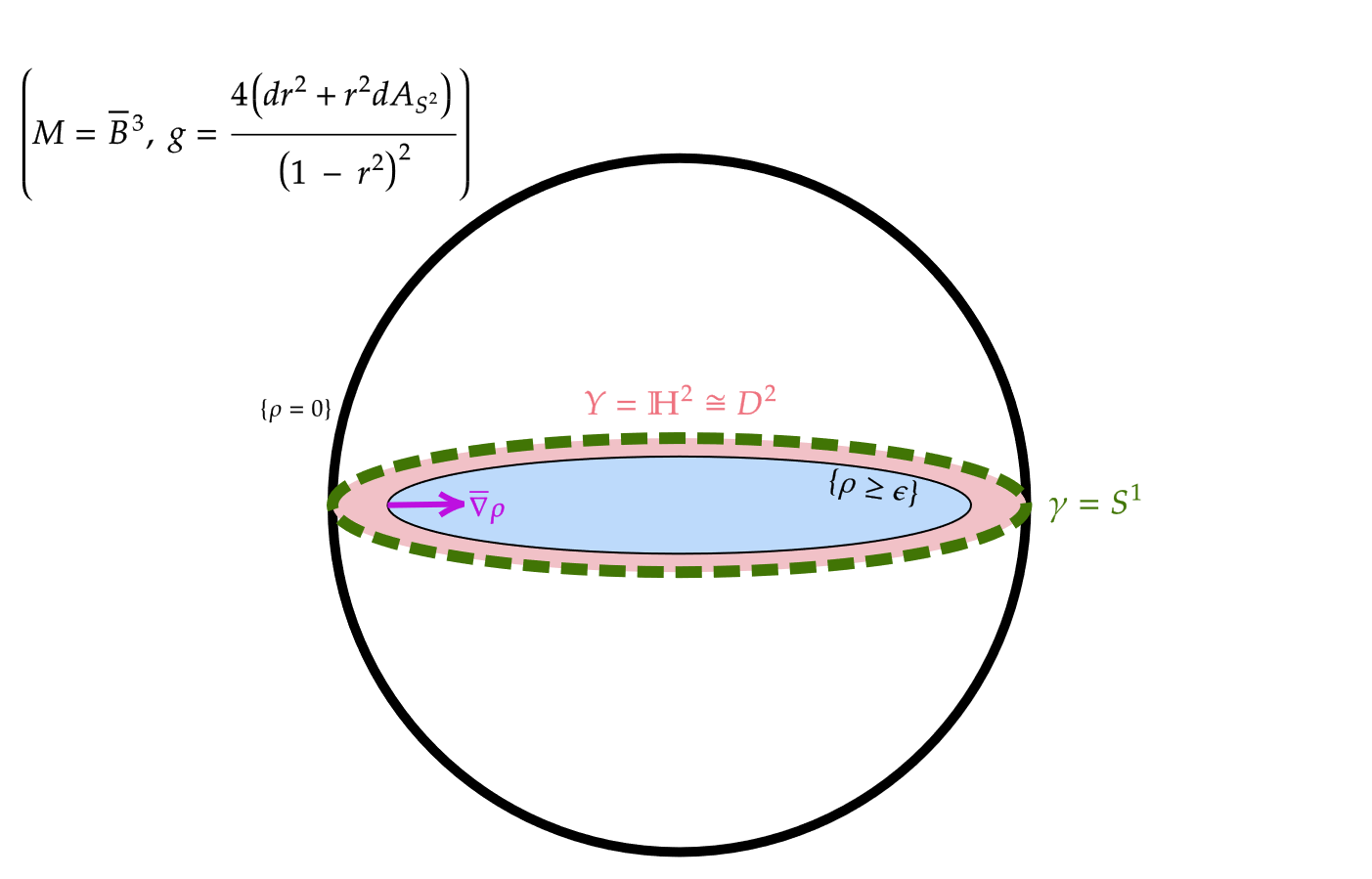}
\caption{Poincare Ball model framed as PE manifold, with $\H^2$ submanifold}
\label{fig:rvpictureballexample}
\end{figure}
\nl\noindent
With this, we can compute the renormalized area of $\H^2 \subseteq \H^3$
\begin{align*}
\RV(\H^2) & = \FPeps \int_{\rho > \eps} dA = \FPeps \int_{\rho > \eps} \frac{4r}{(1 - r^2)^2} dr d \theta \\
&= \FPeps 4 \pi \int_{r = 0}^{(2 - \eps)/(2 + \eps)} \frac{d}{dr} \frac{1}{1 - r^2} dr
\end{align*}
since $\rho > \eps \leftrightarrow \frac{2 - \eps}{2 + \eps} > r$. Integrating, we get 
\[
\int_{\rho > \eps} dA = 4 \pi \left[ (1 - r^2)^{-1} \right]_{r = 0}^{(2 - \eps)/(2 + \eps)} = 4 \pi\left[\frac{4 + 4 \eps + \eps^2}{8 \eps} - 1\right] = 4 \pi \left[ \frac{1}{2 \eps} - \frac{1}{2} + \frac{\eps}{8} \right]
\]
Taking the constant term in $\eps$ then yields
\[
\FPeps \int_{\rho > \eps} dA = 4 \pi \cdot \frac{-1}{2} = \boxed{-2 \pi}
\]
\end{example}
\noindent This example generalizes to higher dimensions as well.
\subsection{Model Case: Half Space Model of $\H^{n+1}$}
Consider $\H^{n+1}$ now with the half-space structure.  The metric is 
\[
g = \frac{dx^2 + (dy_1^2 + \dots + dy_n^2)}{x^2}
\]
so that $k(s,x)$ is the standard Euclidean metric on the first $n$ coordinates, which is even in $x$ as there is no $x$ dependence. Clearly the metric splits in the desired form, and 
\[
||d \log(x) ||_{g}^2 = g^{xx} \partial_x(\log(x))^2 = x^2 \cdot (1/x)^2 = 1 
\]
Moreover the chosen representative of the conformal infinity is
\[
k_0 = dy_1^2 + \dots + dy_n^2 = x^2 g \Big|_{x = 0}
\]
where we take $\R^n = \partial \H^{n+1}$. The issue is that \textbf{the boundary is not compact}. In order for this to be a \textit{conformally compact} manifold, we need to consider the one point compactification of $\R^n$ as the boundary, i.e. $S^n$, and redefine the metric appropriately. Under this compactification, $x$ is no longer a bdf because of the added point at infinity which would have $x = + \infty$ as opposed to $x = 0$. \nl \nl
\underline{\textit{Conformally compact minimal submanifolds of $\H^{n+1}$}}: Despite the above, our analysis in this paper is motivated by and includes $M = \H^{n+1}$ with the half space model. Though $x$ is not a valid bdf for $\H^{n+1}$ itself, it can be used to define renormalized volume for minimal submanifolds with compact boundary (which \textit{are} conformally compact) that are smoothly embedded in a neighborhood of the boundary. We have $\partial Y = \gamma \subset \subset \R^n = \partial \H^{n+1}$ so $x \Big|_Y(p) = 0 \iff p \in \gamma$. This means that $x \Big|_Y$ does define the boundary. However, even if 
\[
||d \log(x)||_{g} = 1
\]
it is not usually true that $||d \log(x)||_h = 1$ for the induced metric $h = g \Big|_Y$. Moreover, the metric may not split, i.e. $h(\partial_{s_a}, \partial_x) \neq 0$. To get around this, the idea is as follows: $Y$ is quadratic and even to high order as we approach the boundary (see \ref{AsymptoticExpansion}). As a result, if we consider a special bdf for $Y$, call it $x_Y$, we can write it in a neighborhood of the boundary as
\[
x_Y = x e^{\omega(s,x)}
\]
where $\omega(s,0) = 0$ and $\omega(s,x)$ has an even expansion up to order $m+2$ (see \S \ref{SpecialBDFonY}). Consequently $k(s,x)$ still has an even expansion up to order $m$ in equation \eqref{splits}, so it makes sense to define 
\[
\RV(Y) = \FPeps \int_{ \{x_Y > \epsilon\} \cap Y} dA
\]
and it turns out that 
\[
\FPeps \int_{ \{x > \epsilon\} \cap Y} dA = \FPeps \int_{ \{x_Y > \epsilon\} \cap Y} dA
\]
by parity considerations. To formally show this, we first introduce Riesz regularization in \S \ref{Riesz}, we then reprove the fact that Riesz regularization produces the result as Hadamard regularization in \S \ref{equivOfDef}, and finally, we show that the usage of $x$ vs. $x_Y$ is irrelevant in defining renormalized volume for minimal submanifolds in $\H^{n+1}$ in \S \ref{SpecialBDFonY}. It is also worth noting that while we consider $M$ a PE space more generally, our analysis of $Y \subseteq M$ is local near a point $p \in \partial Y = \gamma$, for which we can choose coordinate charts resembling hyperbolic space.
%
%
%
\begin{Example}
Consider the geodesic copy of $\H^2$ as a hemisphere of radius $1$ inside $\H^3 = \{(x,y,z) \; | \; x \geq 0\}$ with the metric $\frac{dx^2 + dy^2 + dz^2}{x^2}$. The boundary is a circle of radius $1$, and we parameterize $\H^2$ as 
\[
f(x, \theta) = (x, \sqrt{1 - x^2} \cos \theta, \sqrt{1 - x^2} \sin \theta)
\]
we compute the induced metric
\[ 
h_{xx} = \frac{1}{x^2(1 - x^2)}, \qquad h_{x \theta} = 0, \qquad h_{\theta \theta} = \frac{1 - x^2}{x^2}
\]
so that 
\[
dA_{\H^2} = \sqrt{\det g} dx d \theta = \frac{1}{x^2} dx d\theta
\]
we now compute
\[
\RV(\H^2) = \FPeps \int_{x > \epsilon} dA_{\H^2}
\]
we integrate
\[
\int_{x > \epsilon} dA_{\H^2} = \int_{\theta = 0}^{2\pi} \int_{x = \eps}^1  \frac{1}{x^2} dx d\theta = 2 \pi \left[-\frac{1}{x}\right]_{\eps}^1 = -2\pi + \frac{2 \pi}{\eps}
\]
and so $\RV(\H^2) = \boxed{- 2\pi}$, which is the same result as if we computed the renormalized volume in the ``proper" setting, i.e. the ball model.
\end{Example}

\subsection{Riesz Regularization}
\label{Riesz}
Having defined special bdfs, we can define renormalized volume in an alternate way with Riesz regularization: given an asymptotically hyperbolic manifold, $M$, and a special bdf, $x_M$, on $M$, consider the following function
\begin{align*}
f&: \{ \text{Re}(z) > m \} \to \C \\
f(z) &= \int_M x_M^z dA_M 
\end{align*}
As with Hadamard regularization, the quantity $x_M$ seems unmotivated, however $x_M$  being a special bdf gives $f(z)$ geometric meaning. $f(z)$ is holomorphic for $\text{Re}(z) > m$, and can be extended to be meromorphic on $\C$ with poles at $ z \in \{- \infty, \dots, -1, 0, 1, \dots, m\}$ (see \cite{melrose1996homology}\cite{paycha2003heat} for a proof). We define
\[
\RV(M):= \FPz \int_M x^{z} dA_M = \FPz f(z)
\]
Computing $\FPz f(z)$ amounts to subtracting off the pole at $z = 0$ (if it exists) and evaluating the remaining difference. This process is known as Riesz regularization, and the equivalence of these two definitions is given in the appendix \S \ref{equivOfDef}.
As mentioned before, one can show that for $Y^m \subseteq M^{n+1}$ with $Y$ conformally compact and $m < n + 1$:
\begin{equation} \label{BDFEquivRV}
\FPz \int_Y x_Y^z dA_Y = \FPz \int_Y x^z dA_Y
\end{equation}
On the left hand side, we are using $x_Y$, a special bdf for $Y$ considered as its own asymptotically hyperbolic manifold. On the right hand side, we use $x$, which is a special bdf on $M^{n+1}$. This equation holds for $m$ even, and it holds up to a boundary error for $m$ odd (see \S \ref{SpecialBDFonY}). The latter is expected, as renormalized volume in odd dimensional manifolds depends on the choice of special bdf (\cite{albin2009renormalizing}).
\begin{Example}
We compute $\RV(\H^2)$ for $\H^2 \subseteq \H^3$ using Riesz regularization in the half space model (we leave it to the reader to compute this for the poincare ball model).
\[
\zeta(z) = \int_{\H^2} x^z dA_{\H^2} = \int_{x = 0}^{1}\int_{\theta = 0}^{2\pi} x^{z - 2} dx d\theta = 2\pi \left[\frac{x^{z-1}}{z - 1}\right]_{x = 0}^{x = 1} = 2 \pi \frac{1}{z - 1}
\]
Again, when we find the meromorphic extension, we first assume $\text{Re}(z) \gg 0$ so that $0^{z-1} = 0$. There is no pole at $z = 0$ in this extension, so 
\[
\RV(\H^2) = \FPz \zeta(z) = \zeta(0) = \boxed{-2 \pi}
\]
\end{Example}

\subsection{Parity of functions} \label{ParityDefinition}
Throughout this paper, we will use $x$ to denote a special bdf on our ambient PE space $(M^{n+1}, g)$ and identify a neighborhood of the boundary, $U(\partial M)$, with $\partial M \times [0, \eps)$. When $M = \H^{n+1}$, $x$ is the distinguished direction in the decomposition of $\H^{n+1} \cong \R^n \times \R^+$. Let $f$ be a function defined on $\Gamma \subseteq M$ in coordinates of $(s,x)$. Further assume that $f$ is polyhomogeneous and can be expanded as 
\begin{equation} \label{FunctionExpansion}
f(s,x) = \begin{cases}
f_0(s) + f_1(s) x + \dots + f_m(s) x^m + O(x^{m+1}) & m \text{ is even } \\
f_0(s) + f_1(s) x + \dots + f_{m+1}(s) x^{m+1} + F(s) x^{m+1} \log(x) + O(x^{m+2}) & m \text{ is odd } 
\end{cases}
\end{equation}
then we define for $m < n$ even
\begin{equation} \label{FParityDefEven}
\F(f) = \begin{cases}
0 & f(s,x) \text{ is } O(x^{m+1}) \\
1 & f(s,x) \text{ is even up to } x^m \text{ and not } O(x^{m+1})\\
-1 & f(s,x) \text{ is odd up to } x^{m} \text{ and not } O(x^{m+1}) \\
\text{undefined} & \text{else}
\end{cases}
\end{equation}
When $m = n$, we assume that
\begin{equation} \label{FunctionExpansionCodim1}
f(s,x) = \begin{cases}
f_0(s) + f_1(s) x + \dots + f_m(s) x^m + F(s) x^m \log(x) +  O(x^{m+1}) & m \text{ is even } \\
f_0(s) + f_1(s) x + \dots + f_{m+1}(s) x^{m+1} + F(s) x^{m+1} \log(x) + O(x^{m+2}) & m \text{ is odd } 
\end{cases}
\end{equation}
For $m = n$ even, we define
\begin{equation} \label{FParityDefEvenCodim1}
\F(f) = \begin{cases}
0 & f(s,x) \text{ is } O(x^{n+1}) \\
1 & f(s,x) \text{ is even up to } x^n \log(x) \text{ and not } O(x^{n+1})\\
-1 & f(s,x) \text{ is odd up to } x^{n} \log(x) \text{ and not } O(x^{n+1}) \\
\text{undefined} & \text{else}
\end{cases}
\end{equation}
Similarly for $m = n$ odd, we define
\begin{equation} \label{FParityDefOdd}
\F(f) = \begin{cases}
0 & f(s,x) \text{ is } O(x^{n+2}) \\
1 & f(s,x) \text{ is even up to order } x^{n+1}\log(x) \text{ and not } O(x^{n+2})\\
-1 & f(s,x) \text{ is odd up to order } x^{n+1}\log(x) \text{ and not } O(x^{n+2}) \\
\text{undefined} & \text{else}
\end{cases}
\end{equation}
We note that $\F$ is multiplicative in the sense that if $f$ and $g$ both satisfy equation \eqref{FunctionExpansion} then
\[
\F(fg) = \F(f) \F(g)
\]
We may explicitly write that a given function is ``even/odd up to" a given order when relevant. We are primarily interested in the case of $\F \neq 0$ for all usages of the $\F$ functional. Thus, throughout this paper, any computation of $\F = 1$ signifies $\F \in \{0,1\}$, and similarly $\F = -1$ signifies $\F \in \{0, -1\}$. We adopt this convention for brevity at the expense of some clarity. If there are asymptotics to show that $\F \neq 0$, we will write these explicitly. \nl \nl
\rmk \; The case of $m = n$ even is special because of \eqref{kEquation} for $M^{n+1}$
\begin{align*}
\text{$n+1$ even} \implies k(s,x) &= k_0 + x^2 k_2 + \cdots + k_{n-1} x^{n-1} + k_{n} x^{n} + k_{n+1} x^{n+1} + O(x^{n+2}) \\
\text{$n+1$ odd} \implies k(s,x) &= k_0 + x^2 k_2 + \cdots + k_{n} x^{n} + Kx^{n} \log(x) + O(x^{n+2}) 
\end{align*}
When $Y^m \subseteq M^{n+1}$ with $m < n$, we expect the presence of $k_{n-1} x^{n-1}$ in the even case and $K x^{n} \log(x)$ in the odd case to not affect our formulation of even expansions up to order $m$. However, when $m = n$ even, we expect $K = K(s)$ to give rise to $x^n \log(x)$ terms. We note that when $K \equiv 0$, this separate definition for $\F$ when $m = n$ even is unecessary. In particular, for $M = \H^{n+1} / \Gamma$ for $\Gamma$ a coconvex compact subgroup, no $x^n \log(x)$ term is present and \eqref{FParityDefEven} applies for all $m < n+1$ even. \nl \nl
%
%
We also define a parity preserving first order linear operator $L$ as 
\begin{align} \label{ParityPreserving}
L &:= d^a(s,x) \partial_{s_a} + d^x(s,x) \partial_x \\ \nonumber
\F(d^a) &= 1 \\ \nonumber
\F(d^x) &= -1 \\
\implies \F(L) &= 1
\end{align}
We similarly define a parity preserving first order quadratic differential functional, $Q$, as
\begin{align*}
Q(f, g) &:= d^{ab}(s,x) f_a g_b + d^{ax} f_a g_x + d^{xx} f_x g_x \\ \nonumber
\F(d^{ab}) &= 1 \\ \nonumber
\F(d^{ax}) & = -1 \\ \nonumber
\F(d^{xx}) &= 1 \\
\implies \F(Q) &= 1
\end{align*}
Higher order parity preserving linear operators and quadratic functionals are defined analogously.
\subsection{Variation of Renormalized Volume} \label{VariationsOfRV}
When computing the variation of the renormalized volume, we consider $\{Y_t^m\} \subseteq M^{n + 1}$, a one-parameter family of minimal submanifolds with $Y_{t = 0} = Y$ our designated submanifold. We require that each $Y_t$ be embedded in some neighborhood of the boundary $\mathcal{U} \cong \partial M \times [0, \epsilon)$. Define
\begin{align*}
S_t&: Y \to N(Y), \qquad &\dot{S} &:= \partial_t S_t \Big|_{t = 0} \\
F_t&: Y \to Y_t, \qquad &F_t(p) &:= \exp_p(S_t(p))
\end{align*}
Given that this is a variation among minimal submanifolds, we know that $\dot{S}$ lies in the kernel of the Jacobi operator of $N(Y)$, i.e.
\[
J_Y^{\perp}(\dot{S}) = \Delta_Y^{\perp} (\dot{S}) + \tilde{A}(\dot{S}) + \text{Tr}[R_M(\cdot, \dot{S}) \cdot] = 0
\]
where $\tilde{A}(X) = g((\nabla_{v_{\alpha}} v_{\beta})^N, X) g^{\alpha \gamma} g^{\beta \delta} (\nabla_{v_{\gamma}} v_{\delta})^N$ is the Simons operator and $\text{Tr}[R_M(\cdot, \dot{S})]$ denotes the trace of the ambient Riemann curvature tensor, $R_M$, taken over $TY$, applied to $\dot{S}$. As a result, $\dot{S}$ satisfies a regularity theorem stated in full in \S \ref{Variation}. In codimension $1$, $\dot{S} = \dot{\phi}(s,x) \bnu(s,x) = [x^{-1} \dot{\phi}(s,x)] \nu$ for $\nu$ a normal to $Y$. Then 
%
\[
\forall i, \quad \F(\dot{\phi}^i) = 1
\]
i.e. $\dot{\phi}$ is even to $x^n$ or $x^{n+1}$ with the presence of a log term when $n$ is odd. Similarly, we show in the appendix that $\ddot{S}$ satisfies an equation of the form
\[
J_Y^{\perp}(\ddot{S}) = Q^{\perp}(\dot{S}, \dot{S})
\]
where $Q^{\perp}$ is a quadratic functional in $\{\dot{S}^i, \dot{S}^i_{\alpha}, \dot{S}^i_{\alpha \beta}\}$ valued in $NY$. This establishes regularity in a very similar manner and proves 
\[
\forall i, \quad \F(\ddot{\phi}^i) = 1
\]
In the codimension $1$ even case, neither $\dot{\phi}$ nor $\ddot{\phi}$ have $x^n\log(x)$ terms and the details are done in section \S \ref{SecondVarDetails}. Having established regularity of $\dot{S}$ and $\ddot{S}$, we define
\[
\RV(Y_t) = \FPz \int_{Y_t} x^z dA_t = \FPz \int_Y F_t^*(x)^z F_t^*(dA_t) 
\]
Computing variations of this amounts to differentiating the integrand and interpreting it geometrically in terms of $u(s,x)$, $\dot{S}$, and $\ddot{S}$. With this, we state formulae for the first and second variations of renormalized volume in codimension $1$. The full theorem is stated in \S \ref{MainResult}, theorem \ref{MainThm}.
\begin{theorem*}
For $\{Y_t^n\} \subseteq M^{n+1}$ a one-parameter family of hypersurfaces satisfying mild geometric constraints, suppose $Y_{t = 0} = Y$ is minimal with $\partial Y = \gamma$. Then
\begin{empheq}[box=\fbox]{align*}
n \; \text{even} \; & \implies \frac{d}{dt} \RV(Y_t)\Big|_{t = 0} = -(n+1) \int_{\gamma} \dot{\phi}_0(s) u_{n+1}(s) \; dA_{\gamma}(s) \\
n \; \text{odd} \; & \implies \frac{d}{dt} \RV(Y_t)\Big|_{t = 0} = -(n+1) \int_{\gamma} \left[ \dot{\phi}_0(s) u_{n+1}(s) + F(\dot{\phi}_0, u_2)(s)\right] dA_{\gamma}(s)
\end{empheq}
where $F$ is a polynomial in $\dot{\phi}_0$, $u_2$, and their higher derivatives. If in addition each $Y_t$ is minimal, then we have 
\begin{empheq}[box=\fbox]{align*}
n \; \text{even} \;  \implies \frac{d^2}{dt^2} \RV(Y_t)\Big|_{t = 0} &= \int_{\gamma} -(n+1) \ddot{\phi}_0 u_{n+1} + (1 - n) \dot{\phi}_0(s) \dot{\phi}_{n+1}(s) \\
& \quad -8n\dot{\phi}_0(s)^2 u_2 u_{n+1}(s) \; dA_\gamma(s)  \\
n \; \text{odd} \;  \implies \frac{d^2}{dt^2} \RV(Y_t)\Big|_{t = 0} &= \int_{\gamma} -(n+1) \ddot{\phi}_0 u_{n+1} + (1 - n) \dot{\phi}_0(s) \dot{\phi}_{n+1}(s) \\
& \quad + \dot{\phi}_0(s)^2 \left[-8n u_2 u_{n+1}(s) + \text{Tr}_{T\gamma}(k_{n+1,0}) \right] \\
& - \dot{\phi}_0(s) \left[4(n+2) \dot{\phi}_0(s) u_2(s) U(s) + \dot{\Phi}(s) \right] + F_2(\ddot{\phi}_0, \dot{\phi}_0, u_2) \; dA_\gamma(s) 
\end{empheq}
where $k_{n+1,0}(s) = k_{n+1}(s,0)$ in \eqref{kEquation}.
\end{theorem*}

\section{Graphical Asymptotic Expansion}
\label{Graphical}

\subsection{Results about $u$}
In this section, we leverage the fact that $Y$ is minimal and smoothly embedded in a neighborhood of the boundary to get a polyhomogeneous expansion of each $u^i$ for $u(s,x) = (u^1(s,x), \dots, u^{n-m + 1}(s,x)) \in N(\Gamma)$. Recall that $u$ is polyhomogeneous if 
\[
u(s,x) \sim \sum_{\text{Re}(z_j) \to \infty} \sum_{t = 0}^{N_j}  x^{z_j} \log(x)^t a_{j,t}(s) \; \; \st \; \; a_{j,t} \in C^{\infty}(\gamma)
\] 
To show polyhomogeneity we establish some initial regularity. We assume that as $x \to 0$, the \textit{blown up localized mass} of $Y$ approaches $1$. Formally, let $x_0 \ll 1$, $s_0 \in \gamma$, and define
\begin{equation} \label{BlowUpMap}
F_{s_0, x_0} := (s,x,z) \to (\sigma, \xi, \eta) := \left(\frac{s - s_0}{x_0}, \frac{x}{x_0}, \frac{z - u(s_0, x_0)}{x_0}\right)
\end{equation}
When $M = \H^{n+1}$, $F_{s_0, x_0}$ is an isometry. 
\begin{Assumption}
For $Y^m \subseteq M^{n+1} $ minimal, let $Y_0 = F_{s_0, x_0}(Y)$. Assume
\begin{align} \nonumber
\forall \delta > 0, \quad \exists x^* > 0, \quad \st \quad & \forall x_0 < x^* \\ \label{MassAssumption}
\omega_m^{-1} \rho^{-m} ||V_{Y_0, x_0^{-2} F_{s_0, x_0}^*(\bg)}||(B_{\rho}(a)) &\leq 1 + \delta
\end{align}
for all $s_0 \in \gamma$ and for all $a \in B_{x_0}(s_0, x_0)$. Here, $\omega_m$ is the volume of the $m$-dimensional Euclidean ball of radius $1$, and $||\overline{V}_{Y_0, x_0^{-2} F_{s_0, x_0}^*(\bg)}||(B_{\rho}(a))$ denotes the mass of the varifold intersected with a small ball with respect to the metric $x_0^{-2} F_{s_0, x_0}^*(\bg)$. 
\end{Assumption}
%
\noindent This geometric constraint requires that our minimal surfaces ``flatten" out as we blow up near the boundary. This restriction is stronger than what is needed to apply Allard regularity, but it gives the correct $C^{1,\alpha}$ norm bounds. The author hopes that this can be proven with a weaker assumption. With this, we state our regularity theorem.
\begin{theorem} \label{AsymptoticExpansion}
Suppose $Y^m \subseteq M^{n+1}$ minimal satifying equation \eqref{MassAssumption} and $\gamma = \partial Y = Y \cap \partial M^{n+1}$ is a $C^{m+1,\alpha}$ embedded submanifold in $\partial M^{n+1}$. Further suppose that $Y$ is embedded and graphical in some neighborhood of the boundary $\mathcal{U} \cong \partial M \times [0, \epsilon)$. Let $u(s,x) = u^i(s,x) \partial_{z_i} \in N(\Gamma)$, which describes $Y$ as in \S \ref{FermiCoords}. Then, each $u^i(s,x)$ is polyhomogeneous and even to order $m$ ($m+1$) for $m$ even (odd). 
\begin{empheq}[box=\fbox]{equation*}
u^i(s,x) = \begin{cases}
	u_2^i(s) x^2 + u_4^i(s) x^4 + \dots + u^i_m(s) x^m + u^i_{m+1}(s) x^{m+1} + O(x^{m+2} \log(x)) \text{m even}\\[2ex]
	u_2^i(s) x^2 + u_4^i(s) x^4 + \dots + u^i_{m+1}(s) x^{m+1} + U^i(s) x^{m+1} \log(x) + O(x^{m+2} \log(x)) & \text{m odd }
\end{cases}
\end{empheq}
Here, $\{\partial_{z_i}\}$ is a coordinate basis for $N(\Gamma)$, and $u^i_{k}(s)$, $U^i(s)$ are $C^{m+1, \alpha}$ functions on $\gamma$.
\end{theorem} 
\noindent \rmk \; This theorem justifies the existence of an asymptotic expansion for $u$, the graphical function of $Y$, as in \cite{graham1999conformal}. \nl \nl
There are several steps to the proof, which we carry out in the following sections: 
\begin{enumerate}
\item In \S \ref{Maximum}, we use the maximum principle and the fact that $Y$ is minimal to show that $u$ is $O(x^2)$.
\item In \S \ref{ZeroSpace}, we use Allard's regularity theorem and assumption \eqref{MassAssumption} to establish $u \in C_0^{1,\alpha}$. We then use the theory of edge operators as in \cite{rafe1991elliptic} to prove that $u$ is infinitely regular with respect to edge operators.
\item In \S \ref{Revamped} we note that $Y$ is minimal so $u$ also satisfies a degenerate elliptic PDE. We reframe the PDE in terms of the $0$-operators, $(x \partial_x)$ and $(x \partial_{s_a})$.
\item In \S \ref{Parametrix}, we upgrade regularity in $0$-operators to regularity in $b$-operators, $\{x \partial_x, \partial_{s_a}\}$.
\item In \S \ref{IterationArgument}, we upgrade regularity in $b$-operators to $u$ having a polyhomogeneous expansion using a power series iteration in $x$. This follows by linearizing the minimal surface system about successive iterations of $u$, i.e. $u = 0$, $u = u_2 x^2$, $u = u_2x^2 + u_4 x^4 + \dots$. 
\end{enumerate}
As remarked in the previous section, the regularity of $Y$ allows us to formally define renormalized volume
\begin{corollary} \label{RVWellDefined}
For $Y^m \subseteq M^{n+1}$ as above with $m$ even, the renormalized volume
\[
\RV(Y) := \FPeps \int_{x > \eps} dA_Y = \FPz \int_Y x^z dA_Y
\]
is formally defined and independent of the special bdf $x$. For $m$ odd, $\RV(Y)$ is defined as above, but it depends on the choice of $x$.
\end{corollary}

\subsection{Coordinates and Notation}
\label{FermiCoords}
We coordinatize our space as follows: let $p \in \gamma$ be labeled by geodesic normal coordinates on $\gamma$ about some base point $p_0$, i.e.
\begin{equation} \label{gammaMap}
p = f(s) := \overline{\exp}_{p_0}^{\gamma} (s^a E_a)
\end{equation}
where $\{E_a\}$ is an ONB at $p_0$ spanning $T_{p_0} \gamma$. We then map to the cylinder
\[
R(s,x) = (f(s), x) \in \partial M \times [0, \eps)
\]
where we implicitly use the diffeomorphism of $\partial M \times [0, \eps) \cong U \subseteq \overline{M}$ for an open neighborhood of the boundary. We define 
\[
F(s, x, z) = \overline{\exp}^{\Gamma}_{R(s,x)} (z^i X_i)
\]
where $z = (z_1, \dots, z_{n-m+1})$ are coordinates for the normal bundle, $N \Gamma$, and $\{X_i\}$ is an ONB at $p = R(s,x)$. Note that in both instances, $\overline{\exp}$ denotes the exponential map with respect to the compactified metric, $\bg$, restricted to $\gamma$ and $\Gamma$, respectively. We coordinatize $Y$, in some neighborhood of the cylinder $\Gamma$, via 
\[
Y \in q = F(s,x,u(s,x)) \leftrightarrow (s, x, z = u(s,x))
\]
for  $||s||$ close to $0$ and $x < \epsilon$. This is the definition of the function $u(s,x)$ as an $m$-vector in $N(\Gamma)$, and we investigate this function in the next section. \nl \nl
Finally, we will use $v_{(\cdot)}, \; \partial_{(\cdot)}$ to denote a variety of vectors in $T M = T Y \oplus N Y = T \Gamma \oplus N \Gamma$. Here, we notate
\begin{equation} \label{IndexNotation}
\begin{gathered}
a,b,c,d \leftrightarrow s_a, s_b, s_c, s_d \\
i,j,k, \ell \leftrightarrow z_i, z_j, z_k, z_{\ell} \\
i,j,k,\ell \leftrightarrow w_i, w_j, w_k, w_{\ell} \\
\alpha, \beta, \gamma, \delta \leftrightarrow \{y_{\alpha}, y_{\beta}, y_{\gamma}, y_{\delta} \} \subseteq \{s_a, x\} \\
\sigma, \mu, \nu, \tau, \omega \leftrightarrow \{y_{\sigma}, y_{\mu}, y_{\nu}, y_{\tau}, y_{\omega}\} \subseteq \{s_a, x, z_i\} 
\end{gathered}
\end{equation}
We recognize the abuse of notation between the $i,j,k,\ell$. The context will be clear when using these indices to refer to the fermi normal frame off of $\Gamma$, i.e. $\{\partial_{z_i}\}$, vs. the normal frame off of $Y$, $\{w_i\}$, defined in section \S \ref{NormalFrame}

\subsection{Metric on $Y$}
We have the coordinate representation
\begin{equation} \label{GMap}
G(s,x) := (F(s,u(s,x)), x) = (s,z = \vec{u}(s,x),x) \leftrightarrow p \in Y
\end{equation}
for $G: \Gamma \to Y$. We define
\begin{align*}
v_{a} &= G_*(\partial_{s_a}) \\
\vx &= G_*(\partial_x) \\
\va, \vb, \vg, \vd & \in \{v_a, v_x\} \\
\vs, v_{\mu}, v_{\nu}, \vt, \vo & \in \{v_a, \vx, w_i\}
\end{align*}
where $\{w_i\}$ is the aforementioned normal frame. We also define
\begin{equation} \label{SigmaNotation}
\sigma(\omega) := \begin{cases}
0 & \omega \neq x \\
1 & \omega = x
\end{cases}
\end{equation}
be an operator on indices. \nl \nl
We now define $h$, the induced metric on $Y$ by nature of being embedded in $M^{n+1}$, as well as $\overline{h} = x^2 h$. Assuming $u^i = O(x^2)$ and $\F(u^i) = 1$ (verified in the next section \S \ref{Maximum}), we have from \S \ref{MetricTY}
\begin{align*} 
\bh_{ab} & := \overline{g}(G_*(\partial_{s_a}), G_*(\partial_{s_b})) \\
&= \delta_{ab} + O(x^2) \\
\F(\bh_{ab}) &= 1 \\
\bh_{ax} &:= \overline{g}(G_*(\partial_{s_a}), G_*(\partial_{x}) ) \\
&= O(x^3) \\
\F(\bh_{ax}) &= -1 \\
\bh_{xx} &:= g(G_*(\partial_{x}), G_*(\partial_{x}) )\\
& = 1 + O(x^2) \\
\F(\bh_{xx}) &= 1
\end{align*}
Note that $h = g \Big|_Y$ is the \textit{complete metric} for $Y$, while $\bh = \bg\Big|_Y$, is the \textit{compactified metric} (we use $x$, not $x_Y$ here). Moreover, $\{v_{\alpha}\} = \{v_a, v_x\}$ is a basis for $TY$, with $\alpha$ taking on any of the $x$ and $a$ subscripts.
\begin{example}
The compactified metric on $\H^{n+1}$ is $x^2 \frac{dx^2 + dy_1^2 + \dots + dy_n^2}{x^2} = dx^2 + dy_1^2 + \dots + dy_n^2$ which is just the standard Euclidean metric. 
\end{example}

\subsection{Maximum principle argument}
\label{Maximum}
For $\gamma \embed \partial M$ a compact $C^{m+1,\alpha}$ embedded submanifold, consider an $\epsilon$-tubular neighborhood $N_{\epsilon}(\gamma) \subseteq \partial M^n$. WLOG assume that $\gamma$ is connected, and localize about some $p \in \gamma$. The goal is to show that for $u(s,x) = u^i(s,x) \bN_i(s,x)$ and $\forall x < x_0$ sufficiently small, we have $C = C(x_0)$ such that
\[
|u^i(s,x)| \leq C x^2
\]
\subsubsection{Model Case: $M = \H^{n+1}$}
\begin{figure}[h!]
\centering
\includegraphics[scale=0.35]{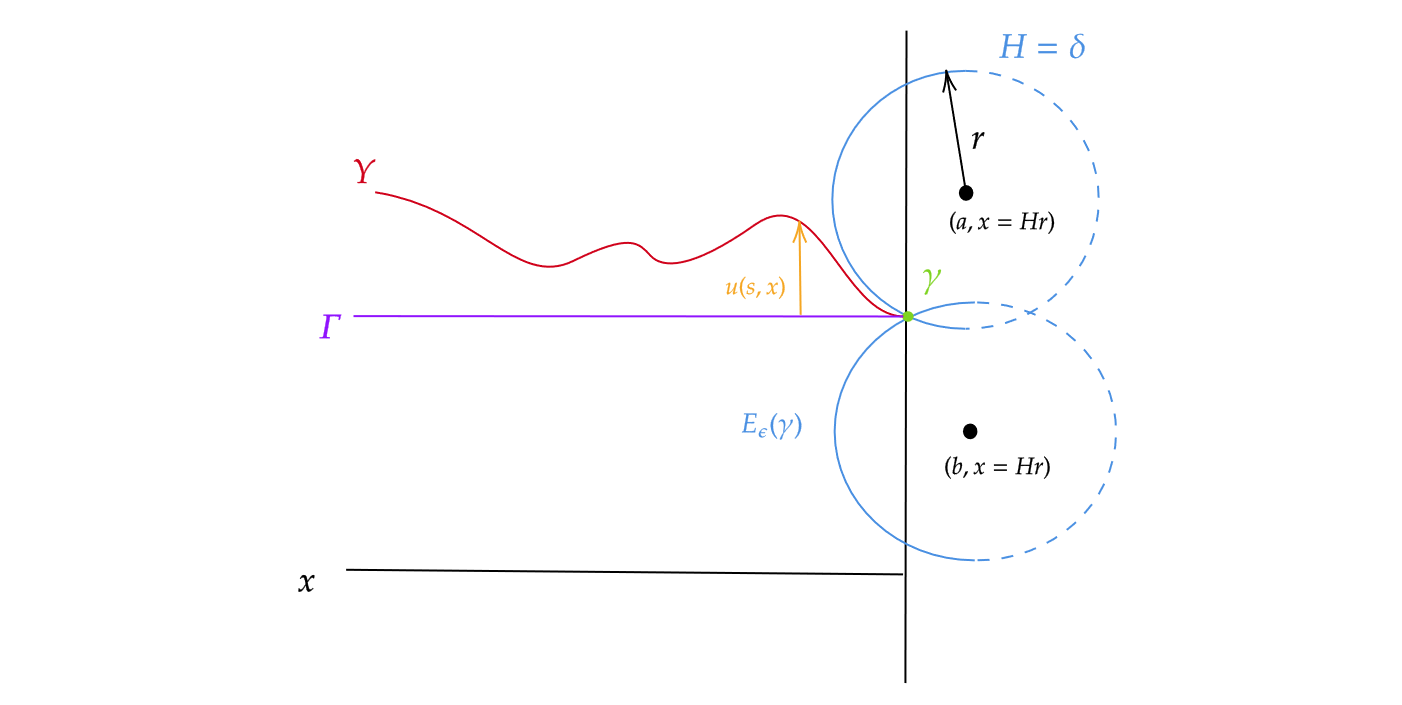}
\caption{Picture of barrier and envelope argument}
\label{fig:modifiedbarrierargument}
\end{figure}

In this case, one can form an envelope of geodesic copies of $\H^n$ as hemispheres to act as a boundary. This argument is historic, originally due to Anderson. We choose to present another argument inspired by \cite{guan2000hypersurfaces}. Let $HS^{n}(R)\subseteq \H^{n+1}$ be the half-sphere of radius $R$ which is a geodesic copy of $\H^n \subseteq \H^{n+1}$. Imagining $\H^{n+1} \subseteq \R^{n+1}$ and hence $HS^n \subseteq \R^{n+1}$, we can shift the center to the right by $x$ to make a new surface, 
\[
HS^n(R, x) = \left(x + HS^n(R) \right) \cap \H^{n+1}
\]
For $x = \delta R$, this is a hypersurface with $H = \delta$ lying inside $\H^{n+1}$. In fact, each of the principle curvatures of this surface is equal to $\frac{\delta}{n}$, so $HS^n(R, \delta R)$ is in fact an $m$-convex surface with 
\[
K_m = \kappa_1 + \dots + \kappa_m = \frac{m}{n} \delta
\]
for any of the $m$ principal curvatures. We use $HS^n(R, \delta R)$ as a barrier around $Y$ (see picture \ref{fig:modifiedbarrierargument}).
Represent $HS^n(R, \delta R)$ graphically over the boundary cylinder $\Gamma_{S^{n-1}} = \partial HS^n(R, \delta R) \times \R^+$ as 
\[
HS^n(R, \delta R) = \exp_{\bg}(v(s,x) N^{n-1}(s))
\]
where $\overline{N}^{n-1}(s)$ is a normal to $\partial HS^{n}(R, \delta R) \subseteq \R^n = \partial \H^{n+1}$.  In these coordinates we have 
\begin{align*}
v(s,x) &= \sqrt{R^2 - (x + \delta R)^2} - R \sqrt{1 - \delta^2} \\
&= -\frac{\delta }{ R^{2}(1 - \delta^2)^{3/2}} x + O(x^2)
\end{align*}
and such a construction holds for any $\delta > 0$. We can repeat this construction about any $p \in \gamma$ such that $\partial HS^n(R, \delta R)$ lies tangent to $\gamma$ for $R = R(p)$ sufficiently small. Now consider the envelope 
\[
E = \partial \left(\bigcup_p HS^n(R(p), \delta R(p), p) \right)
\]
$E$ is now a barrier for $Y$. Let $u(s,x) = u^i(s,x) \bN_i(s,x)$. The maximum principle for $m$-mean convex submanifolds (cf. \cite{jorge2003barrier}, \cite{white2009maximum}) then gives that about any $p$,
\begin{align*}
|u^i(s,x)| & \leq |v(s,x)| \\
&\leq C \delta x
\end{align*}
where $C = C(p)$. In particular at $x = \delta$, we get
\[
|u^i(s,\delta)| \leq C \delta^2
\]
Noting that $\gamma$ is compact and repeating this construction for all $\delta > 0$ sufficiently small, we have
\[
|u^i(s,x)| \leq C x^2
\]
for all $x$ sufficiently small, and some $C$ uniform in $p \in \gamma$.
\subsubsection{$M$ general PE manifold} 
We outline the argument as follows:
\begin{itemize}
\item Find $\eps_0$ sufficiently small so that when we expand
\begin{align} \nonumber
g &= \frac{dx^2 + k(s,x)}{x^2} \\ \nonumber
k(s,x) &= k_0(s) + R \\ \label{tensorBound}
||R|| &\leq C x^2
\end{align}
for all $x < \eps_0$. \eqref{tensorBound} is a tensor bound in $C^{1}_{\bg}$ (see \S \ref{ZeroSpace})

\item Let $\rho$ be the radius such that $N_{\rho}(\gamma) \subseteq \partial M$ is embedded, i.e. the normal bundle is embedded. Let $R = \rho/2$. 
%
Consider $Z := HS^n(R, \delta R)$ for $\delta < \min(\eps_0, \rho/2)$ and note that by \eqref{tensorBound}, we have that each of the principle curvatures satisfy
\begin{align*}
\kappa_{i}(Z) & = \frac{\delta}{n} + E(\eps) \\
\implies K_m(Z) &= \sum_{i = 1}^m \kappa_i(Z) = \frac{m}{n}\delta + E(\eps) \\
E(\eps) & \leq K \eps^2 \\
\implies K_m(Z) \Big|_{x = \delta} &\geq c_0 \delta
\end{align*}
The idea being that because $k(s,x)$ is even up to order $m \geq 2$, $K_m(Z)$ is the same up to quadratic error. Thus, a barrier which is $m$-mean strictly convex with $M = \H^{n+1}$ is still $m$-mean strictly convex for $M$ a general PE manifold.

\item Consider the envelope $E_{\delta, R}(\gamma)$ defined by
\[
E_{\delta, R}(\gamma) := \partial \left( \bigcup_{p \in \gamma} HS^n(R(p), \delta R(p), p) \right)
\]
where $HS^n(R, \delta R, p)$ denotes the above construction based at a point $p \in \gamma$. The same $m$-mean convex maximum principle tells us that $E_{\delta, R}(\gamma)$ is a barrier for $Y$

\item Let $v(s,x)$ be the graphical height function for the envelope over its boundary cylinder. As before,
\begin{align*}
v(s,x) &= -\frac{\delta }{ R^{2}(1 - \delta^2)^{3/2}} x + O(x^2)
\end{align*}
Then we have by the barrier arguments that 
\[
|u^i(s,x)| \leq C \frac{\delta}{R^2} x
\]
Choosing $x = \delta$ (recalling that $R$ independent of $\delta$), we have 
\[
|u^i(s,\delta)| \leq C \delta^2
\]
for $C$ independent of $\delta$ and $p \in \gamma$. Repeat for all $\delta > 0$ sufficiently small to get 
\[
\exists \delta_0 \;\; \st \;\; \forall x < \delta_0, \qquad |u^i(s,x)| \leq C x^2
\]
\end{itemize}

\subsection{Showing $v \in x^2 \bigcap_k C_0^{k, \alpha}$}
\label{ZeroSpace}
In this section, we demonstrate that $u \in x^2 \bigcap_k C_0^{k,\alpha}$ i.e. $u$ is smooth and $u^i = x^2 f^i$ for some $\{f^i\}$ such that 
\[
\forall j, i, \beta \quad \exists \quad C_{j \beta}^i < \infty \quad \st \quad ||(x \partial_x)^j(x \partial_{s_a})^{\beta} f^i||_{C^{0,\alpha}_0} \leq C_{j \alpha}^i
\]
%
for $j$, $\alpha$ arbitrary. Here, $C_0^{k, \alpha}$ is the H\"older space of functions in terms of the edge operators, $\{x \partial_x, x \partial_{s_a}\}$, and
\[
||f||_{C^{k,\alpha}_0} := \sum_{j + |\beta| \leq k} ||(x \partial_x)^j(x \partial_{s_a})^{\beta} f^i||_{0,\alpha, 0}
\]
where $||\cdot ||_{0, \alpha, 0}$ denotes the geometric H\"older norm on $U$ given by 
\[
||f||_{0,\alpha, 0} = \sup_{(s,x) \in U} |f(s,x)| + \sup_{(s,x) \neq (\tilde{s}, \tilde{x}) \in U} \frac{|f(s,x) - f(\tilde{s}, \tilde{x})|(x + \tilde{x})^{\alpha}}{(|x - \tilde{x}|^{\alpha} + ||s - \tilde{s}||_{\bg}^{\alpha})}
\]
where $(s,x)$ are fermi coordinates and $||s - \tilde{s}||_{\bg}$ denotes the distance with respect to the compactified metric. We use $C^{k,\alpha}$ to denote the standard H\"older space with respect to the Euclidean metric. Finally, for any metric space, $(M,g)$, we denote
\[
||f||_{C^{0,\alpha}_g} := \sup_{p \neq q} \frac{|f(p) - f(q)|}{||p - q||_{g}^{\alpha}}
\]
\subsubsection{Showing $u \in C^{1,\alpha}_0$} \label{RescalingArgument}
\begin{figure}[h!]
\centering
\includegraphics[scale=0.3]{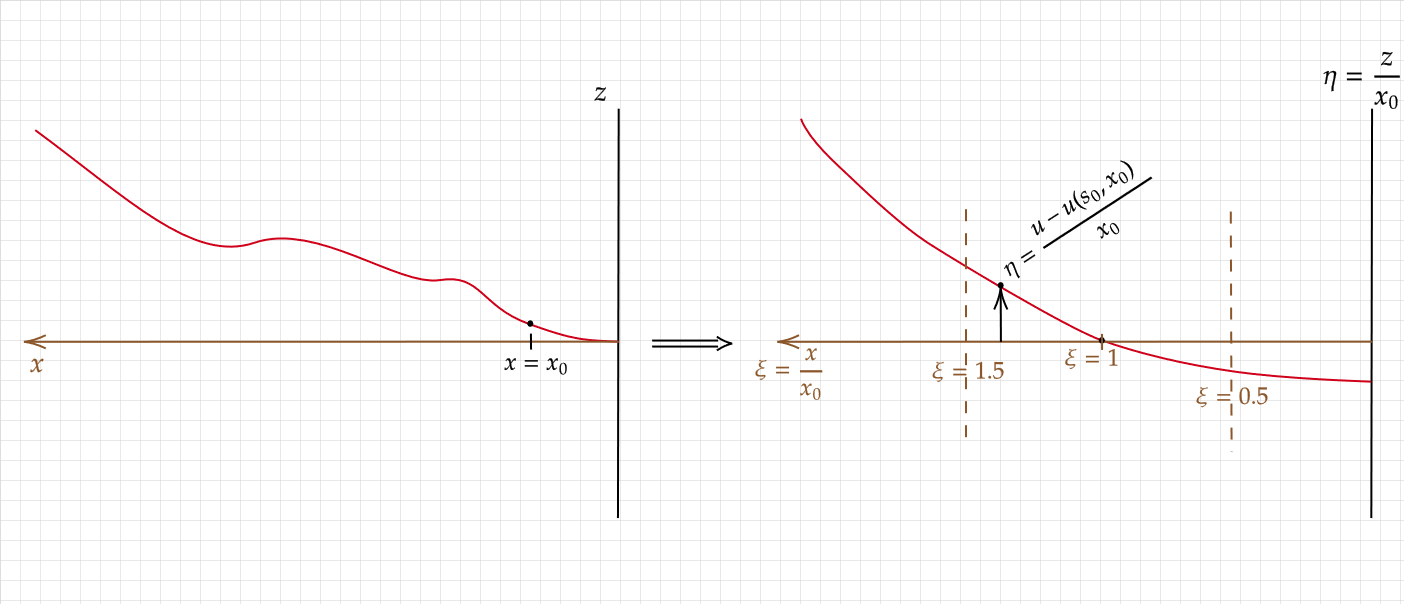}
\caption{Visualization of rescaling with $u = u(\sigma x_0 + s_0, x_0 \xi)$}
\label{fig:rescalingimage}
\end{figure}

Let $p_0 = (s,z,x) = (s_0, u(s_0, x_0), x_0)$, with $x_0$ sufficiently small. We consider rescaled minimal graphs (see figure \ref{fig:rescalingimage}) by changing coordinates
\[
F_{s_0, x_0}: (s,x,z) \mapsto (\sigma, \xi, \eta) := \left(\frac{s - s_0}{x_0}, \frac{x}{x_0}, \frac{z - u(s_0, x_0)}{x_0} \right)
\]
Pulling back the metric by this diffeomorphism, we have 
\begin{align*}
F_{s_0, x_0}^*(g) &= F_{s_0, x_0}^* \left( \frac{dx^2 + k(s,x,z)}{x^2} \right) \\
&= \frac{d \xi^2 + x_0^2 k(x_0 \sigma + s_0, x_0 \xi, x_0 \eta + u(s_0, x_0))}{\xi^2}
\end{align*}
We expand
\begin{align*}
k(s,x,z) &= k_{ab}(s,x,z) ds^{a} ds^b + k_{ai}(s,x,z) ds^a dz^i + k_{ij}(s,x,z) dz^i dz^j \\
x_0^2 k(x_j \sigma + s_0, x_0 \xi, x_0 \eta+ u(s_0, x_0)) &= k_{ab}(x_0 \sigma + s_0, x_0 \xi, x_0 \eta + u(s_0, x_0)) d\sigma^a d \sigma^b \\
&+ k_{ai}(x_0 \sigma + s_0, x_0 \xi, x_0 \eta + u(s_0, x_0)) d \sigma^a d \eta^i  \\
&+ k_{ij}(x_0 \sigma + s_0, x_0 \xi, x_0 \eta + u(s_0, x_0)) d\eta^i d \eta^j
\end{align*}
So that for values of $||\sigma|| \leq 1$, $\frac{1}{2} \leq \xi \leq \frac{3}{2}$, $||\eta|| < \frac{1}{2}$, and $x_0 < \eps$, we have
\begin{align*}
F_{s_0, x_0}^*(g) &= \frac{d \xi^2 + k(s_0,0,0)}{\xi^2} + O(x_0)
\end{align*}
here, we've used that $x_0 < \eps$ and $|z| \leq \frac{1}{2} x_0$, which allows for the above expansion. In particular, we note that 
\begin{align*}
\xi^2 F_{s_0, x_0}^*(g) &= x_0^{-2} F_{s_0, x_0}^*(\bg) \\
&= d \xi^2 + k(s_0,0,0) + O(x_0) \\
&= d \xi^2 + d \sigma_1^2 + \dots + d \sigma_{m-1}^2 + O(x_0)
\end{align*}
The minimal surface then becomes
\begin{align*}
Y_0 &:= F_{s_0, x_0}(Y) \\
(s, x, u(s,x)) \mapsto (\sigma, \xi, \eta) &= \left( \frac{s - s_0}{x_0}, \frac{x}{x_0}, \frac{u(x_0 \sigma + s_0,x_0 \xi) - u(s_0, x_0)}{x_0} \right) \\
&= \left( \frac{s - s_0}{x_0}, \frac{x}{x_0}, O(x_0) \right)
\end{align*}
since $\xi$ and $\eta$ are bounded. Recall the statement of Allard's regularity theorem: 
\begin{theorem}[Allard]
Suppose we have a varifold $V = v(Y^m, \theta) \subseteq \R^{n+1}$, $U$ an open set in $\R^{n+1}$, $\eta > 0$, $\rho_0 > 0$, and $p > 0$, such that for all $a \in \text{spt} ||V||$ and $B_{\rho}(a) \subseteq U$ with $\rho < \rho_0$, we have 
\begin{align*}
1 &\leq \theta  \qquad \mu \; \text{a.e.} \\
\omega_m^{-1} \rho^{-m} ||V||(B_{\rho}(a)) &\leq 1 + \eta \\
\left( \rho^{p-m} \int_{B_{\rho}(a)} |H|^p \right)^{1/p} &\leq \eta
\end{align*}
where $H$ is the generalized mean curvature of the varifold and $p > n$. Then up to a linear isometry of $\R^{n+1}$, $V$ is given graphically by $F = (F^1, \dots, F^{n+1 - m})$ on with
\[
\rho^{-1} \sup |F| + \sup |DF| + \rho^{1-m/p} \sup_{a \neq b} \frac{|DF(a) - DF(b)|}{|a - b|^{1 - m/p}} \leq C \eta^{1/(2m + 2)}
\]
with $C = C(m, n, p)$
\end{theorem}
\noindent \rmk \; Allard Regularity is truly a euclidean theorem, so in order to apply it, we must compute the mean curvature and varifold density in coordinates, with respect to the euclidean metric on the $(\sigma, \xi, \eta)$ coordinates, which we denote as 
\[
g_{euc} = d \sigma_1^2 + \dots + d \sigma_{m-1}^2 + d \xi^2 + d \eta_1^2 + \dots + d \eta_{n+1-m}^2
\]
\noindent We verify the three conditions:
\begin{itemize}
\item $\theta \geq 1$ due to $Y$ being graphical
\item By our assumption \eqref{MassAssumption}
\begin{align*}
\forall \delta > 0, \qquad \exists x^* > 0, \qquad \st \qquad & \forall x_0 < x^* \\
\omega_m^{-1} \rho^{-m} ||V_{Y_0, x_0^{-2} F_{s_0, x_0}^*(\bg)}||(B_{\rho}(a)) &\leq 1 + \delta
\end{align*}
But we've seen that 
\[
x_0^{-2} F_{s_0, x_0}^*(\bg) = g_{euc} + O(x_0)
\]
in $(s,x,z)$ coordinates. Thus
\begin{align*}
\omega_m^{-1} \rho^{-m} ||V_{Y_0, g_{euc}}||(B_{\rho}(a)) &\leq \omega_m^{-1} \rho^{-m}||V_{Y_0, x_0^{-2} F_{s_0, x_0}^*(\bg_0)}||(B_{\rho}(a)) + O(x_0) \\
&\leq 1 + \delta	
\end{align*}
for $x_0$ (and hence $x^*$) sufficiently small
\item We note that
\begin{align*}
H_{Y_0, euc} &= H_{Y_0, x_0^{-2} F_{s_0, x_0}^*(\bg)} + O(x_0) \\
&= H_{Y_0, \xi^2 F_{s_0, x_0}^*(g)} + O(x_0) \\
&= \frac{1}{\xi} \left[H_{Y_0, F_{s_0, x_0}^*(g)} - m \Pi^{NY_0} \nabla(
\ln(\xi) ) + O(x_0)\right] \\
&= \frac{1}{\xi} \left[ -m \Pi^{NY_0} \nabla(
\ln(\xi) ) + O(x_0) \right] \\
\implies |H_{Y_0, euc}| & \leq C \frac{1}{\xi^2} \\
& \leq \tilde{C}
\end{align*}
having applied the formula for (generalized) mean curvature under a conformal change of metric. Here we noted that $H_{Y_0, F_{s_0, x_0}^*(g)} = 0 + O(x_0)$ and $\xi$ is bounded. Also $\Pi^{NY_0}$ denotes the projection onto the normal bundle of $Y_0$ with respect to $F_{s_0, x_0}^*(g)$. Thus $H_{Y_0, euc}$ is bounded, in this rescaled graphical representation. This tells us that any $p$
\begin{align*}
\left( \rho^{p-m} \int_{B_{\rho}(a)} |H|^{p} \right)^{1/p} &\leq \left( \rho^{p-m} \rho^m C^{p} \right)^{1/p} \\
& \leq C \rho
\end{align*}
so choosing $\rho = \delta/C$ gives the desired bound.
\end{itemize}
\noindent Thus, Allard applies and we get the existence of a function
\[
W = (W^1, \dots, W^{n+1-m}) \in C^{1,\alpha}
\]
for $\alpha = 1 - m/p > 0$ with the above bounds. Of course, we already have a graphical description of $Y_0$. Letting 
\[
u_0(\sigma, \xi) := \frac{u(x_0 \sigma + s_0, x_0 \xi) - u(s_0, x_0)}{x_0}
\]
Then, up to an isometry (of euclidean space), $q$, we have $W_0 = q \circ u_0$ and we get the same $C^{1,\alpha}$ bounds for $u_0$. Note that we applied Allard with respect to $(\sigma, \xi)$ coordinates. This gives
\[
||u_0(\sigma, \xi)||_{C^{1,\alpha}(\sigma, \xi)} \leq \delta
\]
but because we're working with $\frac{1}{2} \leq \xi \leq \frac{3}{2}$, the norm computed with $(\sigma, \xi)$ is comparable to the $C^{1,\alpha}_0$ and
\[
||u_0(s,x)||_{C^{1,\alpha}_0} \leq \delta
\]
In our ball corresponding to $\frac{1}{2} \leq \xi \leq \frac{3}{2}$.
\subsubsection{Revamped Schauder Bootstrapping} \label{Revamped}
From the previous section, we have a graphical representation of $Y_0$ in Fermi coordinates. We now consider the metric induced on $Y_0$ as a submanifold of $M^{n+1}$. \textbf{Dropping the $0$ subindex for brevity}, we use the previous sections to write the metric under the $F_{s_0, x_0}$ diffeomorphism as
\[
\overline{h} = Id + M(\nabla u)
\]
where all of the entries of $M$ can be as small as needed by choosing $\delta$ appropriately in our $C^{1,\alpha}$ bounds from \S \ref{RescalingArgument}. Therefore
\[
\overline{h}^{-1} = Id - M + M^2 + \dots = Id - M + O(\delta^2)
\]
Let $H(u)$ denote the mean curavture of the surface given by the graph of $u$ in the $(\sigma, \xi, \eta)$ coordinates. Recall the minimal surface system from Graham and Witten \cite{graham1999conformal}, adapted to $(\sigma, \xi, \eta)$. Here, an $a$ subindex denotes $\sigma_a$, $\xi$ denotes $\xi$, and $i$ denotes $\eta^i$: 
\begin{align} 
0 = \xi H(u)^{k} &= \xi \left[\xi \partial_{\xi} - m + \frac{1}{2} {\xi} L_{\xi}\right] \left[\overline{h}^{\xi\xi} \overline{g}_{ik} u^i_{\xi} + \overline{h}^{a \xi} (\overline{g}_{ak} + \overline{g}_{ik} u^{i}_a ) \right] \label{MinimalEq}\\
& + \xi^2 \left[ \partial_b + \frac{1}{2} L_b\right] \left[ \overline{h}^{\xi b} \overline{g}_{ik} u^{i}_{\xi} + \overline{h}^{ab} ( \overline{g}_{ak} + \overline{g}_{ik} u^i_a) \right]  \nonumber \\
& - \frac{1}{2} x^2 \overline{h}^{ab} \left[\overline{g}_{ab, k} + 2 \overline{g}_{a i, k} u^{i}_b + \bg_{ij, k} u^{i}_a u^{j}_b\right] \nonumber  \\
& -\xi^2 \overline{h}^{a\xi} [\overline{g}_{a i, k} u^i_{\xi} + \overline{g}_{ij,k} u^i_a u^j_{\xi} ] \nonumber  \\
& - \frac{1}{2} \xi^2 \overline{h}^{\xi\xi} [\bg_{ij, k} u^i_{\xi} u^b_{\xi}]\nonumber 
\end{align}
for $L = \log (\det \bh)$. Note that we have multipled by $\xi$ in order to make this a $0$ order differential equation (i.e. can be written in terms of edge operators $(\xi \nabla_{\sigma_a})$ and $(\xi \partial_{\xi})$). Heavily referencing \S \ref{MetricTY}, we can write the above as a quasilinear system of PDEs of the form
\begin{align*}
\xi H(u)^k&=  a_{\alpha \beta i}^k(\xi,\sigma) (\xi \partial_{y_{\alpha}})(\xi \partial_{y_{\beta}}) u^i + g(\xi\nabla u, \xi,\sigma)
\end{align*}
where $\{y_{\alpha}\}$ denote any of $\{\sigma_a, \xi\}$, $g$ is some smooth function, and $\{a_{\alpha \beta i}^k\}$ are uniformly elliptic. In particular, 
\[
a_{\alpha \beta i}^k = \delta_{\alpha \beta} \delta_i^k + O(\delta)
\]
for $\frac{1}{2} \leq \xi \leq \frac{3}{2}$. Alternatively, we frame this as 
\begin{align*}
0 &= \xi^2 \partial_{\alpha}^2 u^k + K_{i}^{\alpha \beta}(\nabla u, \xi, \sigma) \partial_{\alpha} \partial_{\beta} u^i + b_i^{\gamma}(\nabla u, \xi, \sigma) (\partial_{\gamma} u^i) + F(\nabla u, \xi, \sigma) u + G(\nabla u, \xi, \sigma) \\
&= [\delta_{\alpha \beta}\xi^2 + K_i^{\alpha \beta}] \partial_{\alpha} \partial_{\beta} u^i + b_i^{\gamma} \partial_{\gamma} u^i + F u + G
\end{align*}
for some coefficients $K_{i}^{\alpha \beta} = O(\delta)$ and $b_i^{\gamma} = O(\delta)$ for each $k$. Here, $\{K_{i}^{\alpha \beta}\}$ and $\{b_i^{\gamma}\}$ are both $O(\delta)$ in $C^{\alpha}_{(\sigma, \xi)}$ because of their dependence of $\nabla u$ and the fact that $||u(\sigma, \xi)||_{C^{1,\alpha}} = \delta$. Moreover, $F = G = 0$ because the minimal surface system is a divergence system, i.e. it can be written as a collection of equations each of the form 
\[
\text{div} \left(\vec{A}_{\alpha} \cdot \frac{\partial u^{i}}{\partial y_{\alpha}} \right) = 0
\] 
We now apply Schauder estimates using that the functions $K_i^{\alpha \beta}(\nabla u, \xi, \sigma)$, $B_i^{\gamma}(\nabla u, \xi, \sigma)$, and $F(u, \xi, \sigma)$ are all in $C^{\alpha}$
\begin{align*}
||u(\xi, \sigma)||_{C^{2,\alpha}} &\leq \Tau \left( ||u||_{C^{0,\alpha}} + ||G||_{C^{0, \alpha}}\right) = \Tau \left( ||u||_{C^{0, \alpha}} \right)
\end{align*}
See \cite{giaquinta2013introduction} section 5 or \cite{simon1997schauder} for the proof of Schauder estimates in the systems case. We can further improve this:
\[
||u||_{C^{0,\alpha}} \leq ||u||_{C^{1,\alpha}} \leq \tilde{T} ||u||_{C^0} = O(x_j)
\]
having used first order Schauder estimates in the second inequality. We now iterate this argument to get bounds on higher derivatives in terms of the rescaled variables, $(\sigma, \xi)$. This ensures smoothness away from $x = 0$ as well as bounds on higher derivatives in terms of constants independent of $j$. Note that in the above, we've been working with standard H\"older norms in the $(\sigma, \xi)$ variables and the $C^{k,\alpha}$ H\"older norms. But again, because $1/2 \leq \xi \leq 3/2$ we get comparable bounds for the $C^{k,\alpha}_0$ norms for the $(s,x)$ coordinates. \nl \nl
We proved that for $u_0$, there exists a $\delta$ (and hence $x_0$) sufficiently small so that
\[
u_0(\sigma, \xi) \in x_0 \bigcap_k C_0^{k, \alpha} \qquad \forall k \geq 0
\]
Undoing the definition of $u_0(\sigma, \xi)$ in terms of the original function $u$, we get
\[
u(x_0 \sigma + s_0, x_0 \xi) = x_0 u_0(\sigma, \xi) + u(s_0, x_0)
\]
Thus we actually have that $u$ is regular at $x = x_0$ in a neighborhood of radius $x_0$. This construction holds for all $x_0$ sufficiently small, so we conclude
\[
u \in x^2 \bigcap_k C_0^{k, \alpha} \qquad \forall k \geq 0
\]
i.e. 
\[
u = x^2 f, \qquad \forall j, \beta, \quad \exists \; C_{j \beta}^i < \infty \quad \st \quad ||(x \partial_x)^j(x \partial_{s_a})^{\beta} f^i||_{C^{0,\alpha}} \leq C_{j \beta}^i
\]
\subsection{Parametrix Argument} \label{Parametrix}
Having shown that
\[
u \in x^2\bigcap_k C_0^{k, \alpha} 
\]
we now want to show that 
\[
u \in x^2\bigcap_{k} C^{k, \alpha}_b
\]
i.e. for $u = x^2 f$, we have
\[
\forall j, |\beta| \; \in \Z^+, \qquad \exists \; \tilde{C}_{j \beta}^i \; \; \st \;  ||(x \partial_x)^j (\partial_{s_a})^{\beta} f^i||_{C_0^{0,\alpha}} \leq \tilde{C}_{j \beta}^i
\]
(note that we use $\partial_{s_a}$, not $(x \partial_{s_a})$!). To show this, we briefly recall relevant facts from microlocal analysis and the theory of edge operators from \cite{rafe1991elliptic} 
\begin{itemize}
\item The space of conormal functions is
\[
\mathcal{A}:= \bigcap_{k} C^{k, \alpha}_b
\]
\item The space of polyhomogeneous function is
\[
\mathcal{A}_{phg}:= \{u(s,x) \; | \; u(s,x) \sim \sum_{\text{Re}(z_j) \to \infty} \sum_{t = 0}^{N_j}  x^{z_j} \log(x)^t a_{j,t}(s) \; \; \st \; \; N_j < \infty, \qquad a_{j,t} \in C^{\infty} \}
\]
where for $m > |\beta| > 0$
\[
\forall \eps > 0, \qquad \Big| D^{\beta} \left( u(s,x) - \sum_{\text{Re}(z_j) < m} \sum_{t = 0}^{N_j} x^{z_j} \log(x)^t a_{j,t}(s) \right) \Big| = O(x^{m - |\beta| - \eps})
\]
i.e. the remainder and its derivatives (in $x$ and $s_a$ variables) decay at a faster rate. In practice, we'll be dealing with $z_j$ real, positive, and integer valued.

\item We denote the space of edge operators, $V_e$, and the space of $b$-operators, $V_b$ as 
\begin{align*}
V_e &= \text{span}_{C^{\infty}(\mathcal{U})} \{(x \partial_x), (x \partial_{s_a})\} \\
V_b &= \text{span}_{C^{\infty}(\mathcal{U})} \{(x \partial_x), \partial_{s_b}\}
\end{align*}
for $\mathcal{U}$ the neighborhood of $\Gamma$ as defined in \S \ref{VariationsOfRV}

\item The weighted H\"older space of orders $\ell \in \N$, $\alpha \in (0,1)$, and $\delta \in \R$ are 
\[
x^{\delta} \Lambda^{\ell, \alpha}_0:= \{u = x^{\delta} v \; \st \; V_1 \cdots V_j v \in \Lambda^{0, \alpha}_0, \qquad \forall V_i \in V_e, \;\; j \leq \ell\}
\]
where $\Lambda^{0, \alpha}_0 = C^{0,\alpha}$ is the geometric H\"older space with norm
\[
||f||_{0, \alpha, 0} = \sup |f| + \sup_{(s,x) \neq (\tilde{s}, \tilde{x})} \frac{ (x + \tilde{x})^{\alpha} |f(s, x) - f(\tilde{s}, \tilde{x})|}{|x - \tilde{x}|^{\alpha} + ||s - \tilde{s}||_{\bg}^{\alpha}} 
\]
We also define
\[
x^{\delta} \Lambda^{\ell, \alpha, m}_0 = \{u = x^{\delta} v \; \st \; (V_1 \cdots V_j)(\tilde{V}_1 \cdots \tilde{V}_k) v \in \Lambda^{0, \alpha}_0, \qquad  V_i \in V_e, \; \tilde{V}_i \in V_b, \qquad \;\; j \leq \ell, \; k \leq m\}
\]
\item $\Psi^{m, \mathcal{E}}$ denotes pseudodifferential operators in the small calculus. $\Psi^{m, \mathcal{E}}_0$ denotes the large calculus. $\Psi^{m, \mathcal{E}}_b$ denotes the analogous calculus but with respect to the $b$-operators

\item For $L \in \Psi^{m, \mathcal{E}}_0$ an elliptic pseudodifferential operator, a parametrix, $G \in \Psi^{-m, \mathcal{E}}_0$, exists such that 
\begin{align*}
LG &= I - R_1 \\
GL &= I - R_2
\end{align*}
where $I$ is the identity and $R_1, R_2 \in \Psi^{-\infty, \mathcal{E}}$ are ``residual" operators. Here, $m$ is the order of the principal symbol of $L$. Roughly speaking, $R_i$ sends functions of any regularity into polyhomogeneous functions.

\end{itemize}
We also recall a few relevant propositions from \cite{rafe1991elliptic} adopted for our case of the index set $\mathcal{E} = \{0, m + 1\}$
\begin{itemize}
\item (Proposition $3.27$) For $A \in \Psi^{m, \mathcal{E}}_0$, suppose $\ell \geq m$ and $\delta > -1$, then 
\[
A: x^{\delta} \Lambda^{\ell, \alpha}_0 \to x^{\delta} \Lambda^{\ell - m, \alpha}_0
\]

\item (Proposition $3.28$) For $f \in \mathcal{A}_{phg}$, $A \in \Psi_0^{m, \mathcal{E}}$, we have $Af \in \mathcal{A}_{phg}$

\item (Proposition $3.30$) For $v \in V_b$, $A \in \Psi_0^{m, \mathcal{E}}$, we have $[v, A] \in \Psi_{e}^{m, \mathcal{E}}$
\end{itemize}
\noindent \rmk \; Hereafter, we use $O(x^k)$ to denote a remainder term which lies in $x^k \mathcal{A}$, and $o(x^k)$ to a remainder term, $f$, such that $\lim_{x \to 0} x^{-k} f = 0$ with convergence in $\mathcal{A}$. \nl \nl
We now prove that $u \in x^2 \mathcal{A}$. We argue as follows: for $\bN_i = \partial_{z_i}$ a basis for $N(\Gamma)$ ($z_i$ are Fermi coordinates for the normal bundle), we have $u = u^i \bN_i$ and
\begin{align*}
0 &= H(u) = H(0) + L_0(u) + Q_0(u) \\
\implies L_0(u)^i &= -H(0)^i - Q_0(u)^i
\end{align*}
where $L_0 = J_{\Gamma}$ is the Jacobi operator on $\Gamma \subseteq M$ and also the linearization of the mean curvature functional about $u = 0$. $Q_0$ is the quadratic remainder from the linearization and depends on $\{x \nabla u, x, s \}$, which are parameters in the coefficients for our elliptic system of equations which we've bounded in the previous section. Here $Q_0(u) = O(x^4)$ because $u = O(x^2)$, and the superscript denotes the $\bN_i$th component. We have
\[
L_0(u) = \Delta_{\Gamma}(u) + \tilde{A}_{\Gamma}(u) - (m + E(x)) u
\]
Here, $\Delta_{\Gamma}$ is the Laplacian on $\Gamma$ on the \textit{normal bundle} computed with respect to the $(\xi, \sigma)$ variables, $\tilde{A}$ is the Simons operator, a $0$th order operator that is $O(x^2)$ (see \S \ref{AsymptoticVariational}), and $E(x) = O(x^2)$ is an error term coming from the computation of $\text{tr}_{\Gamma} [R_M(\cdot, u) \cdot]$ as in the standard Jacobi operator. Let $G$ be a parametrix for this operator
\begin{align*}
GL &= I - R \\
GL(u)^i &= u^i - (Ru)^i \\
& = - (G H(0))^i + (G Q)^i \\
\implies u^i &= -(G H(0))^i + (GQ)^i + (R u)^i \\ 
& = GQ^i + F^i
\end{align*}
where $Ru$ is a residual term. One can compute $H(0) = x^2 H_{\gamma} + O(x^4)$ analytic in $x$ and $s$. By Propositions $3.27$ and $3.28$, we have that $G H(0)$ is $O(x^2)$ and polyhomogeneous. Moreover, because $R$ is residual, we have that $Ru$ is polyhomogeneous. Finally by $3.27$, we know that $GQ = O(x^4)$. With this, we can write 
\[
u^i + (G H(0))^i - (GQ)^i = (Ru)^i
\]
and note that the left hand side is $O(x^2)$, so $(Ru)^i$ must be both $O(x^2)$. This tells us that $F^i := -(G H(0))^i + (R u)^i$ is $O(x^2)$ and polyhomogeneous.  
We now differentiate this equation to get
\[
\partial_{s_a} u^i = \partial_{s_a} F^i +  G( x^{-1} (x\partial_{s_a}) Q) + [\partial_{s_a}, G] Q
\]
Again, $\partial_{s_a} F^i$ is $O(x^2)$ by polyhomogeneity. From our initial estimates, we have that $x^{-1} (x\partial_{s_a}) Q = O(x^3)$, and by Proposition $3.27$, $G (x^{-1} (x\partial_{s_a}) Q) \in x^3 \cap_k C^{k, \alpha}_0$. Similarly, by $3.30$ and $3.27$, we have that $[\partial_{s_a}, G] Q = x^4 \bigcap_{k} C^{k,\alpha}_0$ because $Q = O(x^4)$. This shows that 
\[
u \in \bigcap_k x^2 \Lambda^{k, \alpha, 1}_0
\] 
We now proceed by induction. Assume that
\[
u \in x^2 \bigcap_k \Lambda^{k, \alpha, j}_0
\]
For $\alpha$ a multi-index of order $j + 1$, we write
\[
\partial_{s_{a_1}} \cdots \partial_{s_{a_{j+1}}} u = \partial_s^{\alpha} u = \partial_s^{\alpha} (GQ^i + F^i) = \partial_{s}^{\alpha} (GQ^i) + \partial_s^{\alpha} (F^i)
\]
we automatically have that $\partial_s^{\alpha} F^i \in x^2 \Lambda_0^{k,\alpha, j}$ for any $k$ and $j$ since $F^i$ is polyhomogenous. For the first term, we write 
\begin{align*}
\partial_s^{\alpha} &= \partial_{s_a} \partial_s^{\beta}, \qquad |\beta| = j\\
\partial_s^{\beta} (GQ^i) & = \sum_{|\gamma| + |\delta| = j} c_{\gamma} [ \partial_{s_{\gamma_1}}, \cdots, [\partial_{s_{\gamma_{j}}}, G] ] \partial_s^{\delta} Q^i
\end{align*}
where $c_{\gamma}$ is some integer valued coefficient reflecting the combinatorics of how many commutator terms we get. By induction and the chain rule, we know that $\partial_s^{\delta} Q^i$ is $O(x^4)$ for all $\delta$. By repeated application of Proposition $3.30$, we know that the nested commutator term lies in $\Psi^{-2, \mathcal{E}}_e$. Then by Proposition $3.27$ and $3.28$, we can conclude that 
\[
c_{\gamma} [ \partial_{s_{\gamma_1}}, \cdots, [\partial_{s_{\gamma_{j}}}, G] ] \partial_s^{\delta} Q^i \in x^4 \bigcap_k C^{k ,\alpha}_0
\]
so that 
\[
\partial_{s_a} \left( c_{\gamma} [ \partial_{s_{\gamma_1}}, \cdots, [\partial_{s_{\gamma_{j}}}, G] ] \partial_s^{\delta} Q^i \right) = x^{-1} (x \partial_{s_a}) \left( c_{\gamma} [ \partial_{s_{\gamma_1}}, \cdots, [\partial_{s_{\gamma_{j}}}, G] ] \partial_s^{\delta} Q^i \right) \in x^3 \bigcap_k C^{k ,\alpha}_0
\]
adding the $F^i$ term we have
\[
u \in x^2 \bigcap_k \Lambda_0^{k,\alpha, j + 1}
\]
This completes the induction and we get
\begin{align} \nonumber
u &\in x^2 \bigcap_{k,j} \Lambda_0^{k,\alpha, j} = x^2 \bigcap_k C^{k, \alpha}_b \\ \label{xsquareddecay}
\implies u & \in x^2 \mathcal{A} 
\end{align}
%

%
\subsection{Expanding Mean Curvature Functional} \label{ExpandMC}
%
%
We now compute the linearization of the mean curvature functional, $H$, on graphical submanifolds of the form $\{\vec{u}(s,x)\}$ from before. We first linearize about $u_0 = 0$:
\[
H(u) = H(0) + L(u) + Q(u)
\]
so that in \eqref{MinimalEq}, we set $0 = z^i = p^i_{\alpha}$ when evaluating $\bh^{\alpha\beta}$ as abstract functions of $(x,\{z^i\},\{p^i_{\alpha}\})$ to be set equal to $(x, \{u^i\}, \{ u^i_{\alpha}\})$. Here, $Q(u)$ is an expression that's at least quadratic in the components of  $\{u, x \nabla u\}$ and depends smoothly on $s$. Because $u \in x^2 \mathcal{A}$, we have $Q(u) \in x^4 \mathcal{A}$. Note that $H(0)$ is the mean curvature of the graph corresponding to $\vec{u}(s,x) = 0$, which is just $\Gamma = \R^+ \times \gamma$. A short computation gives
\begin{align*}
H(0) &= H_{\{u =0\}} \\
& = [x^2 H_{\gamma}^i + R^i] \bN_i \\
R^i &= O(x^4) \\
\F(R^i) &= 1
\end{align*}
where $H_{\gamma}^i$ are the components of the mean curvature of the boundary submanifold, computed with the compactified metric restricted to the boundary. In particular, we note that 
\begin{equation} \label{MCErrorTerm}
R^i = \begin{cases}
R^i_4 x^4 + \dots + R^i_{n+2} x^{n+2} + R^i x^{n+2} \log(x) + O(x^{n+3} \log(x)) & n \text{ even } \\
R^i_4 x^4 + \dots + R^i_{n+2} x^{n+2} +  O(x^{n+3}) & n \text{ odd}
\end{cases}
\end{equation}
With this, we note that 
\begin{align*}
L & = (\text{Jacobi operator evaluated at $\{u = 0\}$ graph}) \\
& = \Delta_{\{ u = 0\}} + \tilde{A}_{\{u = 0\}} + \text{Tr}[R_M(\cdot, \_) \cdot] = \Delta_{\Gamma} + \tilde{A}_{\Gamma} + \text{Tr}[R_M(\cdot, \_) \cdot] 
\end{align*}
for $\tilde{A}$ the Simons operator. Here, let 
\[
\begin{gathered}
u = u^i \bN_i = \hat{u}^i \hat{N}_i \\
\hat{u}_i = x^{-1} u^i, \qquad \hat{N}_i := x \bN_i
\end{gathered}
\]
Note that $g(\hat{N}_i, \hat{N}_i) = 1 + O(x^2)$ on $\Gamma$. This choice of notation is so that geometric operators with respect to $g$ act more naturally on $\hat{u}^i$ as opposed to $u^i$ due to the choice of normalization. Following the work of \cite{graham1991einstein} (corollary 2.8) along with \S \ref{Simons}, we can write 
\[
L(u) = [(x \partial_x)^2 - (m-1)(x \partial_x) - m](\hat{u}^i) (\hat{N}_i) + E(x^{-1} u)^i (x\hat{N}_i)
\]
where $E$ is an error term that is at most second order in $V_b$ operators and has $O(x^2)$ coefficients and analogously to \eqref{MCErrorTerm} can only have $x^{k} \log(x)$ terms at $k \geq n+2$. Thus $E(u)^i = O(x^3)$.
Via \S \ref{Simons}, we see that $\tilde{A}_{\{u = 0\}}(u) = f_i \hat{N}_i$ where $f_i = O(x^3)$. With this, we begin our iteration at 
\begin{align} \nonumber
H(u) & = 0 = H(0) + L(u) + Q(u) \\ \label{MCLinearization}
&= x^2 H_{\gamma} + [(x \partial_x)^2 - (m-1)(x \partial_x) - m](\hat{u}^i) \hat{N}_i  + O(x^4)
\end{align}
\subsection{Iteration Argument}
\label{IterationArgument}
To begin the proof of the theorem \ref{AsymptoticExpansion}, we first establish the lowest order coefficient in the expansion
\begin{lemma} \label{u2Lemma}
The minimal submanifold, $Y^m$, can be described as a graph over $\Gamma = \partial Y \times [0, \epsilon)$ via 
\[
Y \cap \{0 \leq x < \epsilon\} = \{ \overline{\exp}_{(\gamma(s), x)}(u(s,x)) \quad | \quad 0 \leq x < \epsilon, \qquad \gamma(s) \in \gamma\}
\]
where $u(s, x) = u^i(s,x) \bN_i = u^i(s,x) \partial_{z_i}$ and
\[
\boxed{ u^i(s,x) = \frac{1}{2(m-1)} H_{\gamma}^i x^2 + O(x^3) }
\]
Where $\partial_{z_i}$ is a Fermi coordinate basis for the normal bundle with respect to the compactified metric and $H_{\gamma}$ is the mean curvature of $\gamma \subseteq \partial M$, and $m \geq 2$.
\end{lemma}
\noindent \Pf \; Having extracted the linear term in equation \eqref{MCLinearization} and shown that the remainder is $O(x^2)$, we write
\begin{align*}
L(u) &= -H(0) - Q(u) \\
&= - x H_{\gamma}^i \hat{N}_i - Q(u) - R_0 \\
&= [-x H_{\gamma}^i + O(x^3) ]\hat{N}_i
\end{align*}
Hence
\begin{align*}
((x \partial_x)^2 - (m-1) (x \partial_x)  - m) (\hat{u}^i)\hat{N}_i & = - x H_{\gamma}^i \hat{N}_i + [-Q(u) + O(x^4)] \hat{N}_i \\
\implies  ((x \partial_x)^2 - (m-1) (x \partial_x) - m) (\hat{u}^i) &  = -x H_{\gamma}^i + O(x^3)
\end{align*}
we can write this second order operator in $x \partial_x$ as a composition of first order operators:
\[
((x \partial_x)^2 - (m-1) (x \partial_x) - m) = (x \partial_x + 1) (x \partial_x - m)
\]
We then use integrating factors as follows: define
\[
f^i = (x \partial_x - m) \hat{u}^i \in x^2 \mathcal{A}
\]
which means that 
\begin{align*}
((x \partial_x)^2 - (m-1) (x \partial_x) - m) (\hat{u}^i) &  = -x H_{\gamma}^i + O(x^3) \\
\implies (x \partial_x + 1) f^i &= -x H_{\gamma}^i + O(x^4) \\
\partial_x(x f) &= -x H_{\gamma}^i + O(x^3)
\end{align*}
now integrating from $0$ to an arbitrary value of $x$ and dividing by $x$, we have
\begin{align*}
f^i &= - \frac{x}{2} H_{\gamma}^i + O(x^4)
\end{align*}
Repeating the same procedure on $f$, but first multiplying by $x^{-m-1}$, we conclude
\begin{align*}
\hat{u}^i &= \frac{1}{2(m-1)} H_{\gamma}^i x + O(x^3)
\end{align*}
%
This is valid when $m \geq 3$ since we absorb $K x^{m}$ into $O(x^3)$. Converting back to $u^i = x \hat{u}^i$, this process gives an explicit formula for $u_2^i(s)$. \qed \nl \nl
\noindent We now want to iterate this argument to get an even expansion up to $m$, with a potential log term when $m$ is odd. \nl \nl
\textbf{Proof of Theorem \ref{AsymptoticExpansion}}: \nl
\noindent We first do $m$ even. Assume the inductive hypothesis of
\[
\hat{u}^i = p_{2k-1}^i(x) + f^i(x), \qquad \qquad f^i(x) = o(x^{2k-1}) \qquad 2k < m
\]
where $p_{2k-1}^i$ is an odd polynomial in $x$ of order $2k-1$ with coefficient depending smoothly on $s$ and $f^i \in \mathcal{A}$. Further assume
\[
H(p_{2k-1}^i \hat{N}_i) = O(x^{2k+1}) = [a_{2k+1}^i x^{2k+1} + o(x^{2k+1})]\hat{N}_i
\]
We have established the base case, $k = 0$ with $p_0^i = 0$ and $a_2^i= \frac{1}{2(m-1)} H_{\gamma}^i $. For higher values of $k$, we can expand
\[
H([p_{2k-1}^i + f^i]\hat{N}_i ) = H((p_{2k-1}^i \hat{N}_i) + L_{p_{2k-1}^i \hat{N}_i}(f^i \hat{N}_i) + Q_{p_{2k-1} \hat{N}_i}(f^i \hat{N}_i )
\]
Abbreviate $L_{2k-1} := L_{p_{2k-1}^i \hat{N}_i}$. This is the linearized operator (i.e. Jacobi Operator) corresponding to the graph of $\{u = p_{2k-1}^i \hat{N}_i\}$. Then using the fact that $p_{2k-1}^i = O(x)$ produces a graphical asymptotically hyperbolic manifold with odd coefficients up to at least order $x^2$, we have as before
\begin{align*}
L_{2k-1} &= I_{L_{2k-1}} + T_{L_{2k-1}}  \\
I_{L_{2k-1}} &= [(x \partial_x)^2 - (m-1)(x \partial_x) - m] \\
T_{L_{2k-1}}&: x^k \mathcal{A} \to x^{k+2} \mathcal{A} \\
\implies L_{2k-1}(f^i \hat{N}_i) &= [(x \partial_x)^2 - (m-1) (x \partial_x) - m](f^i) \hat{N}_i + o(x^{2k + 1})
\end{align*}
where ``$I_{L_{2k-1}}$" stands for the indicial operator of the linearization at $p_{2k-1}$ and $T_{L_{2k-1}}$ is the remainder. Finally, $Q_{2k-1}(f^i \hat{N}_i) := Q_{p_{2k-1}^i \hat{N}_i}(f^i \hat{N}_i)$ will be at least $x^2$ times the order of $f^i$ and hence of order $o(x^{2k + 1})$. Thus
\begin{align*}
0 & = H(\hat{u}^i \hat{N}_i) = H([p_{2k-1}^i + f^i] \hat{N}_i ) \\
& = H(p_{2k-1}^i \hat{N}_i) + L_{2k-1}(f^i \hat{N}_i) + Q_{2k-1}(f^i \hat{N}_i) \\
& = a_{2k+1}^i x^{2k+1} \hat{N}_i + [(x \partial_x - m) (x \partial_x + 1)f^i] \hat{N}_i + o(x^{2k + 1}) 
\end{align*}
Rearranging and matching vector components, we get
\[
(x \partial_x + 1) (x \partial_x - m) f^i = -a_{2k+1}^i x^{2k+1} + o(x^{2k + 1})
\]
as before, we perform an integrating factor for $(x\partial_x + 1)$ first and then $(x \partial_x - m)$
\begin{align} \nonumber
x^{-m-1} (x \partial_x - m) f^i &= -\frac{a_{2k+1}^i}{2k+1}x^{2k-m} + o(x^{2k-m}) \\ \label{IterationIntermediate}
\partial_x (x^{-m} f^i) &= -\frac{a_{2k+1}^i}{2k+1} x^{2k-m}  + o(x^{2k-m})\\ \nonumber
f^i &= -\frac{a_{2k+1}^i}{(2k-m+1)(2k + 1)}x^{2k+1} +  K x^{m} + o(x^{2k+1})
\end{align}
we see that the denominators are never $0$ when $m$ is even. Note that $K$ is the constant from evaluating $x^{-m} f$ at some point $x = x_0$ small but non-zero. This shows that we can continue to induct and get the next even term in our expansion as long as  $2k < m$.  \nl \nl
When $2k = m$, the above process shows that $f^i = K' x^{m+1} + O(x^{m + 1})$ and we can continue the expansion but the expansion is no longer even. Converting back to $u^i = x \hat{u}^i$, we have
\begin{align*}
m \text{ even } \implies u^i &= \frac{1}{2(m - 1)} H_{\gamma}^i x^2 + \dots + u_{m}^i x^{m} + u_{m + 1}^i x^{m+1} + O(x^{m+2})
\end{align*}
i.e. $\F(u^i) = 1$. In particular our remark about \eqref{MCErrorTerm} shows that when $m = n$ even, there is no $x^n \log(x)$ term because $x^k \log(x)$ error terms occur for $k \geq n+ 1$ in the iteration. \nl \nl
When $m$ is odd, most of the proof remains the same. However, when $2k = m - 1$, we see that \ref{IterationIntermediate} becomes
\begin{align*}
\partial_x (x^{-m} f^i) &= -\frac{a_{m}^i}{m} x^{-1}  + o(x^{-1}) \\
x^{-m} f^i &= K - \frac{a_{m}^i}{m} \log(x) + o(\log(x)) \\
f^i &= K x^m - \frac{a_{m}^i}{m} x^{m} \log(x) + o(x^{m} \log(x) ) 
\end{align*}
so a log term appears in this case. After setting 
\[
p_{2k + 1}^i = p_{2k-1}^i -\frac{a_{m}^i}{m}x^{m} \log(x) + K x^{m}
\]
we can continue the iteration without $\log$ terms but we lose evenness of the expansion. Converting back to $u^i$, we conclude
\begin{align*}
m \text{ odd } \implies u^i &= \frac{1}{2(m - 1)} H_{\gamma}^i x^2 + \dots + u_{m+1}^i x^{m+1} + U^i x^{m+1} \log(x) + u_{m+2} x^{m+2} + O(x^{m+3})
\end{align*}
which is again, the statement of $\F(u^i) = 1$. \qed 
\begin{Remark} \label{LocalTermsuExpansion}
In the above, we note that the log coefficient, $\{u_{2k}\}_{k = 1}^{(m-1)/2}$ and $U^i(s)$, are \textit{locally} determined, i.e. they are polynomials of derivatives of $u_2(s)$, which is a direct consequence of the above computation. \nl \nl
\end{Remark}
\noindent Our proof also gives the following result about the induced metric on $Y$
\begin{corollary} \label{MetricExpansion}
The induced metric on $Y$ satisfies for
\begin{align*}
\F(\overline{h}_{\alpha \beta}) &= (-1)^{\sigma(\alpha) + \sigma(\beta)} \\
\overline{h}_{xx} &= 1 + O(x^2) \\
\overline{h}_{xa} &= O(x^3) \\
\overline{h}_{ab} &= \delta_{ab} + O(x^2)
\end{align*}
where $\sigma(\alpha)$ is as in equation \eqref{SigmaNotation}.
\end{corollary}
\noindent \rmk \; We can take the analysis further for $m = n$ even:
\begin{corollary} \label{NoLogTerms}
For $m = n$ even,
\begin{align*}
[\bh_{ab}]^{\log, n +1} &= 0 \\
[\bh_{xx}]^{\log, n} &= 0 \\
[\bh_{xx}]^{\log, n + 1} &= 0 \\
[\sqrt{\det \bh}]^{\log, n + 1} & = 0
\end{align*}
\end{corollary} 
\noindent This follows from the explicit formula for $v_a$ and $v_x$ in \eqref{TYMetricAsymptotics} and holds for all coefficients $K$ in the $K x^n \log(x)$ term of \eqref{kEquation}. Note that $\bh_{ab}$ and $\sqrt{\det \bh}$ may have an $x^n \log(x)$ term.
\section{Parity of Second Fundamental Form} \label{SFFParity}
In this section we aim to prove the following theorem:
\begin{theorem} \label{SFFParityMinimal}
Suppose that $Y^m \subseteq M^{n+1}$ minimal with $\overline{h} = \overline{g} \Big|_{TY}$ even up to order $x^m$. Let $p \in Y$ and $\bA: \text{Sym}^2(TY) \to N(Y)$ denote the second fundamental form, and let $\{w_i(s,x)\}$ be the frame for the normal bundle described in \S \ref{NormalFrame}. Define
\[
T_{\gamma \delta i; \alpha_1 \cdots \alpha_p} := \overline{g}\left((\bnabla_{v_{\alpha_1}} \cdots \bnabla_{v_{\alpha_p}} \bA)(v_{\gamma}, v_{\delta}), w_i \right)
\]
where $\alpha_i$ can take on any of the indices $\{s_{1}, \dots, s_{m-1}, x\}$. Let $q$ denote the number of ``$x$"s among the indices $\{\gamma, \delta, \alpha_1, \dots, \alpha_p\}$, then we have 
\[
\forall i, \qquad \F(T_{\gamma \delta i; \alpha_1 \cdots \alpha_p}) = (-1)^q
\]
\end{theorem}
\noindent We notate the following
\begin{align} \label{SFFDefinition} 
\bA_{\alpha \beta} &= \bg(\bnabla_{F_{\alpha}} F_{\beta}) \\ \label{ChristoffelDefinition}
\tilde{\bGamma}_{\sigma \tau \omega} &:= g(\bnabla_{v_{\sigma}} v_{\tau}, v_{\omega}) \\ \label{MetricDefinitions}
\bh_{\alpha \beta}(t) &:= \bg(F_{\alpha}, F_{\beta}) 
\end{align}
%
where $\bnabla$ is the connection with respect to $\bg$. 
We also define $\bA_{\alpha \beta i}$, the components of $\bA_{\alpha \beta}$, and $\bA_{\alpha \beta}^i$ and $\tilde{\bGamma}_{\sigma \tau}^{\omega}$ by raising the tensors appropriately. Finally, we let the indices $\{\sigma, \tau, \omega, \mu\}$ denote any vector in the basis for $TM = T Y \oplus N Y$, i.e. $v_{\sigma}, v_{\tau}, v_{\omega} \in \{v_{s_a}, v_x, w_i\}$ and similarly with $F_{\alpha} \in \{F_a, F_x\} \subseteq TY_t$. For example 
\[
\tbGamma_{\alpha \beta i} = g(\bnabla_{v_{\alpha}} v_{\beta}, w_i)
\]
Note that 
\[
\bA_{\alpha \beta} = \tbGamma_{\alpha \beta}^j w_j
\]
for $\alpha$, $\beta$ denoting $\va, \vb \in TY$.
\subsubsection{Lemmas for theorem \ref{SFFParityMinimal}}
We start by writing the tangent basis for $Y$ in fermi decomposition
\begin{align*}
v_a &= G_*(\partial_{s_a}) = \partial_{s_a} + u_a^i \partial_{z_i} \\
v_x &= G_*(\partial_x) = \partial_x + u_x^i \partial_{z_i}
\end{align*}
Similarly, recall the parity of the coefficients for our normal frame: (cf. section \S \ref{MetricTY} and lemma \ref{NormalFrame})
\begin{align*}
w_i &= c^{a}_i(s,x) \partial_{s_a} + c^{x}_i(s,x) \partial_x + c^{j}_i(s,x) \partial_{z_j} \\
\F(c^{a}_i) &= 1 \\
\F(c^{x}_i) &= -1 \\
\F(c^{j}_i) &= 1
\end{align*}
Let $\bGamma_{\sigma \tau \omega}$ denote the christoffel symbols in the basis of $\{\partial_{s_a}, \partial_x, \partial_{z_i}\}$, with respect to $\bg$, in a tubular neighborhood of $\Gamma = \gamma \times \R^+$, parameterized by $(s,x, z)$. 
\begin{lemma} \label{ChristoffelCylinderLemma}
We have
\begin{align} \label{ChristoffelParityCylinder}
\forall p \in Y, \qquad \F\left(\bGamma_{\sigma \tau\omega}\Big|_{p \in Y}\right) = (-1)^q
\end{align}
where $q$ is the number of indices among $\sigma$, $\tau$, $\omega$ that are equal to $x$
\end{lemma}
\noindent \Pf First note that via the fermi coordinate decomposition
\[
\bGamma_{xx b} = \bGamma_{xxi} = 0
\]
And also
\[
\bGamma_{\alpha i\beta} = \bg(\bnabla_{\partial_{y_{\alpha}}} \partial_{z_i}, \partial_{y_{\beta}}) = - \bg(\partial_{z_i}, \bnabla_{\partial_{y_{\alpha}}} \partial_{y_{\beta}}) = - \bGamma_{\alpha \beta i}
\]
and so it suffices to consider $\bGamma_{\cdot \cdot i}$, $\bGamma_{axb}$, $\bGamma_{abx}$, and $\bGamma_{abc}$. The proof of the result comes from the splitting of the ambient metric under Graham-Lee Normal form in a tubular neighborhood of the boundary (i.e. on $\partial M \times [0, \eps)$): 
\[
\bg = dx^2 + k(x,s,z)
\]
Here, $k(x,s,z)$ is a $2$-tensor such that $k(x,s,z)(\partial_x, \cdot) \equiv 0$. Moreover $k(s,x,z)$ expands as
\begin{align*}
\text{$n+1$ even} \implies k(s,x) &= k_0 + x^2 k_2 + \cdots + k_{n-1} x^{n-1} + k_{n} x^{n} + k_{n+1} x^{n+1} + O(x^{n+2}) \\
\text{$n+1$ odd} \implies k(s,x) &= k_0 + x^2 k_2 + \cdots + k_{n} x^{n} + Kx^{n} \log(x)  + O(x^{n+2}) 
\end{align*}
where each $k_i = k_i(s,z)$ is a $2$-tensor. $k_i(s,z)$ can also be expanded up to order $z^m$ in $z$ since $z = (z^1, \dots, z^{n+1-m})$ is a system of fermi coordinates and $\gamma$ is $C^{m+1, \alpha}$ embedded in $\partial M$. so we can expand 
\begin{align*}
k_i(s,z) &=  k_{i,0}(s) + k_{i,1}(s) z + \dots + k_{i,p}(s) z^m + O(z^{m+1})
\end{align*}
for any $z$. Evaluating this on $Y$ (i.e. $z = u(s,x)$), we see that $k(s,x,u(s,x))$ is a tensor that is even in $x$ up to order $x^{m+1}$. With this and the Koszul formula, one can directly show the parity statements in \eqref{ChristoffelParityCylinder} hold in a tubular neighborhood of $\gamma \times [0, \eps) \subseteq \partial M \times [0, \eps)$. \qed \nl \nl 
We now extend this to compute the Christoffel's in the basis of $\{v_a, v_x, w_i\}$: Let $\tbGamma$ be as in \eqref{ChristoffelDefinition} and $\tbGamma^{\omega}_{\sigma \tau}$ the raised versions of these christoffels by the induced metric on $Y$. 
\begin{lemma} \label{YChristoffelParity}
For $\tbGamma_{\sigma \mu \omega}$ as above evaluated on $Y$, we have that 
\[
\forall p \in Y, \qquad \F\left(\tbGamma_{\sigma \mu \omega}\Big|_{p \in Y}\right) = (-1)^q
\]
where $q$ the number of $x$'s among the indices $\sigma$, $\mu$, $\omega$.
\end{lemma}
\noindent \Pf Again, this boils down to recording parity of the coefficients of $\{v_a, v_x, w_i\}$ in the basis of $\partial_{s_a}, \partial_x, \partial_{z_i}$. We'll compute the first christoffel and leave the remainder to the reader
\begin{align*}
\bnabla_{v_a} v_b &= \bnabla_{\partial_{s_a}}(\partial_{s_b} + u_b^j \partial_{z_j}) \\
& + u_a^i \bnabla_{\partial_{z_i}} (\partial_{s_b} + u_b^j \partial_{z_j}) \\
&= A_1 + A_2 
\end{align*}
We have
\begin{align*}
A_1 &= \bGamma_{ab}^{\omega} \partial_{y_{\omega}} + u_{ab}^j \partial_{z_j} + u_b^j \bGamma_{bj}^{\sigma} \partial_{y_{\sigma}}\\
A_2&= u_a^i \Big[ \bGamma_{ib}^{\omega} \partial_{y_{\omega}} + u_b^j \bGamma_{ij}^{\omega} \partial_{y_{\omega}}\Big] 
\end{align*}
One can now compute using lemma \ref{ChristoffelCylinderLemma} that 
\begin{align*}
\F(g(A_i, \partial_{s_c})) &= 1 \\
\F(g(A_i, \partial_{z_j})) &= 1  \implies \F(u^j g(A_i, \partial_{z_j})) = 1\\
\F(g(A_i, \partial_{y_{\mu}})) &= (-1)^{\sigma(\mu)} \implies \F(u \Gamma_{ci}^{\mu} g(A_i, \partial_{y_{\mu}})) = 1 \\
\implies \F(g(A_i, v_c))&= 1
\end{align*}
which verifies that $(\Gamma_{abc}) = 1$. The remaining symbols proceed similarly. \qed \nl \nl
Finally, we establish a short lemma about the metric in the basis of $\{v_a, v_x, w_i\}$:
\begin{align*}
\bW_{\alpha \beta} &= \bg(v_{\alpha}, v_{\beta}) \\
\bW_{\alpha i} &= \bg(v_{\alpha}, w_i) \\
\bW_{ij} &= \bg(w_i, w_j)
\end{align*}
and $\bW^{\gamma \delta}$ is defined as the inverse.
\begin{lemma} \label{WMetricParity}
For $\bW$ as above, we have 
\[
\F(\bW_{\gamma \delta}) = \F(\bW^{\gamma \delta}) = (-1)^p
\]
where $p$ is the number of $x$'s among $\gamma$ and $\delta$
\end{lemma}
\noindent \Pf This comes from taking the decomposition of the normal, $\{w_i\}$ and tangent frames, $\{\va, v_x\}$ for $NY$, $TY$, as given in section \S \ref{MetricTY} and section \S \ref{NormalFrame}, and then noting that parity is preserved under inversion. \qed  \nl \nl
As a result of this, we define
\[
\bGamma_{\sigma \mu}^{\omega} := \bW^{\omega \nu} \bGamma_{\sigma \mu \omega}
\]
and conclude
\begin{corollary} \label{YChristoffelParityUpper}
\[
\F(\bGamma_{\sigma \mu}^{\omega}) = (-1)^p
\]
where $p$ the number of $x$'s among the indices $\sigma$, $\mu$, $\omega$
\end{corollary}

\subsection{Proof of theorem \ref{SFFParityMinimal}}
We prove theorem \ref{SFFParityMinimal}
\begin{theorem*} 
Suppose that $Y^m \subseteq M^{n+1}$ minimal with $\overline{h} = \overline{g} \Big|_{TY}$ even up to order $x^m$. Let $p \in Y$ and $\bA: \text{Sym}^2(TY) \to N(Y)$ denote the second fundamental form, and let $\{w_i(s,x)\}$ be the frame for the normal bundle described above. Define
\[
T_{\gamma \delta i; \alpha_1 \cdots \alpha_p} := \overline{g}\left((\bnabla_{v_{\alpha_1}} \cdots \bnabla_{v_{\alpha_p}} \bA)(v_{\gamma}, v_{\delta}), w_i \right)
\]
where $\alpha_i$ can take on any of the indices $\{s_a, x\}$. Let $q$ denote the number of ``$x$"s among the indices $\{\gamma, \delta, \alpha_1, \dots, \alpha_p\}$, then we have 
\[
\forall i, \qquad \F(T_{\gamma \delta i; \alpha_1 \cdots \alpha_p}) = (-1)^q 
\]
\end{theorem*}
\noindent\Pf The base case is an application lemma \ref{YChristoffelParity} as 
\begin{align*}
\bA_{\alpha \beta i} &= \tbGamma_{\alpha \beta i} \\
\implies \F(\bA_{\alpha \beta i}) &= \F(\tbGamma_{\alpha \beta i}) = (-1)^{\sigma(\alpha) + \sigma(\beta)}
\end{align*}
Now from here, we prove the theorem by induction for 
\[
g((\bnabla_{v_{\alpha_1}} \cdots \bnabla_{v_{\alpha_n}} \bA)(v_{\gamma}, v_{\delta}), w_i)
\]
where $\{v_{\delta}, v_{\gamma}, v_{\alpha_i} \} \in \{\va, v_x\}$. Assume the parity statement holds for $n - 1$. We compute
\begin{align*}
(\bnabla_{v_{\alpha_1}} \cdots \bnabla_{v_{\alpha_n}} \bA)(v_{\gamma}, v_{\delta}) &= \bnabla_{v_{\alpha_1}}^{\perp} [(\bnabla_{v_{\alpha_2}} \cdots \bnabla_{v_{\alpha_n}} \bA)(v_{\gamma}, v_{\delta})]  \\
& - (\bnabla_{v_{\alpha_2}} \cdots \bnabla_{v_{\alpha_n}} \bA)(\bnabla_{v_{\alpha_1}}^{\|} v_{\gamma}, v_{\delta})\\
& - (\bnabla_{v_{\alpha_2}} \cdots \bnabla_{v_{\alpha_n}} \bA)(v_{\gamma}, \bnabla_{v_{\alpha_1}}^{\|} v_{\delta}) \\
& = I_1 + I_2 + I_3
\end{align*}
here $\bnabla^{\|} = \bnabla^Y$ is the connection on $TY$ and $\bnabla^{\perp}$ is the connection on $NY$ (both using $\bh$). Let $p_{n-1}$ denote the number of $x$'s among the indices $\{\alpha_2, \dots, \alpha_n, \gamma, \delta\}$. For any index, $\omega$, recall the $\sigma(\omega)$ notation \ref{SigmaNotation}. We have
\begin{align*}
\F(g(I_1, w_i)) &= \F\left( \bg(\bnabla_{v_{\alpha_1}}^{\perp} [(\bnabla_{v_{\alpha_2}} \cdots \bnabla_{v_{\alpha_n}} \bA)(v_{\gamma}, v_{\delta})], w_i) \right) \\
&= \F \left( v_{\alpha_1} g((\bnabla_{v_{\alpha_2}} \cdots \bnabla_{v_{\alpha_n}} \bA)(v_{\gamma}, v_{\delta}), w_i) - g((\bnabla_{v_{\alpha_2}} \cdots \bnabla_{v_{\alpha_n}} \bA)(v_{\gamma}, v_{\delta}), \bnabla_{v_{\alpha_1}}^{\perp} w_i) \right) 
\end{align*}
By the inductive hypothesis, we have 
\[
\F \left( v_{\alpha_1} g((\bnabla_{v_{\alpha_2}} \cdots \bnabla_{v_{\alpha_n}} \bA)(v_{\gamma}, v_{\delta}), w_i) \right) = (-1)^{\sigma(\alpha_1) + p_{n-1}}
\]
And similarly
\[
\bnabla_{v_{\alpha_1}}^{\perp} w_i = \tbGamma_{\alpha_1 i}^k w_k 
\]
so that
\begin{align*}
\F(g((\bnabla_{v_{\alpha_2}} \cdots \bnabla_{v_{\alpha_n}} \bA)(v_{\gamma}, v_{\delta}), \bnabla_{v_{\alpha_1}}^{\perp} w_i)) &= \F(\tilde{\bGamma}_{\alpha_1 i}^k) \F(g((\bnabla_{v_{\alpha_2}} \cdots \bnabla_{v_{\alpha_n}} \bA)(v_{\gamma}, v_{\delta}), w_k)) \\
&= (-1)^{\sigma(\alpha_1) + p_{n-1}}
\end{align*}
again by the inductive hypothesis and lemma \ref{YChristoffelParity}. $I_2$ and $I_3$ proceed analogously, completing the induction. \qed
\section{Asymptotics for the variational vector field}
\label{AsymptoticVariational}
Having established an expansion for $u$, we want to show the analogous expansion for our variational vector fields $\dot{S}:= F_*(\partial_t) \Big|_{t = 0}$ and $\ddot{S} = \nabla_{F_*(\partial_t)} F_*(\partial_t) \Big|_{t = 0}$. We first need to fix a frame for the normal bundle.
\begin{lemma}
For any $p \in \gamma$ and a neighborhood $N(p) \subseteq M$, there exists a frame $\{w_1, \dots, w_{n+1-m}\}$ for $N(Y)$ which is orthonormal on $(\gamma, \bg\Big|_{\gamma})$ such that 
%
\begin{empheq}[box=\fbox]{gather*}
\overline{g}(w_i, \partial_x) = \text{$O(x)$ and odd up to order $m - 1$} \\
\overline{g}(w_i, \partial_{s_a}) = \text{$O(x^2)$ and even up to order $m + 2$} \\
\overline{g}(w_i, \partial_{r_j}) = \text{$\delta_{ij} + O(x^2)$ and even up to order $m$}
\end{empheq}
%
\end{lemma}
\noindent Alternatively, we phrase this as 
\begin{align*}
w_i &= c^{a}_i(s,x) \partial_{s_a} + c^{x}_i(s,x) \partial_x + c^{j}_i(s,x) \partial_{z_j} \\
\F(x^{-2} c^{a}_i) &= 1 \\
\F(c^{x}_i) &= -1 \\
\F(c^{j}_i) &= 1
\end{align*}
This is done by taking a normal frame for $\gamma \subseteq \partial M$, translating it to the interior so that the frame is constant in $x$, and projecting onto $TY^{\perp}$. See \eqref{NormalFrame} in the appendix. With this frame, we prove
\begin{theorem} \label{Variation}
Consider $\dot{S} = \dot{\phi}^i(s,x) w_i(s,x)$ and $\ddot{S} = \ddot{\phi}^i(s,x) w_i(s,x)$ be the first and second variational vector fields for a family of minimal submanifolds $\{Y_t\}_{t \geq 0}$ with $Y$ as in \S \ref{prelims}. Then
\begin{align*}
\F(\dot{\phi}^i) &= 1 \\
\F(\ddot{\phi}^i) &= 1
\end{align*}
Moreover, when $m = n$ even, there are no $x^n \log(x)$ or $x^{n+1} \log(x)$ terms.
\end{theorem} 
\noindent The theorem says that in a good ($x$-dependent!) frame for the normal bundle, we have a polyhomogeneous expansion to all orders which is even up to order $m$ ($m+1$) for $m$ even (odd). The idea is that $\dot{S}$ satisfies a homogeneous Jacobi equation since $Y_0$ is minimal, and $\ddot{S}$ satisfies an inhomogeneous Jacobi Equation since $\{Y_t\}$ is a variation through minimal submanifolds. We leverage these equations to deduce a polyhomogeneous expansion of $\dot{\phi}^i$ and $\ddot{\phi}^i$ by doing the analogous PDE analysis for the minimal surfaces system as in section \S \ref{Graphical}. 

\subsection{Jacobi Operator in full codimension}
By definition, $Y_t: = \exp_Y(S_t)$ and $\dot{S} = \partial_t S_t \Big|_{t = 0}$. Given that $\{Y_t\}$ is a family of minimal submanifolds, $\dot{S}$ lies in the kernel of the Jacobi operator
\[
J(X) = \Delta_Y^{\perp} X + \tilde{A}(X) + \text{Tr}[R_M(\cdot, X) \cdot]
\]
Here $\Delta_Y^{\perp}$ denotes the laplacian on the normal bundle, $\tilde{A}(X)$ denotes the Simons' operator on $Y$, and $\text{Tr}[R_M(\cdot, X) \cdot]$ is a trace of the \textit{ambient} Riemann curvature tensor over $TY$. As we showed in section \S \ref{ExpandMC}, we have
\begin{proposition}
For $Y^m \subseteq M^{n+1}$ as in our setup, the Jacobi operator decomposes as
\[
J(\phi^i w_i) = [(x \partial_x)^2 - (m-1) (x \partial_x) - m](\phi^i) \; w_i + R^i(\{\phi^j\}) w_i
\]
where 
\[
R:  x^{\delta} C^{k+2, \alpha}_0(Y) \to x^{\delta + 2} C^{k, \alpha}_0(Y)
\]
is an error term
\end{proposition}
\noindent In particular, if we expand $R$
\[
R = \sum_{p, \beta} r_{p, \beta}(s,x) (x \partial_x)^p (x \partial_s)^{\beta}
\]
for $\beta$ a multi-index, then $r_{p, \beta } = O(x^2)$ and $\F(r_{p,\beta}) = 1$. Because we have $\dot{\phi}^i = O(1)$, we see that the same PDE analysis and iteration argument as in section \S \ref{IterationArgument} gives
\[
\F(\dot{\phi}^i) = 1
\]
%
as desired. \qed \nl \nl
\rmk \; Note that if we choose to expand $\dot{S}$ in powers of $x$, we lose parity
\begin{align*}
\dot{S} &= \dot{S}_0(s) + x \dot{S}_1(s) + x^2 \dot{S}_2(s) + \dots \\
\dot{S}_i &\in T_p M, \qquad p \in Y
\end{align*}
i.e. both even and odd terms appear! This is because $N(Y)$ ``tilts" with $x$ so a priori, we have no parity of $\dot{S}$ in powers of $x$ with $x$-independent vectors, $\{\dot{S}_k(s)\}$. 
%
%

\subsection{Regularity and Parity of $\ddot{S}$} \label{SecondVariationRegAndPar}
\begin{proposition} \label{SecondVariationOfMC}
Let $\{Y_t\} \subseteq M^{n+1}$ be a family of minimal of $m$-dimensional minimal submanifolds. Let $Y = Y_0$ and $\overline{h}$ denote a compactified metric on $Y$. Then for
\[
Y_t = \{ \exp_{\overline{h}, p}(S_t(p) \bnu(p)) \; | \; p \in Y\}
\]
The second variation of mean curvature is given by
\[
\frac{d^2}{dt^2} H_t = J_Y^{\perp}(\ddot{S}) + Q^{\perp}(\dot{S})
\]
where $Q^{\perp}$ is a quadratic differential functional in $\dot{S}$ and $\F(Q^{\perp}(\dot{S})) = 1$
\end{proposition}
\noindent The details and the verification that 
\begin{align*}
Q^{\perp}(\dot{S}) &= Q^i(s,x) w_i \\
\F(Q^i) &= 1
\end{align*}
are shown in the appendix \S \ref{SecondVarDetails}. By the same work with $\dot{S}$, this immediately gives
\begin{theorem} \label{SecondDerivParity}
Let $\{Y_t\} \subseteq M^{n+1}$ be a family of minimal of $m$-dimensional minimal submanifolds. Let $Y = Y_0$ and $\overline{h}$ denote a compactified metric on $Y$. Then for
\[
Y_t = \{ \exp_{\overline{h}, p}(S_t(p)) \; | \; p \in Y\}
\]
with $S_t(p) \in NY$ for all $t$ and $\{w_i\}$ the normal basis described in \ref{NormalFrame}, we have 
\[
\ddot{S} = \frac{d^2}{dt^2} \Big|_{t = 0} S_t = \ddot{S}^i w_i
\]
Moreover
\[
\F(\ddot{S}^i) = 1
\]
and when $m = n$ even, there are no $x^n \log(x)$ or $x^{n+1} \log(x)$ terms.
\end{theorem}
\noindent \rmk \; In codimension $1$, we write 
\begin{align} \label{codim1Expansions}
\dot{S} &= \dot{\phi}(s,x) \bnu(s,x) \\ \nonumber
\ddot{S} &= \ddot{\phi}(s,x) \bnu(s,x)
\end{align}
where $\bg(\bnu, \bnu) \equiv 1$. In parallel with remark \ref{LocalTermsuExpansion}, we note that by the same proof mechanics in \S \ref{IterationArgument}, we have 
\begin{Remark} \label{LocalTermsPhiExpansion}
For $Y^n \subseteq M^{n+1}$ with $n$ odd, we have the expansions of 
\begin{align*}
\dot{\phi}(s,x) &= \dot{\phi}_0(s) + \dot{\phi}_2(s) x^2 + \dots + \dot{\phi}_{n-1}(s) x^{n-1} + \dot{\Phi}(s) x^{n+1} \log(x) + \dot{\phi}_{n+1}(s) x^{n+1} + O(x^{n+2} \log(x)) \\
\ddot{\phi}(s,x) &= \ddot{\phi}_0(s) + \ddot{\phi}_2(s) x^2 + \dots + \ddot{\phi}_{n-1}(s) x^{n-1} + \ddot{\Phi}(s) x^{n+1} \log(x) + \ddot{\phi}_{n+1}(s) x^{n+1} + O(x^{n+2} \log(x))
\end{align*}
and $\{\dot{\phi}_{2k}, \ddot{\phi}_{2k}\}_{k = 0}^{(n-1)/2}$, $\dot{\Phi}, \ddot{\Phi}$ are local, i.e. polynomials of derivatives of $\dot{\phi}_0(s)$, $\ddot{\phi}_0(s)$ (respectively).
\end{Remark}
\noindent Having shown that $\{u^i\}$, $\{\dot{\phi}^j\}$, and $\{\ddot{\phi}^k\}$ have polyhomogeneous expansions, we want to the variations of renormalized volume. To do this, we first review the mechanics of finite part evaluation
\section{Mechanics of Finite Part Evaluation} \label{mechanics}
When computing variations of renormalized volume, we encounter integrals of the form 
\[
I(z) = \int_Y z^p b(x,s) x^{z-j} dA_Y
\]
for $b$ having a polyhomogeneous expansion in $x$ (after pulling back to $\Gamma$) and $i \geq 0$, $j \in \{0,1,2\}$. We write 
\[
I(z) = \left(\int_{Y \cap \{x < \delta\}} + \int_{Y \cap \{x \geq \delta\}}\right) z^p b(x,s) x^{z-j} dA_Y = I_1(z, \delta) + z^p I_2(z, \delta)
\]
for some $1 \gg \delta > 0$, where we've pulled out the factor of $z^p$ in the $\{x \geq \delta\}$. As before, $I_2(z)$ is holomorphic because the integral is over $x \geq \delta$. In particular
\[
\FPz z^p I_2(z, \delta) = \begin{cases}
I_2(0, \delta) & p = 0\\
0 & p \geq 1
\end{cases}
\]
We further assume the following expansions (after pulling back to Fermi coordinates)
\begin{align} \nonumber
dA_Y &= \frac{\sqrt{\det \overline{h}}}{x^m} dx \wedge dA_{\gamma} \\ \nonumber
\sqrt{\det \overline{h}} &= \sum_{k = 0}^{m+2} \overline{q}_j(s) x^k + \tilde{Q}(s) x^m \log(x) + \overline{Q}(s) x^{m+1} \log(x) + O(x^{m+2} \log(x)) \\ \label{bExpansion}
b(x,s) &= \sum_{k = 0}^{m+2} b_j(s) x^k + \tilde{B}(s)x^m \log(x) + B(s) x^{m+1} \log(x) + O(x^{m+2} \log(x))
\end{align}
i.e. if a $x^d \log(x)^q$ term manifests, it can occur only when $d \geq m$. This accounts for both even and odd expansion of $u(s,x)$ and $\{\overline{h}^{\alpha \beta}\}$ as in \S \ref{Graphical}. $I_1$ expands as 
\begin{align*}
I_1(z, \delta) &= z^p \int_{0}^{\delta} \int_{\gamma} x^{z-m-j}\left[\sum_{\ell = 0}^{m+2} \sum_{k + j = \ell} b_{k} \overline{q}_j x^{\ell} \right] ds dx  \\
& + z^p \int_{0}^{\delta} \int_{\gamma} x^{z-m-j}\left[ (\overline{q}_0 \tilde{B} + \tilde{Q} b_0 ) x^m \log(x) + (b_0 \overline{Q} + b_1 \tilde{Q} + \overline{q}_0 B + \overline{q}_1 \tilde{B}) x^{m+1} \log(x)\right] ds dx  \\
&+ z^p \int_{0}^{\delta} \int_{\gamma} x^{z-m-j}\left[ (b_1 \overline{Q} + b_2 \tilde{Q} + \overline{q}_1 B + \tilde{B} \overline{q}_2)x^{m+2} \log(x) \right] ds dx \\
& + z^p\int_0^{\delta} \int_{\gamma} O(x \log(x)) ds dx
\end{align*}
Observe that
\[
\FPz z^p \int_{0}^{\delta} \int_{\gamma} O(x \log(x)) dx = \begin{cases}
C(\delta) & p = 0 \\
0 & p \geq 1
\end{cases}
\]
for some finite constant $C(\delta)$. It remains to compute
\begin{align} \label{FPExpansion}
\FPz I_1(z, \delta) &= \FPz z^p \sum_{k = 0}^{m+2} c_k \int_{0}^{\delta} x^{z+ k - m - j} dx  \\ \nonumber
& + \FPz z^p \sum_{k = 0}^2 c_{m + k}^* \int_0^{\delta} x^{z + k - m - j} \log(x) dx
\end{align}
for 
\begin{align*}
c_k &= \int_{\gamma} \left[\sum_{\ell + j = k} b_{\ell} \overline{q}_j  \right] dV_{\gamma}, \qquad 0 \leq k \leq m + 2 \\
c_m^* &= \int_{\gamma} \left[\overline{q}_0 \tilde{B} + \tilde{Q} b_0 \right] dA_{\gamma} \\
c_{m+1}^* &= \int_{\gamma} \left[ b_0 \overline{Q} + b_1 \tilde{Q} + \overline{q}_0 B + \overline{q}_1 \tilde{B}\right] dA_{\gamma} \\
c_{m+2}^* &= \int_{\gamma} \left[ b_1 \overline{Q} + b_2 \tilde{Q} + \overline{q}_1 B + \tilde{B} \overline{q}_2\right] dA_{\gamma}
\end{align*}
Integrating,
\[
I_1(z, \delta) = z^p \left(c_{m+j-1}\frac{\delta^z}{z} + c^*_{m+j-1}\frac{\delta^z ((z \log(\delta) - 1)}{z^2}+ F(\delta, z)\right) 
\]
where $F(\delta, z)$ is a holomorphic near $0$ and computable in terms of the remaining $c_k$ coefficients through equation \eqref{FPExpansion}. In particular
\begin{align*}
\FPz z^p F(\delta, z) &= \begin{cases}
F(\delta, 0) & p = 0 \\
0 & p \geq 1
\end{cases} \\
\FPz z^p \frac{\delta^z}{z} &= 
\begin{cases}
\log(\delta) & p = 0\\
1 & p = 1 \\
0 & p > 1
\end{cases} \\
\FPz z^p \frac{\delta^z (z \log(\delta) - 1)}{z^2} &=
\begin{cases}
\frac{\log(\delta)^2}{2} & p = 0 \\
0 & p = 1 \\
-1 & p = 2 \\
0 & p \geq 3
\end{cases}
\end{align*}
Note that a $x^d \log(x)$ term in the expansion of $I_1$ leads to higher order poles. We summarize this work as
\begin{lemma} \label{FiniteEvaluation}
Consider integrals of the form
\[
I(z) = \int_Y z^p b(x,s) x^{z-j} dA_Y 
\]
for $b$ and $dA_Y$ having polyhomogeneous expansions in $x$ and $p \geq 0$, $j \in \{0,1,2\}$. Moreover, assume that $x^d \log(x)^q$ terms only manifest when $d \geq m$ and $q = 1$, or $d \geq m + 3$. Then we have that 

\begin{empheq}[box=\fbox]{equation*}
\FPz I(z) = \begin{cases} 
	C(\delta) + F(\delta, 0) + I_2(0, \delta) + c_{m+j-1}\log(\delta) + c^*_{m+j-1}\frac{\log(\delta)^2}{2} & p = 0 \\[1ex]
	c_{m+j-1} & p = 1 \\[1ex]
	-c^*_{m+j-1} & p = 2 \\[1ex]
	0 & p \geq 3
\end{cases}
\end{empheq}
for the coefficients $\{c_{k}\}$ listed above
\end{lemma}
\noindent \rmk:
\begin{itemize}
\item This calculation illustrates the following key point: \textbf{when $p \geq 1$, the finite part can be expressed as an integral over the boundary}. We will refer to this process from here on as \textbf{localization}. 
\item In the future we write 
\[
B(s) = [b(x,s)]^{\log, k}
\] 
to indicate the $x^{k} \log(x)$ term in the expansion of any function, $b(x,s)$, as in equation \eqref{bExpansion}.
\item Taking $b(x,s) = 1$ and $p = 0$ demonstrates how to compute the renormalized volume of $Y$ via Riesz regularization

\item While the result for $p = 0$ seems to depend on $\delta$, one can show that by changing $\delta \to \delta'$ and keeping track of boudary terms from the intermediate integral $\int_{x = \delta}^{\delta'}$, the result is independent of $\delta$. This is shown for $b(x,s) = 1$ in \S \ref{equivOfDef} but holds for any $b(x,s)$ of the form \eqref{bExpansion}.
\end{itemize}

\section{Renormalized Volume for $Y \subseteq M^{n+1}$}
\label{SpecialBDFonY}
Let $x_Y$ be a special bdf on $(Y, h)$ considered as its own asymptotically hyperbolic manifold with metric even to high order. In this section, we prove the following:
\begin{theorem}
\label{EquivofRVonY}
Let $Y^m \subseteq M^{n+1}$ minimal, satisfying the conditions \S \ref{Graphical} and $x_Y$ a special bdf on $Y$. Let $x$ a special bdf on $M$, inducing the same conformal infinity on $Y$. We have that 
\begin{empheq}[box=\fbox]{equation*}
\FPz \int_Y x_Y^z dA_Y = \begin{cases} \FPz \int_Y x^z dA_Y & \text{ m even }\\[2ex] 
	\int_{\gamma} p(s) dA_{\gamma}(s) + \FPz \int_Y x^z dA_Y & \text{ m odd }
\end{cases}
\end{empheq}
where $p(s) = F(u_2(s))$ is a computable function of $u_2(s) = \frac{1}{2(m-1)} H_{\gamma}(s)$ and higher derivatives. 
\end{theorem}
\noindent This theorem says that for even dimensional manifolds, we can use either a special bdf on $Y$, which is labeled as $x_Y$, or the almost special bdf, $x$, on $M^{n+1}$. Recall that a special bdf, $x$, satisfies
\[
||dx||_{\overline{g}}^2 = \overline{g}^{xx} = 1
\]
where $\overline{g} = x^2 g$. We want to find a special bdf, $x_Y$, for $Y$, such that 
\[
||d \log(x_Y)||_{h}^2 = h^{\alpha \beta} x_Y^{-2} dx_Y(v_{\alpha}) dx_Y(v_{\beta}) = 1
\]
where $\{v_{\alpha}\}$ is the pushforward of the coordinate basis for $TY$ defined in \ref{MetricTY} and $h = g \Big|_Y$. As in \cite{albin2009renormalizing}, \cite{graham1991einstein} and \cite{graham1999conformal}, we begin with a bdf on $Y$ written as
\[
x_Y(x) = x e^{\omega(s,x)}
\]
where
\[
\omega(s,x) = \sum_{k = 0}^{\infty} \omega_k(s) x^k
\]
such an expansion was shown in \cite{graham1991einstein}. We now enforce $1 = ||d \log(x_Y)||_{h}^2$:
\begin{align} \label{OmegaEquation}
1 &= h^{\alpha \beta} e^{-2\omega} x^{-2} [e^{\omega} dx(v_{\alpha}) + x e^{\omega} d\omega(v_{\alpha})] [e^{\omega} dx(v_{\beta}) + x e^{\omega} d\omega(v_{\beta})] \\ \nonumber
&= \overline{h}^{\alpha \beta} [dx(v_{\alpha}) + x d\omega(v_{\alpha})] [dx(v_{\beta}) + x d\omega(v_{\beta})]
\end{align}
As in \cite{graham1991einstein}, the above equation shows that $\omega_1 = 0$, and in general that $\omega$ has an even expansion to high order. When $m$ is even, the first non-trivial odd coefficient occurs at $x^{m+1}$, with potentially an $x^m \log(x)$ in the codimension $1$ case. When $m$ is odd, there may be $x^{m} \log(x)$ and $x^{m+1} \log(x)$ terms. In both cases, the first odd order term in \eqref{OmegaEquation} comes from the first odd order terms in the expansion of $\{\overline{h}^{\alpha \beta}\}$. To summarize
\begin{lemma}
\label{OmegaExpansion}
Let $Y^m \subseteq M^{n+1}$ be a minimal submanifold. There exists a bdf $x_Y: Y \to \R^+$ such that 
\[
x_Y(s,x) = x e^{\omega(s,x)}
\]
with 
\[
\F(\omega) = 1
\]
\end{lemma}
\noindent We now prove theorem \ref{EquivofRVonY}
\subsection{Equivalence of Renormalized Volume of $Y^n \subseteq M^{n+1}$, $m$ even}
\noindent Let $x$ be a special bdf on $M$ and $x_Y$ a special bdf for $Y$. We compute the following difference
%
\begin{align*}
\int_Y (x^z - x_Y^z) dA_Y &= \FPz \int_Y x^z(1 - e^{\omega(s,x) z}) dA_Y = \int_Y (z \omega + O(z^2)) dA_Y \\
&= \FPz \int_Y z \omega dA_Y = \int_{\gamma} [\omega(s,x) \sqrt{\det \bh}]^{(m-1)} \\
&= \int_{\gamma} 0 = 0
\end{align*}
having used that 
\[
\F(\sqrt{\det \bh}) = 1, \quad \F(\omega) = 1, \quad \implies \quad \F\left(\omega(s,x) \sqrt{\det \bh}\right) = 1
\]
so $[\omega(s,x) \sqrt{\det \bh}]^{(m-1)} = 0$. Note that the first and second variation of renormalized volume can also be computed with $x$ instead of $x_Y$. The proof uses the same techniques as showing that these variations are independent of the choice of conformal infinity, which is done in \S \ref{ConformalInvariance}. \qed

\subsection{Anomly for Renormalized Volume of $Y^n \subseteq M^{n+1}$, $m$ odd}
When $m$ is odd, the two definitions of renormalized volume using $x$ and $x_Y$ are not equal. This is discussed in \cite{graham1999conformal} among other sources, but we compute the anomaly here using Riesz Reguarlization:
\[
\FPz \int_Y x^z(1 - e^{\omega(s,x) z}) dA_Y = \int_{\gamma} [\omega(s,x) \sqrt{\det \bh}]^{(m-1)} = \int_{\gamma} \sum_{k = 0}^{(m-1)/2} \omega_{2k}(s) \overline{q}_{(m-1) - 2k}(s)
\]
Because $m$ is odd, this sum may be non-zero and the renormalized volume depends on the choice of bdf. We note that for $2k < m + 1$ the coefficients of $\{\omega_{2k}(s)\}$ and $\overline{q}_{2k}(s)$ are determined by $u_2(s) = \frac{1}{2(m-1)} H_{\gamma}(s)$ via the iterative procedure used to show their existence (see \S \ref{IterationArgument}). As a result,
\[
\V(Y) = \FPz \int_Y x^z dA_Y + \int_{\gamma} p(s) dA_{\gamma}
\]
where $p(s)$ is a function determined by $u_2(s) = \frac{1}{2(m-1)} H_{\gamma}(s)$ and its derivatives. \qed 

\section{Variational Formulae}
\label{MainResult}
We derive formulae for the first and second variations of minimal submanifolds $Y^m \subseteq M^{n+1}$. Following \cite{alexakis2010renormalized}, let $\{Y_t\}$ be a one-parameter family of minimal submanifolds and assume each $Y_t \cap \{x < \delta\}$ is embedded for some small $\delta > 0$. From equation \eqref{prelims} and section \S \ref{Graphical}, we can write for $m$ even
\begin{align*}
u(s,x) &= u^i \bN_i \\
u^i(s,x) &= u_2^i(s) x^2 + u_4^i(s) x^4 + \dots + u_m^i(s) x^m + u_{m + 1}^i(s) x^{m + 1} +  O(x^{m+2}) \\
\dot{S} &= \dot{S}^i w_i \\
\dot{S}^i(s,x) &= \dot{S}^i_0(s) + \dot{S}^i_2(s) x^2 + \dots + \dot{S}^i_m(s) x^m + \dot{S}^i_{m+1} x^{m+1} + O(x^{m+2}) 
\end{align*}
and analogously for $m$ odd with some $x^{m+1} \log(x)$ terms potentially.
%
%
\begin{theorem}
\label{MainThm}
Let $\{Y_t\} \subseteq M^{n+1}$ be a one-parameter family of $m$-dimensional minimal submanifolds for $m < n + 1$ and with $Y := Y_0$. Further suppose that for some $\delta > 0$, for all $ t > 0$ sufficiently small, $Y_t \cap \{x < \delta\}$ is embedded in $\{x < \delta\}$, and that 
\[
Y_t \cap \mathcal{U} = \text{Im}\left(\exp_Y(S_t) \Big|_{x < \delta}\right)
\]
for $S_t \in N(Y)$. If $0$ is not in the $L^2$ spectrum of the Jacobi operator, $J_Y$, and $\dot{S} = \frac{d}{dt} S_t \Big|_{t = 0}$ is a bounded Jacobi field (w.r.t. $\overline{g}$), then the first variation of renormalized volume is given by
\[
D \mathcal{V}\Big|_Y(\dot{S}) = \int_{\gamma}  [dx(\dot{S}) \sqrt{\det \overline{h}} ]^{(m)} dA_\gamma
\]
where $dA_Y = \frac{\sqrt{\det \overline{h}}}{x^m} (d A_\gamma \wedge dx)$. Furthermore,
\begin{equation} \label{SecondVarEq}
D^2 \mathcal{V} \Big|_{Y} (\dot{S}, \ddot{S}) = \begin{cases} \int_{\gamma} \Big(-[dx(\dot{S})^2 \sqrt{\det \overline{h}}]^{(m+1)} + [dx(\ddot{S}) \sqrt{\det \bh}]^m  & \;\; \text{ $m$ even }\\[1ex] 
	- \frac{1}{2}\left[ ||\dot{S}||^2 \cdot ||\nabla x||^2 \sqrt{\det \overline{h}} \right]^{(m + 1)} + \frac{1}{2}\left[ ||\dot{S}||^2 \Delta x \sqrt{\det \overline{h}} \right]^{(m)} \Big) \; d A_\gamma &  \\[1ex] 
	& \\
	\int_{\gamma} \Big (- [dx(\dot{S})^2 \sqrt{\det \overline{h}}]^{(m + 1, \log)} + [dx(\ddot{S}) \sqrt{\det \bh}]^m & \; \; \text{$m$ odd} \\[1ex] 
	- [dx(\dot{S})^2 \sqrt{\det \overline{h}}]^{(m + 1)} + \frac{1}{2} (||\dot{S}||^2 \Delta x \sqrt{\det \overline{h}})^{(m)}  & \\[1ex] 
	- \frac{1}{2}(||\dot{S}||^2 ||\nabla x||^2 \sqrt{\det \bh})^{(m+1)} - (||\dot{S}||^2 ||\nabla x||^2 \sqrt{\det \bh})^{(m + 1, \log)} \Big) dA_{\gamma} & 
\end{cases}
\end{equation}
\end{theorem}
\noindent \rmk 

\begin{itemize}
\item These formulae show that variations of renormalized volume depend only on the following geometric quantities: the volume form, the special bdf, $x$, and the variational vector fields.
\item The condition of $0 \not \in \sigma(J_Y)$ guarantees that the moduli space of smooth minimal submanifolds with smooth boundary curves is a Banach space. The proof is analogous to the one in \cite{alexakis2010renormalized}, assuming that $Y$ is embedded in a neighborhood of the boundary.

\item The first variation formula holds as long as $Y = Y_0$ is minimal, and the remaining $\{Y_t\}$ have the same embedding and asymptotic expansion properties, i.e. they are not required to be minimal, as long as we have parity results for $\dot{S}$. The second variation formula requires minimality
%
\end{itemize}
\begin{corollary}[Codimension $1$] 
For $Y^{n} \subseteq M^{n+1}$ with $n$ even, $\dot{S} = \dot{\phi}\bnu$ for $\bnu$ a unit normal to $Y$ w.r.t $\bg$, the formulae above become
\begin{empheq}[box=\fbox]{align} \label{codim1EvenVariation}
D A\Big|_Y(\dot{\phi}) &=  -(n+1) \int_{\gamma} \dot{\phi}_0(s) u_{n+1}(s) dA_\gamma(s)\\ \nonumber
D^2 A \Big|_Y(\dot{\phi}) &= \int_{\gamma} -(n+1) \ddot{\phi}_0 u_{n+1} + (1 - n) \dot{\phi}_0(s) \dot{\phi}_{n+1}(s) \\
& \quad - 8n \dot{\phi}_0(s)^2 u_2 u_{n+1}(s) dA_\gamma(s) \nonumber
\end{empheq}
When $n$ is odd, we have
\begin{empheq}[box=\fbox]{align} \label{codim1OddVariation}
D A\Big|_Y(\dot{\phi}) &=  \int_{\gamma} \Big[ -(n+1) \dot{\phi}_0(s) u_{n+1}(s) + F_1(\dot{\phi}_0, u_2) \Big] dA_\gamma(s)\\ \nonumber
D^2 A \Big|_Y(\dot{\phi}) &= \int_{\gamma} -(n+1) \ddot{\phi}_0 u_{n+1} + (1 - n) \dot{\phi}_0(s) \dot{\phi}_{n+1}(s)  \\ \nonumber
& + \dot{\phi}_0(s)^2 \left[ -8n u_2 u_{n+1}(s) + \text{Tr}_{T\gamma}(k_{n+1, 0})\right]  \\
&- \dot{\phi}_0(s) \left[4(n+2) \dot{\phi}_0(s) u_2(s) U(s) + \dot{\Phi}(s) \right]  + F_2(\ddot{\phi}_0, \dot{\phi}_0, u_2) \; dA_\gamma(s) \nonumber
\end{empheq}
\end{corollary}
\noindent \rmk \;As we'll see in the proof, $F_1$ and $F_2$ are actually polynomials in the coefficients $\{u_2, \dots, u_{n}\}$, $\{\dot{\phi}_0, \dots, \dot{\phi}_n\}$, and $\{\ddot{\phi}_0, \dots, \ddot{\phi}_n\}$. As shown in \ref{IterationArgument}, these coefficients are determined by the derivatives of $u_2(s)$, $\dot{\phi}_0(s)$, and $\ddot{\phi}_0(s)$, respectively. Thus, $F_1$ and $F_2$ are differential operators that only depend on $\gamma$ (which determines $u_2$), $\dot{\phi}_0$, and $\ddot{\phi}_0$, the ``Dirichlet data" of $Y$, $\dot{S}$, and $\ddot{S}$. \nl \nl
Specializing to the case of $m = 2$ and $n + 1 = 3$ as in \cite{alexakis2010renormalized} we have
\begin{corollary}
For the set up as above with $Y^2 \subseteq \H^3$, we have
\begin{empheq}[box=\fbox]{align*}
D A\Big|_Y(\dot{\phi}) &= -3 \int_{\gamma} \dot{\phi}_0 u_3 dA_{\gamma} \\
D^2 A \Big|_Y(\dot{\phi}) &= \int_{\gamma} \left(-3 \ddot{\phi}_0 u_{3} - \dot{\phi}_0 \dot{\phi}_3 - 16 u_2 u_3 \dot{\phi}_0^2\right) dA_{\gamma}
\end{empheq}
\end{corollary}
\noindent Note that the formula for $D^2 A\Big|_Y(\dot{\phi})$ is a slight correction to the formula in \cite{alexakis2010renormalized}.
\section{Proof of Variational Formulae}
\label{ProofMainResult}
\subsection{First Variation }
We prove the first part of theorem \ref{MainThm}. Recall our set up: For $Y^m_t \subseteq M^{n + 1}$ a family of minimal submanifolds with $\partial Y_t = \gamma_t$, we describe these via Fermi coordinates off of $Y$ with respect to $\overline{g}_Y$:
\[
F(t,p) := \exp_p(S_t(p))
\]
with $S_t \in N(Y)$. We will write $F_t(p) = F(t,p)$ when we want to emphasize that we're working over a fixed $p$. We compute
\begin{align*}
\mathcal{V}(Y) &= \FPz \int_Y x^z dA_Y  \\ 
D \mathcal{V}\Big|_{Y}(\dot{S}): &= \frac{d}{dt} \mathcal{V}(Y_t) \Big|_{t = 0} \\ 
&= \FPz \frac{d}{dt} \int_{Y_0} F_t^*(x^z) F_t^* (dA_t) = \FPz \int_{Y_0} z x^{z-1}(p) dx(\dot{S}(p)) dA_0 \\
& = \FPz \int_{Y_0} z x^{z-1} dx(\dot{S}) dA_0 \\
&= \boxed{\int_{\gamma} [ dx(\dot{S}) \sqrt{\det \overline{h}}]^{(m)}}
\end{align*}
where in the third line we use $\frac{d}{dt} F_t^* (dA_t)\Big|_{t = 0} = 0$ from the minimal surface condition, and 
\[
\frac{d}{dt} F_t^*(x^z) = z x^{z-1} dx(\dot{S})
\]
Note that because we're only taking an $m$th term, the above result holds for both $m$ even and $m$ odd. When $m$ is odd only $x^{m+1} \log(x)$ terms appear, which doesn't affect the $(m)$th order term. \qed
\subsection{Second variation}
The second variation in theorem \ref{MainThm} is derived using the same procedure
\[
D^2\mathcal{V}\Big|_{Y}(X) = \frac{d^2}{dt^2} \FPz \int_{Y_t} x^z dA_t = \FPz \frac{d^2}{dt^2} \int_{Y_0} x(F(t,p))^z F_t^*(d A_t)
\]
Differentiating under the integral, we get 
\[
D^2\mathcal{V}\Big|_{Y}(X) = \FPz \left[ \int_{Y_0} [z (z - 1) x^{z-2} \dot{x}^2 + z x^{z-1} \ddot{x}] dA + 2 \int_{Y_0} z x^{z-1} \dot{x} \frac{d}{dt} F_t^*(dA_t) \Big|_{t = 0} + \int_{Y_0} x^z \frac{d^2}{dt^2} F_t^*(dA_t)\Big|_{t = 0} \right]
\]
where $\dot{x}$ and $\ddot{x}$ are equal to
\[
\frac{d^i}{dt^i} x(F(t,p)), \qquad i = 1, 2
\]
Note that $ \frac{d}{dt} F_t^*(dA_t)\Big|_{t = 0}$ vanishes when $Y_0$ is minimal. Hence
\begin{align*}
D^2\mathcal{V}\Big|_{Y} & = \FPz \left[\int_{Y_0} \left[ z (z - 1) x^{z-2} \dot{x}^2 + z x^{z-1}\ddot{x}\right] dA_Y + \int_{Y_0} x^z \frac{d^2}{dt^2} F_t^*(dA_t) \Big|_{t = 0} \right] \\
& = \FPz [I_1(z) + I_2(z)]
\end{align*}

\subsubsection{$I_1$ Computation}
We compute the finite part of the first integral $I_1 = A_1 + B_1$.
\begin{align*}
\FPz A_1 &= \FPz\int_{Y} z (z-1) x^{z-2} \dot{x}^2 dA  \\
&= \FPz z^2 \int x^{z-2} \dot{x}^2 dA_Y - \FPz \int_Y z x^{z-2} \dot{x}^2 dA_Y \\
&= \int_{\gamma}- \left[(dx(\dot{S}))^2 \sqrt{\det(\overline{h})}\right]^{(m + 1, \log)} - \left[(dx(\dot{S}))^2 \sqrt{\det(\overline{h})}\right]^{(m + 1)} 
\end{align*}
having used that $\dot{x} = dx(\dot{S})$ along with the techniques in \eqref{mechanics}. For $B_1$, we write this as 
\begin{align*}
B_1 &= \int_{Y} z x^{z-1} \ddot{x} = \int_{Y} z x^{z-1} dx(\ddot{S}) \\
& = \int_{\gamma} [dx(\ddot{S}) \sqrt{\det(\bh)}]^{m} 
\end{align*}
Thus we have
\[
\boxed{\FPz I_1 = \FPz (A_1 + B_1) = \int_{\gamma}\left(-[(dx(\dot{S}))^2 \sqrt{\det(\overline{h})}]^{(m + 1, \log)} -[dx(\dot{S})^2 \sqrt{\det \overline{h}}]^{(m + 1)} + [dx(\ddot{S}) \sqrt{\det(\bh)}]^{m} \right) dA_{\gamma}}
\]
\subsubsection{$I_2$ Computation}
We compute 
\[
I_2 = \FPz \int_Y x^z \frac{d^2}{dt^2} F_t^*(dA_t) \Big|_{t = 0}
\]
We know from a variety of sources (e.g. \cite{colding2011course}) that for a geodesic variation
\[
\frac{d^2}{dt^2} F_t^*(dA_t) = \langle \nabla^{\perp} \dot{S}, \nabla^{\perp} \dot{S} \rangle - \Ric^{\perp}(\dot{S},\dot{S}) -  |\langle A(\cdot, \cdot), \dot{S} \rangle|^2 + \div_{Y} (\ddot{S})
\]
where $\nabla^{\perp}$ denotes the connection on the normal bundle, $A$ denotes the second fundamental form, $\Ric^{\perp}$ is the trace over $TY \subseteq TM$ of the ambient Riemann curvature applied to elements in $N(Y)$. We first integrate by parts on the divergence term and get 
\[
\int_Y x^z \div_Y(\ddot{S}) dA_Y = \int_{\gamma} x^z \langle \ddot{S}, \hat{n} \rangle dA_{\gamma} - \int_Y z x^{z-1}  \langle \nabla^Y x, \ddot{S} \rangle
\]
Again, the boundary term vanishes because we first assume $\text{Re}(z) \gg 0$. For the second term, $\nabla^Y x \in TY$ and as we show in \S \ref{SecondVarDetails}, $\ddot{S} \in NY$, so this term vanishes automatically. \nl \nl
We now handle the remaining terms 
\[
\langle \nabla^{\perp} \dot{S}, \nabla^{\perp} \dot{S} \rangle - \Ric^{\perp}(\dot{S},\dot{S}) -  |\langle A(\cdot, \cdot), \dot{S} \rangle|^2
\] 
This is the quadratic form for the corresponding Jacobi operator
\[
J_Y X = \Delta^{\perp} X + \Ric^{\perp}(X, \cdot) + \tilde{A}(X)
\]
where 
\[
\tilde{A}(X) = g((\nabla_{v_{\alpha}} v_{\beta})^N, X) h^{\alpha \gamma} h^{\beta \delta} (\nabla_{v_{\gamma}} v_{\delta})^N
\]
is the Simons operator. Because we consider a \textit{variation among minimal submanifolds}, we have $J_Y(\dot{S}) = 0$. In order to get the integrand in the form of the Jacobi operator, we integrate by parts and gain a boundary term which contributes to our second variation. Thus
\[
\int_Y x^z \left(\langle \nabla^{\perp} \dot{S}, \nabla^{\perp} \dot{S} \rangle - \Ric^{\perp}(\dot{S},\dot{S}) - |\langle A(\cdot, \cdot), \dot{S} \rangle|^2 \right) dA_y = A_1 + A_2 + A_3
\]
We now integrate by parts on the first term in our original expression for $I_2$ and get
\begin{align} \nonumber
\int_Y x^z \langle \nabla^{\perp} \dot{S}, \nabla^{\perp} \dot{S}\rangle & = - \int_Y z x^{z-1} \sum_{i = 1}^m \partial_{i}(x) \langle \dot{S}, \nabla_i \dot{S} \rangle - \int_Y x^z \sum_{i = 1}^m\langle \dot{S}, \nabla_i \nabla_i \dot{S} \rangle + \int_{\gamma} x^z \sum_{i = 1}^m \langle \dot{S}, \nabla_i \dot{S} \rangle \\ \label{IBPFirstLook}
&= - \int_Y z x^{z-1} \sum_{i = 1}^m \partial_{i}(x) \langle \dot{S}, \nabla_i \dot{S} \rangle - \int_Y x^z \langle \dot{S}, \Delta^{\perp} \dot{S} \rangle
\end{align}
where we sum over an orthonormal frame of $TY^m$, $\{e_1, \dots, e_m \}$ with respect to $g\Big|_{Y}$.  The last integral in the first line vanishes because $x \Big|_{\gamma} \equiv 0$ when $Re(z) \gg 0$. To see this, Let
\begin{equation} \label{BoundaryInt}
	\overline{I}(z) := \int_{\gamma} x^z \sum_{i = 1}^m \langle \dot{S}, \nabla_i \dot{S} \rangle
\end{equation}
and note that 
\[
\text{Re}(z) \gg 0 \implies  \overline{I}(z) \equiv 0
\]
Since $x \Big|_{\gamma} = 0$ and $\sum_{i = 1}^m \langle \dot{S}, \nabla_i \dot{S} \rangle = O(x^{-2})$ (i.e. has finite order blow up in $x$). Thus for the purposes of computing the meromorphic extension of $f(z) = \int_Y x^z \frac{d^2}{dt^2} F_t^*(dA_t) |_{t = 0}$ we can ignore $\overline{I}(z)$, i.e. we choose its meromorphic extension to be $0$ on all of $\C$. In particular, this means that $\FPz \overline{I}(z) = 0$. \nl \nl
\noindent The second integral in the second line of equation \eqref{IBPFirstLook} combines with $A_2$ and $A_3$ to yield $0$ because $\dot{S}$ is a Jacobi field. Thus
\[
I_2 = \FPz \int_Y x^z \frac{d^2}{dt^2} F_t^*(dA_t) = - \FPz \int_Y z x^{z-1} \sum_{i = 1}^m \partial_i(x) \langle \dot{S}, \nabla_i \dot{S} \rangle dA_Y 
\]
Integrating the first term by parts again, we get
\begin{align*}
- \int_Y z x^{z-1} \sum_{i = 1}^m \partial_i(x) \langle \dot{S}, \nabla_i \dot{S} \rangle dA_Y &= - \frac{1}{2} \int_Y z x^{z-1} \sum_{i = 1}^m \partial_i (x) \partial_i ||\dot{S}||^2 \\ 
&= \frac{1}{2} \int_Y z ||\dot{S}||^2 \partial_i\left(x^{z-1} \partial_i(x) \right) - \frac{1}{2} \int_{\gamma} z x^{z-1} \sum_i \partial_i(x) ||\dot{S}||^2
\end{align*}
Again, the boundary integral vanishes in the meromorphic extension, with the argument identical to the term \eqref{BoundaryInt}. We take the remaining integral and expand it as 
%
\begin{align*}
\frac{1}{2} \int_Y z ||\dot{S}||^2 \sum_{i}\partial_i\left(x^{z-1} \partial_i(x) \right) &= \frac{1}{2} \int_Y z ||\dot{S}||^2 \sum_i [(z-1) x^{z-2}  (\partial_i(x))^2 + x^{z-1} \partial_i^2(x) ] \\
&= \frac{1}{2}\int_Y z ||\dot{S}||^2 \left[ (z-1) x^{z-2} ||\nabla x||^2 + x^{z-1} \Delta x\right]
\end{align*}
Note that when we write $\Delta x$ and $\nabla x$, we consider $x$ as a function restricted to $Y$ and compute the laplacian and gradients with respect to bases on $TY$ with the complete metric $g$. Now we localize and get 
\begin{align*}
\FPz I_2(z) & = \FPz \int_Y x^z \frac{d^2}{dt^2} F_t^*(dA_t) \\
& = \frac{1}{2} \FPz \int_Y z(z-1) x^{z-2} ||\dot{S}||^2 ||\nabla x||^2 dV_{Y} + \frac{1}{2} \FPz \int_Y z x^{z-1} ||\dot{S}||^2 \Delta x dA_Y 
\end{align*}
%
so that 
\begin{empheq}[box=\fbox]{align*}
\FPz I_2 &=  - \frac{1}{2}\int_{\gamma} \left[ ||\dot{S}||^2 ||\nabla x||^2 \sqrt{ \det \bh} \right]^{m+1, \log} - \frac{1}{2}\int_{\gamma} \left[ ||\dot{S}||^2 \cdot ||\nabla x||^2 \sqrt{\det \overline{h}} \right]^{(m + 1)} dA_{\gamma} \\
& + \frac{1}{2}\int_{\gamma} \left[ ||\dot{S}||^2 \Delta x \sqrt{\det \overline{h}} \right]^{(m)} dA_{\gamma} 
\end{empheq}

\subsubsection{Putting it Together}
We add the two integrals and get for
\begin{empheq}[box=\fbox]{align*}
D^2 \mathcal{V} \Big|_{Y} (\dot{S}, \dot{S}) &= \int_{\gamma}- [dx(\dot{S})^2 \sqrt{\det \overline{h}}]^{(m + 1, \log)} + [dx(\ddot{S}) \sqrt{\det \overline{h}}]^{m} - [dx(\dot{S})^2 \sqrt{\det \overline{h}}]^{(m + 1)}\\
& + \frac{1}{2} \left[(||\dot{S}||^2 \Delta x \sqrt{\det \overline{h}})^{(m)} - (||\dot{S}||^2 ||\nabla x||^2 \sqrt{\det \bh})^{(m+1)} - (||\dot{S}||^2 ||\nabla x||^2 \sqrt{\det \bh})^{(m + 1, \log)} \right] 	
\end{empheq}
This shows that the second variation is computable in terms of the asymptotics of the metric and variational vector fields along $\gamma$. This proves theorem \ref{MainThm} in the $m$ is odd case. If $m < n$ even or $K = 0$ in the expansion \eqref{kEquation} then equation \eqref{SecondVarEq} becomes
\begin{empheq}[box=\fbox]{align*}
D^2 \mathcal{V} \Big|_{Y} (\dot{S}, \dot{S}) &= \int_{\gamma}[dx(\ddot{S}) \sqrt{\det \overline{h}}]^{m} - [dx(\dot{S})^2 \sqrt{\det \overline{h}}]^{(m + 1)}\\
& + \frac{1}{2} \left[(||\dot{S}||^2 \Delta x \sqrt{\det \overline{h}})^{(m)} - (||\dot{S}||^2 ||\nabla x||^2 \sqrt{\det \bh})^{(m+1)} \right] 	
\end{empheq}
This follows by the remarks on the $\F$ functional for $m < n$ or when $K = 0$. In section \ref{codim1}, we show that even in the case of $m = n$ even, the above holds. This will conclude the full statement of theorem \ref{MainThm}. \qed \nl \nl
Having given formulas for first and second variation, we show that these are independent of the choice of special bdf, and hence independent of the choice of representative of the conformal infinity, $k_0 = \overline{g} \Big|_{\gamma}$.
\subsection{Conformal Invariance of Variational Formulae for even dimensional submanifolds}
\label{ConformalInvariance}
In theorem \ref{EquivofRVonY}, we showed that the renormalized volume for an even dimensional submanifold can be computed with either $x$, the ambient special bdf, or $x_Y$, a special bdf for $Y$ itself. For completeness, we show that the first and second variations of the renormalized volume of $Y$ can be computed with $x_Y$, a special bdf on $Y$, or $x$, the ambient special bdf
\begin{proposition}
Let $\{Y_t^m\} \subseteq M^{n+1}$ a family of minimal submanifolds with $m$ even. Let $x: \overline{M} \to \R^{\geq 0}$ a special bdf for $M$ and  $x_t: Y_t \to \R^{\geq 0}$ a family of special bdfs for each $Y_t$, then 
\begin{align*}
\frac{d}{dt} \int_{Y_t} x_t^z dA_{Y_t} \Big|_{t = 0} &= \frac{d}{dt} \int_{Y_t} x^z dA_{Y_t} \Big|_{t = 0} \\
\frac{d^2}{dt^2} \int_{Y_t} x_t^z dA_{Y_t} \Big|_{t = 0} &= \frac{d^2}{dt^2} \int_{Y_t} x^z dA_{Y_t}  \Big|_{t = 0}
\end{align*}
\end{proposition}
%
%
\noindent \Pf \; For the first variation, recall $x_{Y_t} = x e^{\omega_t}$. Here, $\omega_t$ is initially defined on $Y_t$, but we can extend it to a uniform neighborhood of $Y_t$ near the boundary, called $U(Y_t) \cap \{x \leq \delta\}$ for $\delta$ sufficiently small (here $x$ is the ambient special bdf). For all $t < t_0$ sufficiently small (and perhaps $\delta$ smaller), we can extend each $\omega_t$ to a neighborhood of $Y = Y_0$ near the boundary, notated as $U(Y) \cap \{x \leq \delta\}$. Thus, it makes sense to define
\[
\omega: [0, t_0) \times U(Y) \cap \{x \leq \delta\} \to \R
\]
and we further notate $\omega_t = \omega(t, \cdot)$. With this, we compute
\[
\FPz \frac{d}{dt} \int_{Y_t} (x^z - x_{\omega_t}^z) dA_{Y_t} = \FPz \frac{d}{dt} \int_{Y_t} x^z(1 - e^{z \omega_t}) dA_{Y_t}
\]
Note that $dA_{Y_t}$ is the area form with respect to the complete metric, $g$ restricted to $Y_t$, and hence in this form is invariant, i.e. doesn't depend on the choice of boundary representative. Pulling back, we get
\begin{align*}
\frac{d}{dt} \int_{Y_t} x^z(1 - e^{z \omega_t}) dA_{Y_t} &= \frac{d}{dt} \int_{Y_0} F_t^*(x)^z(1 - e^{z F_t^*(\omega_t)}) F_t^*(dA_{Y_t}) \Big|_{t = 0} \\
& = \int_Y [z x^{z-1} \dot{x}(1 - e^{z \omega_0}) + x^z (-z(\dot{\omega} + \dot{\phi} \nu (\omega_0)) e^{z \omega_0}) ]dA_Y
\end{align*}
again, $\frac{d}{dt} F_t^*(dA_{Y_t})\Big|_{t = 0} = 0$. For the first term
\begin{align*}
z x^{z-1} \dot{x} (1 - e^{z \omega_0}) &= z^2 x^{z - 1} \dot{x} (z \omega_0 + O(z^2)) = z^2 x^{z-1} \dot{x} \omega_0 \\
\implies \FPz \int_Y z x^{z-1} \dot{x} (1 - e^{z \omega_0}) &= -\int_{\gamma} [ \dot{x} \omega_0 \sqrt{\det \bh}]^{\log, m} \\
&= 0
\end{align*}
The vanishing of this integral is seen in that $\dot{x} = O(x)$ and $\F(\dot{x}) = -1$ so that 
\[
[ \dot{x} \omega_0 \sqrt{\det \bh}]^{\log, m} = [ \dot{x}]^{\log, m} \omega_0
\]
but then for $m < n$, $[ \dot{x}]^{\log, m} = 0$ by the regularity and expansion of $\dot{\phi}^i$. When $m = n$, one makes explcit use of the expansion of $\bnu$ in \S \ref{Normal} to show that this is $0$. For the second term, we get 
\begin{align} \nonumber
	\FPz z \int_Y x^z (\dot{\omega} + \dot{\phi} \nu (\omega_0)) e^{z \omega_0} dA_Y &= \FPz z \int_Y x^z (\dot{\omega} + \dot{\phi} \nu (\omega_0)) (1 + z \omega_0 + O(z^2))dA_Y \\ \nonumber
	& =  \FPz \int_Y[z x^z (\dot{\omega} + \dot{\phi} \nu (\omega_0)) + z^2 \omega_0 (\dot{\omega} + \dot{\phi} \nu (\omega_0)) x^z + O(z^3)]dA_Y\\ \label{VarEquivFinal}
	& = \int_{\gamma} [(\dot{\omega} + \dot{\phi} \nu (\omega_0)) \sqrt{\det \overline{h}}]^{(m-1)} + [\omega_0 (\dot{\omega} + \dot{\phi} \nu (\omega_0)) \sqrt{\det \bh}]^{\log, m-1}
\end{align}
where we've again noted that all quadratic terms in $z$ vanish under finite part evaluation at $z = 0$ from lemma \ref{FiniteEvaluation}. For the last two integrals, note that $\F\left(\dot{\omega}|_Y\right) = 1$ since $\omega_t$ will be even for each $t$ by the parity of $\omega$ as in theorem \ref{EquivofRVonY} (applied to each $\omega_t$ individually). The fact that 
\[
\F(\omega_0) = \F(\dot{\phi} \nu(\omega_0)) = 1
\]
on $Y$ again follows from the parity of $\omega_0$ as well as the parity of $\dot{\phi}$ and $\nu$. This means that the $m-1$ and $[\log, m-1]$ terms in equation \eqref{VarEquivFinal} vanish. \nl \nl
The second variation proceeds similarly, making use of the $\F$ functional for $m < n$ and the expansion of $\bnu$ when $m = n$ even. For $m$ odd, renormalized volume is not conformally invariant. Using the above process, one could compute how the first and second variations depend on the initial choice of bdf. \qed \nl \nl
Having derived formulae for first and second variations, we simplify it for codimension $1$ submanifolds. In this case, $\sqrt{\det \bh}$ and $\bh^{\alpha \beta}$ are tractable in terms of the coefficients $\{u^i_k\}$. Similarly, terms involving $\dot{S}$ and $\ddot{S}$ simplify with $\dot{S} = \dot{\phi} \bnu$ and $\ddot{S} = \ddot{\phi} \bnu$. In particular, $\bnu(x,s)$, the unit normal to $Y^n \subseteq M^{n+1}$ (with respect to the complete metric), is computable in terms of $\{u^i_k\}$.
\section{Codimension $1$ case}
\label{codim1}
For $Y^n \subseteq M^{n+1}$, we write 
\[
\dot{S} = \dot{\phi} \bnu
\]
where $\dot{\phi} = O(1)$ and $\bg(\bnu, \bnu) = 1$. In this section, we compute $\bnu$ explicitly and show
\begin{theorem} \label{codim1FirstVar}
For $\{Y_t^n\} \subseteq M^{n+1}$ a family of minimal hypersurfaces, we have
\begin{align*}
D \mathcal{V}(Y) &= - (n+1)\int_{\gamma} \dot{\phi}_0(s)  u_{n+1}(s) ds	
\end{align*}
for $n$ even and 
\begin{align*}
D\mathcal{V}(Y) = - \int_{\gamma} \left( [(n+1) u_{n+1}(s) + U(s)] \dot{\phi}_0(s) + F(\dot{\phi}_0(s), u_2(s)) \right) ds
\end{align*}
for $n$ odd.
\end{theorem}
\noindent Now recall $k_{n+1} = [k(s,z=0,x)]^{n+1}$ from \eqref{kEquation}. We have
\begin{theorem} \label{codim1SecondVar}
For $\{Y_t^n\} \subseteq M^{n+1}$ a family of minimal hypersurfaces, we have
\begin{align*}
D^2 \mathcal{V}(Y) &= \int_{\gamma} -(n+1) \ddot{\phi}_0 u_{n+1} + (1-n) \dot{\phi}_0 \dot{\phi}_{n + 1} - 8n\dot{\phi}_0^2 u_2 u_{n+1}
\end{align*}
for $n$ even and
\begin{align*}
D^2 \mathcal{V}(Y) &= \int_{\gamma} -(n+1) \ddot{\phi}_0 u_{n+1} + (1-n) \dot{\phi}_0 \dot{\phi}_{n+1} + \dot{\phi}_0^2 \left[-8n u_2 u_{n+1} + \text{Tr}_{T\gamma}(k_{n+1,0}) \right] \\
& - \dot{\phi}_0\left[ 4(n+2) \dot{\phi}_0 u_2 U + \dot{\Phi}\right] +  F(\ddot{\phi}_0, \dot{\phi}_0, u_2)
\end{align*}
for $n$ odd.
\end{theorem}
\subsection{Computing the normal} \label{Normal}
As in \S \ref{FermiCoords}, let $f: U \subseteq \partial M \to \gamma$ be a map in fermi coordinates. Let $p = f(s)$ and $\{e_a(s) = \gamma_a(s) = \partial_{s_a}\}$ be a frame for $T_p \gamma$ in some neighborhood of $p \times \{0\}$. Translate this to all of $\gamma \times [0, \eps) = \Gamma$. Complete the basis with $\partial_x$ and $N(s,x) = \partial_z$ a normal coordinate for $\Gamma$ such that $N(s,0) = N(s)$ a unit normal for $\gamma \subseteq \partial M$. Then let 
\[
G(s,x) = \overline{\exp}_{\Gamma}(u(s,x) N(s,x))
\]
which in Fermi coordinates of $(s,x,z)$ can be written as 
\[
G(s,x) = (s,x, z = u(s,x))
\]
The tangent space is spanned by
\begin{align*}
v_a &= G_*(\partial_{s_a}) = \partial_{s_a} + u_a(s,x) \partial_z \\
\vx &= G_*(\partial_x) = \partial_x + u_x(s,x) \partial_z 
\end{align*}
having noted that $\bGamma_{xz}^x = 0$ in the second line. Now we can compute the normal to $Y$, $\bnu$, by projecting onto the tangent basis. We have
\begin{align*}
\tilde{\nu} &= \partial_z - \overline{h}^{ab} \overline{g}(\partial_z, v_a) v_b - \overline{h}^{ax} \overline{g}(\partial_z, v_a) v_x \\
& - \overline{h}^{xa} \overline{g}(\partial_z, v_x) v_a - \overline{h}^{xx} \overline{g}(\partial_z, v_x) v_x \\
&= \partial_z - \sum_a (u_a + R_a) v_a - (u_x + R_x) v_x \\
&= (1 - u_x^2 + \tilde{R}_z) \partial_z - \sum_a (u_a + \tilde{R}_a) \partial_a - (u_x + \tilde{R}_x) \partial_x
\end{align*}
where 
\begin{align*}
\tilde{R}_z \Big|_{z = u(s,x)} &= O(x^4), &\F\left( x^{-2}\tilde{R}_z \Big|_{z = u(s,x)} \right) &= 1 \\
\tilde{R}_a \Big|_{z = u(s,x)} &= O(x^4), &\F\left( x^{-2}\tilde{R}_a \Big|_{z = u(s,x)} \right) &= 1 \\
\tilde{R}_x \Big|_{z = u(s,x)} &= O(x^3), &\F\left( x^{-2}\tilde{R}_x \Big|_{z = u(s,x)} \right) &= -1
\end{align*}
The $\F(\cdot)$ statements actually say that $\tilde{R}_z$, $\tilde{R}_a$ are even to order $n+2$, and $\tilde{R}_x$ odd to order $n + 1$. Moreover, note that in this decomposition
\[
\bg(\tilde{\nu}, \tilde{\nu}) = 1 + O(x^2), \qquad \F(\bg(\tilde{\nu}, \tilde{\nu})) = 1
\]
so that if we normalize, we get 
\begin{align*}
\bnu &= \frac{\tilde{\nu}}{\sqrt{\overline{g}(\tilde{\nu}, \tilde{\nu})}} \\
&= c^z \partial_z + c^x \partial_x + c^a \partial_{s_a} 
\end{align*}
such that 
\begin{align*}
c^z &= 1 + \tilde{R}_z,  & \F(x^{-2} \tilde{R}_z) &= 1 \\
c^x &= -u_x + \tilde{R}_x = O(x),  & \F(x^{-2} \tilde{R}_x) &= -1 \\
c^a &= -u_a + \tilde{R}_z = O(x^2),  & \F(x^{-2} \tilde{R}_a) &= 1 \\
\end{align*}
In particular, because $u$ lacks $x^n \log(x)$, $x^{n+1} \log(x)$ terms, we see that $c^z$, $c^x$, and $c^a$ lack $x^n \log(x)$, $x^{n+1} \log(x)$ terms.
\subsection{First variation, codimension $1$, even}
We now apply the first variation formula to prove theorem \ref{codim1FirstVar} in the even case. Consider the expansion for $\dot{S} = \dot{\phi} \bnu$
\[
D\mathcal{V}(Y) = \int_{\gamma} [dx(\dot{S}) \sqrt{\det \overline{h}}]^{(n)} = -\int_{\gamma} [\dot{\phi} (u_x + \tilde{R}_x)\sqrt{\det \overline{h}}]^{(n)}
\]
since $n$ is even and $\F(\dot{\phi}) = \F(\sqrt{\det \overline{h}}) = 1$, while $\F(u_x) = \F(x^{-2} \tilde{R}_x) = -1$, we have that 
\begin{align*}
[\dot{\phi} (u_x + \tilde{R}_x)\sqrt{\det \overline{h}}]^{(n)} &= [\dot{\phi} u_x \sqrt{\det \overline{h}}]^{(n)} + [\dot{\phi} \tilde{R}_x\sqrt{\det \overline{h}}]^{(n)} \\
&= [\dot{\phi}]_0 [u_x]_{n} [\sqrt{\det \overline{h}}]_0 + [\dot{\phi}]_0 [\tilde{R}_x]_{n} [\sqrt{\det \overline{h}}]_0 
\end{align*}
Now we compute
\[
[\dot{\phi}]_0 [u_x]_{n} [\sqrt{\det \overline{h}}]_0 = \dot{\phi}_0 (n+1) u_{n+1} 
\]
Here, we've noted that $[\sqrt{\det \overline{h}}]_0 = 1$. Thus
\[
\boxed{D \mathcal{V}(Y) = - (n+1)\int_{\gamma} \dot{\phi}_0(s)  u_{n+1}(s) ds} 
\]
This finishes the proof of theorem \ref{codim1FirstVar} in the even case. \qed  \nl \nl
Recall that for $Y^n \subseteq M^{n+1}$ \textit{non-degenerate}, any $\dot{\phi}_0(s) \in C_c^{\infty}(\gamma)$ can be extended to a Jacobi field, $\dot{\phi}$, on all of $Y$ (see \S \ref{Degeneracy}) for which we can define $\phi_t = t \dot{\phi}$ and hence $Y_t$. In this case, we have: 
\begin{corollary}
If $Y^n \subseteq M^{n+1}$ is a nondegenerate minimal submanifold and a critical point for renormalized volume with $n$ even, then $u_{n+1}(s) \equiv 0$.
\end{corollary}
\noindent When $Y$ is degenerate, the set of $\{\dot{\phi}_0(s)\}$ which can be extended to an $L^2$ Jacobi field on all of $Y$ are orthogonal to a finite dimensional kernel (cf \S \ref{Degeneracy} and \cite{alexakis2010renormalized} for details). \nl \nl
\noindent \rmk \; As seen in \S \ref{IterationArgument}, $u_2(s)$ determines the coefficients $\{u_{2k}(s)\}_{k = 2}^{n/2}$ via the minimal surface system. We think of these terms as ``local" in the sense that they are determined by the boundary $\gamma$, which determines $u_2(s) = \frac{1}{2(n-1)} H_{\gamma}(s)$. By contrast, $u_{n+1}(s)$ is \textit{not} determined by $u_2(s)$, and hence is ``global". The rest of the expansion of $u(x,s)$ is determined by $\gamma$ and $u_{n+1}(s)$ and continuing the iteration. In a loose sense, $u_{n+1}(s)$ represents the output of a Dirichlet-to-Neumann map that is posed: given $\gamma^{m-1} \subseteq \partial M$, find $Y^m \subseteq M^{n+1}$ minimal with $\partial Y = Y \cap \partial M^{n+1} = \gamma$. The Dirichlet data is $\gamma$, and the image of the Dirichlet-to-Neumann map is the first undetermined term in the series expansion, $u_{n+1}(s)$. Thus, \textit{for nondegenerate critical points of renormalized volume, the Neumann output is exactly $0$}. 
%
%
\subsection{Codimension $1$, First variation, odd}
\label{Codim1FirstOdd}
We demonstrate that the first variation is not as transparent when $n$ is odd. We compute
\[
[\dot{\phi} (u_x + \tilde{R}_x) \sqrt{\det \bh} ]^{(n)}
\]
where 
\begin{align*}
u_x &= 2 u_2(s) x + \dots + (n-1) u_{n-1}(s) x^{n-2} + [(n+1) u_{n+1}(s) + U(s)] x^n + (n+1) U(s) x^n \log(x) + O(x^{n+1}) \\
\dot{\phi} & = \dot{\phi}_0(s) + \dots + \dot{\phi}_{n-1}(s) x^{n-1} + \dot{\phi}_{n+1}(s) x^{n+1} + \dot{\Phi}(s) x^{n+1} \log(x) + O(x^{n+2}) \\
\tilde{R}_x(x,s) &= \tilde{R}_{x,3}(s) x^3 + \dots + \tilde{R}_{x, n+3}(s) x^{n+3} + \mathscr{R}_{x}(s) x^{n+3} \log(x) + O(x^{n+4}) \\
\sqrt{\det \bh}(s,x) & = \overline{q}_0(s) + \dots + \overline{q}_{n-1}(s) x^{n-1} + \overline{q}_{n+1}(s) x^{n+1} + \overline{Q}(s) x^{n+1} \log(x) + O(x^{n+2}) 
\end{align*}
We can already see that there are many combinations that multiply to form an $x^n$ term, e.g. 
\[
2u_2 \dot{\phi}_0 \overline{q}_{n-1}, \quad \tilde{R}_{x,3} \dot{\phi}_0 \overline{q}_{n-3}, \quad\dots
\]
however, we can write the $n$th term as 
\begin{align*}
[\dot{\phi} (u_x + \tilde{R}_x) \sqrt{\det \bh} ]^{(n)} &= [(n+1) u_{n+1}(s) + U(s)]\dot{\phi}_0 + P(\{\dot{\phi}_{2k}\}_{k = 0}^{(n-1)/2}, \{u_{2k}\}_{k = 1}^{(n-1)/2}, \{\tilde{R}_{x, 2k+1}\}_{k = 1}^{(n-1)/2}, \{\overline{q}_{2k}\}_{k = 0}^{(n-1)/2}\})
\end{align*}
Clearly, $P$ is determined by terms of order $n-1$ or lower, so we write
\[
F(\dot{\phi}_0(s), u_2(s)):= P(\{\dot{\phi}_{2k}\}_{k = 0}^{(n-1)/2}, \{u_{2k}\}_{k = 1}^{(n-1)/2}, \{\tilde{R}_{x,2k}\}_{k = 0}^{(n-1)/2}, \{\overline{q}_{2k}\}_{k = 0}^{(n-1)/2}\})
\]
noting implicitly that $\{\tilde{R}_{x, 2k+1} \}_{k=1}^{(n-1)/2}$ is determined by $\{u_{2k}\}_{k = 1}^{(n-1)/2}$, which follows from the construction in \S \ref{Normal}. \nl \nl
\rmk \; From hereon, we will use $F(\dot{\phi}_0, u_2)$ to denote a polynomial function of $\{\dot{\phi}_0, \dot{\phi}_2, \dots, \dot{\phi}_{n-1}\}$ and $\{u_2, u_4, \dots, u_{n-1}\}$. In \S \ref{IterationArgument}, we showed that $\dot{\phi}_{2k}(s)$ and $u_{2k}(s)$ are determined by $\dot{\phi}_0(s)$ and $u_2(s)$, respectively, for $2k\leq n-1$. Because of this, we can think of $F$ as a non-linear differential operator acting on $\dot{\phi}_0$ and $u_2$. We will make a slight abuse of notation and write ``$F(\dot{\phi}_0, u_2)$" wherever such a function appears, as opposed to having a distinct labeling for each such function of $(\dot{\phi}_0, u_2)$. We will make the same convention for functions $R(u_2)$, which are the same as $F(\dot{\phi}_0, u_2)$ when there is no $\dot{\phi}_0$ dependence. We will also use such convention for functions $F(\ddot{\phi}_0, \dot{\phi}_0, u_2)$ when there is dependence on $\{\ddot{\phi}_0, \dots, \ddot{\phi}_n\}$. \nl \nl
%
We conclude 
\[
\boxed{D\mathcal{V}(Y) = - \int_{\gamma} \left( [(n+1) u_{n+1}(s) + U(s)] \dot{\phi}_0(s) + F(\dot{\phi}_0(s), u_2(s)) \right) dA_{\gamma}}
\]
proving theorem \ref{codim1FirstVar} in the odd case. \qed
\subsection{Codimension $1$, Second variation, Even}
The formula of interest is 
%
%
\begin{align*}
D^2 \mathcal{V} \Big|_{Y} (\dot{S}, \dot{S}) &= \int_{\gamma}- [dx(\dot{S})^2 \sqrt{\det \overline{h}}]^{(n + 1, \log)} - [dx(\dot{S})^2 \sqrt{\det \overline{h}}]^{(n + 1)}\\
& + \frac{1}{2} \left[(||\dot{S}||^2 \Delta x \sqrt{\det \overline{h}})^{(n)} - (||\dot{S}||^2 ||\nabla x||^2 \sqrt{\det \bh})^{(n+1)}\right] \\
& - \frac{1}{2} (||\dot{S}||^2 ||\nabla x||^2 \sqrt{\det \bh})^{(n + 1, \log)}  + [dx(\ddot{S}) \sqrt{\det \bh}]^n	\\
&= I_1 + I_2 + I_3 + I_4 + I_5 + I_6
\end{align*}
We'll look at each of the summands individually.
\subsubsection{$I_1$}
As in \S \ref{ConformalInvariance}, we have 
\begin{align*}
dx(\dot{S})^2 &= [-u_x + \tilde{R}_x]^2 \dot{\phi}^2 = O(x^2) \\
\implies [dx(\dot{S})^2 \sqrt{\det \overline{h}}]^{(n + 1, \log)} &= [dx(\dot{S})^2]^{(n + 1, \log)} \cdot [\sqrt{\det \overline{h}}]^0 \\
&= [dx(\dot{S})^2]^{(n + 1, \log)} \\
&= [(-u_x + \tilde{R}_x)^2]^{n+1, \log} \cdot \dot{\phi}_0^2
\end{align*}
having used that $\F(\dot{\phi}) = 1$ and has no $x^n\log$ or $x^{n+1} \log(x)$ terms. The same holds for $u$ and $\F(x^{-2} \tilde{R}_x) = -1$, so we see that 
\[
[(-u_x + \tilde{R}_x)^2]^{n+1, \log} = [u_x^2]^{n+1, \log} = 0
\]
hence
\[
I_1 = 0
\]
\subsubsection{$I_2$}
In this case, similar reasoning holds and we see that 
\begin{align*}
[dx(\dot{S})^2 \sqrt{\det \overline{h}}]^{n + 1} &= [(-u_x + \tilde{R}_x)^2]^{n+1} \cdot \dot{\phi}_0^2 \\
&= [u_x^2]^{n+1} \cdot \dot{\phi}_0^2 \\
&= 4 u_2 u_{n+1} \dot{\phi}_0^2
\end{align*}
so 
\[
I_2 = -\int_{\gamma} 4 (n+1) u_2 u_{n+1} \dot{\phi}_0^2
\]
\subsubsection{$I_3$}
For the second term, we have 
\[
I_2 = \frac{1}{2} \int_{\gamma} [||\dot{S}||^2 \Delta x \sqrt{\det\overline{h}}]^{(n)}
\]
recall the notation of $\{y_{\alpha}\}$ denoting any coordinate of $\{s_a, x\}$, and that $\Delta = \Delta_Y$ represents the laplacian on $Y$ with respect to the complete (induced) metric, $h$.
We decompose
\begin{align*}
\sqrt{\det \overline{h}}\; \Delta x &= x^n \partial_{y_\alpha}(\sqrt{\det h} h^{\alpha \beta} \partial_{y_\beta} x) = x^n \partial_{y_\alpha} (x^{2-n} \sqrt{\det \overline{h}} \; \overline{h}^{\alpha \beta} \partial_{y_\beta} x) \\
&= x^2 \partial_{s_a}(\sqrt{\det \overline{h}} \; \overline{h}^{a x}) + (2-n) x \sqrt{\det \overline{h}} \; \overline{h}^{xx} + x^2 \partial_x(\sqrt{\det \overline{h}} \; \overline{h}^{xx})
\end{align*}
hence
\[
[||\dot{S}||^2 \Delta x \sqrt{\det \overline{h}}]^{(n)} = [\dot{\phi}^2 \partial_{s_a}(\sqrt{\det \overline{h}} \; \overline{h}^{ax})]^{(n)} + (2-n) [\dot{\phi}^2 \sqrt{ \det \overline{h}} \; \overline{h}^{xx}]^{(n+1)} + [\dot{\phi}^2 \partial_x(\sqrt{\det \overline{h}} \; \overline{h}^{xx})]^{(n)}
\]
From our previous work, $\overline{h}^{ax}$ is $O(x^3)$ and is odd up to order $(n+1)$. This tells us that $\partial_{s_a}(\sqrt{\det \overline{h}} \; \overline{h}^{ax})$ is odd up to order $(n + 1)$ and $O(x^3)$ and so because $\F(\dot{\phi}^2) = 1$, we have that 
\[
[\dot{\phi}^2 \partial_{s_a}(\sqrt{\det \overline{h}} \; \overline{h}^{ax})]^{(n)} = 0
\]
For the second term, we do the same analysis as before: $\dot{\phi}^2$, $\overline{h}^{xx}$, $\sqrt{\det \overline{h}}$ all satisfy $\F(\cdot) = 1$, so any $(n+1)$st term must come from the $(n+1)$st term in one of the factors multiplied by the $0$th order term in the remaining factors. We get
\[
(2-n) [\dot{\phi}^2 \sqrt{ \det \overline{h}} \; \overline{h}^{xx}]^{(n+1)} = (2 - n)[ 2 \dot{\phi}_0 \dot{\phi}_{n+1}  + \dot{\phi}_0^2 (\overline{h}_{n+1}^{xx} + \overline{q}_{n+1})]
\]
For the last term of $[\dot{\phi}^2 \partial_x(\sqrt{\det \bh} \bh^{xx})]^{(n)}$, we know that $\F(\partial_x (\sqrt{\det \overline{h}} \; \overline{h}^{xx})) = -1$ because both $\F(\sqrt{\det \overline{h}}) = \F(\bh^{xx}) = 1$ and so the derivative of their product is $O(x)$ and satisfies $\F(\cdot) = -1$. On the other hand $\F(\dot{\phi}^2) = 1$. Thus the $n$th term of the product can only come from the $n$th term of $\partial_x(\sqrt{\det \overline{h}} \; \overline{h}^{xx})$ paired with the $0$th order term of $\dot{\phi}^2$. Recall from corollary \ref{NoLogTerms} that $\sqrt{\det \overline{h}} \; \overline{h}^{xx}$ has no $x^{n+1} \log(x)$ term. Thus
\[
[\dot{\phi}^2 \partial_x(\sqrt{\det \overline{h}} \; \bh^{xx})]^{(n)} = (n+1) \dot{\phi}_0^2 (\overline{q}_{n+1} + \overline{h}_{n+1}^{xx})
\]
since $\overline{q}_0 = \overline{h}_0^{xx} = 1$. Thus
\[
I_3 = \int_{\gamma} (2-n) \dot{\phi}_0 \dot{\phi}_{n+1} + \frac{3}{2} \dot{\phi}_0^2 (\overline{q}_{n+1} + \overline{h}_{n+1}^{xx})
\]

\subsubsection{$I_4$, even}
For the third term 
\[
I_3 = -\frac{1}{2} \int_{\gamma} [||\dot{S}||^2 \cdot ||\nabla x||^2 \sqrt{\det \overline{h}} ]^{(n+1)}
\]
We expand using
\begin{align*}
||\nabla x||^2 &= h^{\alpha \beta} v_{\alpha}(x) v_{\beta}(x) \\
&= h^{xx} v_x(x) v_x(x) \\
&= h^{xx}\\
||\dot{S}||^2 &= x^{-2} \dot{\phi}^2
\end{align*}
From \eqref{ChristoffelParityCylinder} and corollaries \ref{MetricExpansion}, \ref{NoLogTerms}, we see that $\F(x^{-2} ||\nabla x||^2) = 1$. Similarly, $\F(x^2 ||\dot{S}||^2) = 1$. Thus
\begin{align*}
[||\dot{S}||^2 ||\nabla x||^2 \sqrt{\det \overline{h}} ]^{(n+1)} &= [ \dot{\phi}^2 (x^{-2} ||\nabla x||^2) \sqrt{\det \bh}]^{n+1} \\
&= [\dot{\phi}^2]^{n+1} [x^{-2} ||\nabla x||^2]^0 [\sqrt{ \det \bh}]^0 + [\dot{\phi}^2]^{0} [x^{-2} ||\nabla x||^2]^{n+1} [\sqrt{ \det \bh}]^0 \\
&+ [\dot{\phi}^2]^{0} [x^{-2} ||\nabla x||^2]^0 [\sqrt{ \det \bh}]^{n+1}
\end{align*} 
Moreover, by the vanishing orders of $u_i$ and $\Gamma_{ai}^x$ and $\Gamma_{xi}^x$, we have that 
\begin{align*}
[x^{-2}||\nabla x||^2]^{0} &= 1 \\
[x^{-2}||\nabla x||^2]^{n+1} &= [\bh^{xx}]^{n+1} = \bh_{n+1}^{xx} 
\end{align*}
Similarly, 
\begin{align*}
[\dot{\phi}^2]^{0} &= \dot{\phi}_0^2 \\
[\dot{\phi}^2]^{n+1} &= 2 \dot{\phi}_0 \dot{\phi}_{n+1} \\
[\sqrt{ \det \bh}]^{0} &:= 1 \\
[\sqrt{ \det \bh}]^{n+1} &:= \overline{q}_{n+1} 	
\end{align*}
From this, we get that
\[
[||\dot{S}||^2 ||\nabla x||^2 \sqrt{\det \overline{h}} ]^{(n+1)} = 2 \dot{\phi}_0 \dot{\phi}_{n+1} + \dot{\phi}_0^2( \overline{h}_{n+1}^{xx} + \overline{q}_{n+1})
\]
This gives
\[
I_4 = \int_{\gamma} - \dot{\phi}_0 \dot{\phi}_{n+1} - \frac{1}{2}(\overline{h}_{n+1}^{xx} + \overline{q}_{n+1}) \dot{\phi}_0^2
\]

\subsubsection{$I_5$}
As in $I_4$, we have that 
\begin{align*}
||\dot{S}||^2 ||\nabla x||^2 \sqrt{\det \bh} &= \dot{\phi}^2 (x^{-2} ||\nabla x||^2) \sqrt{\det \bh} \\
\implies [||\dot{S}||^2 ||\nabla x||^2 \sqrt{\det \overline{h}} ]^{\log, n+1} &= [\dot{\phi}^2]^{\log, n+1} [x^{-2} ||\nabla x||^2]^0 [\sqrt{\det \bh}]^0  \\
&+ [\dot{\phi}^2]^{0} [x^{-2} ||\nabla x||^2]^{\log, n+1} [\sqrt{\det \bh}]^0 + [\dot{\phi}^2]^{0} [x^{-2} ||\nabla x||^2]^{0} [\sqrt{\det \bh}]^{\log, n+1} 
\end{align*}
We first recall that $[\dot{\phi}]^{\log, n+1} = 0$ which gives 
\[
[\dot{\phi}^2]^{\log, n+1} = 0
\]
we also note that 
\begin{align*}
x^{-2} ||\nabla x||^2 &= \bh^{xx} (1 + u^i \Gamma_{xi}^x)^2 + 2 \bh^{xa} u^i \Gamma_{ai}^x (1 + u^i \Gamma_{xi}^x) \\
& + \bh^{ab} u^i \Gamma_{ai}^x u^j \Gamma_{bj}^x \\
\implies [x^{-2} ||\nabla x||^2]^{\log, n+1} &= [\bh^{xx}]^{\log, n+1} = 0 \\
\implies [x^{-2} ||\nabla x||^2]^{0} &= [\bh^{xx}]^{0} = 1
\end{align*}
by the vanishing order of the other expression and again from corollary \ref{NoLogTerms} that $\bh_{xx}$ has no $x^{n+1} \log(x)$ term. Finally, from the same remark, we see that $[\sqrt{\det \bh}]^{\log, n+1} = 0$. Thus
\[
I_5 = 0
\]
\subsubsection{$I_6$}
We compute
\begin{align*}
[dx(\ddot{S}) \sqrt{\det \bh}]^n &= [\ddot{\phi} dx(\nu) \sqrt{ \det \bh}]^n\\
&= \ddot{\phi}_0 [-u_x + \tilde{R}_x]^n [ \sqrt{ \det \bh}]_0 \\
&= - (n+1) u_{n+1}\ddot{\phi}_0 
\end{align*}
This comes from section \S \ref{Normal} where 
\begin{align*}
\nu &= c^z \partial_z + c^x \partial_x + c^a \partial_{s_a} \\
c^x &= -u_x + R_x \\
\F(x^{-2} R_x) &= 1 
\end{align*}
Since $u_x = O(x)$ and is odd, the $n$th order term of $(\ddot{\phi} dx(\nu) \sqrt{ \det \bh})$ can only come from the $0$th order terms of $\ddot{\phi}$ and $\sqrt{\det \bh}$ and the $n$th order term of $c^x$. This is precisely because $\F(\ddot{\phi}) = \F(\sqrt{\det \bh}) = 1$. Thus
\[
I_6 = \int_{\gamma} -(n+1) \ddot{\phi}_0 u_{n+1} dA_{\gamma}
\]
\subsubsection{The Full Expression, Even Case}
Together, 
\[
\boxed{D^2 \mathcal{V}(Y) = \int_{\gamma} -(n+1) \ddot{\phi}_0 u_{n+1} + (1 - n)\dot{\phi}_0 \dot{\phi}_{n+1} + \dot{\phi}_0^2 \left[ (\overline{h}_{n+1}^{xx} + \overline{q}_{n+1}) - 4(n+1) u_2 u_{n+1}\right]}
\]
in the codimension $1$ case. In \S \ref{SecondVarComp}, we show by a more careful analysis that
\[
\overline{h}^{xx}_{n + 1} + \overline{q}_{n + 1} = -4 (n-1) u_2 u_{n+1}
\]
and so the final formula is
\[
\boxed{D^2 \mathcal{V}(Y) = \int_{\gamma} -(n+1) \ddot{\phi}_0 u_{n+1} + (1-n) \dot{\phi}_0 \dot{\phi}_{n + 1} - 8n \dot{\phi}_0^2  u_2 u_{n+1}} 
\]
%
%
In particular, noting that $I_1 = I_5 = 0$, this finishes the proof of theorem \ref{MainThm} and theorem \ref{codim1SecondVar} in the even case.
\subsection{Codimension $1$, Second Variation, Odd} \label{Codim1SecondOdd}
Note that when $n$ is odd, there is no $x^n \log(x)$ or $x^{n+1} \log(x)$ terms in equation \eqref{kEquation}. So terms of the form $[ \cdot ]^{\log, n}, [\cdot]^{\log, n+1}$ will only come from $[u]^{\log, n}, [u]^{\log, n+1}$ or $[\dot{\phi}]^{\log, n}, [\dot{\phi}]^{\log, n+1}$. Recall that the formula is given by
\begin{align*}
D^2 \mathcal{V} \Big|_{Y} (\dot{S}, \dot{S}) &= \int_{\gamma}- [dx(\dot{S})^2 \sqrt{\det \overline{h}}]^{(n + 1, \log)} - [dx(\dot{S})^2 \sqrt{\det \overline{h}}]^{(n + 1)}\\
& + \frac{1}{2} \left[(||\dot{S}||^2 \Delta x \sqrt{\det \overline{h}})^{(n)} - (||\dot{S}||^2 ||\nabla x||^2 \sqrt{\det \bh})^{(n+1)}\right] \\
& - \frac{1}{2} (||\dot{S}||^2 ||\nabla x||^2 \sqrt{\det \bh})^{(n + 1, \log)}  + [dx(\ddot{S}) \sqrt{\det \bh}]^n	\\
&= I_1 + I_2 + I_3 + I_4 + I_5 + I_6
\end{align*}
where $[f(x,s)]^{(n+1, \log)}$ denotes the coefficient of the $x^{n+1} \log(x)$ term for $f(x,s)$. 
\subsubsection{$I_1$}
Similar to the even case,
\[
I_1 = - \int_{\gamma} [dx(\dot{S})^2 \sqrt{\det \bh}]^{(n+1, \log)} dA_{\gamma} = - \int_{\gamma} [\dot{\phi}^2 (u_x + \tilde{R}_x)^2 \sqrt{\det \bh}]^{(n+1, \log)}
\]
From the expansions used in the first variation formula, adapted to the odd case, we have
\begin{align*}
\dot{\phi}^2 &= \dot{\phi}_0(s)^2 + \dots + 2 \dot{\phi}_0(s) \dot{\phi}_{n+1}(s) x^{n+1} + 2\dot{\Phi}(s) \dot{\phi}_0(s) x^{n+1} \log(x) + \dots \\
(u_x + \tilde{R}_x)^2 &= 4 u_2(s)^2 x^2 + \dots + 4[(n+1) u_{n+1}(s) + U(s) + \tilde{R}_{x,n-1}] u_2(s) x^{n+1} + 4 (n+1) u_2(s) U(s) x^{n+1} \log(x) + \dots \\
\sqrt{ \det \bh}(x,s) & = \overline{q}_0(s) + \dots + \overline{q}_{n+1}(s) x^{n+1} + \overline{Q}(s) x^{n+1} \log(x) + \dots 
\end{align*}
so the $x^{n+1} \log(x)$ coefficient of the product is
\[
[\dot{\phi}^2 (u_x + \tilde{R}_x)^2 \sqrt{\det \bh}]^{(n+1, \log)} = \dot{\phi}_0(s)^2 4(n+1) u_2(s) U(s) \overline{q}_0(s) = 4(n+1) u_2(s) U(s) \dot{\phi}_0(s)^2
\]
Together,
\[
\boxed{I_1 = - \int_{\gamma} 4(n+1) u_2(s) U(s) \dot{\phi}_0(s)^2 dA_{\gamma}(s) }
\]

\subsubsection{$I_2$}
Again, similar to the even case:
\[
I_2 = - \int_{\gamma} [dx(\dot{S})^2 \sqrt{\det \bh}]^{(n+1)} dA_{\gamma} = - \int_{\gamma} [\dot{\phi}^2 (u_x + \tilde{R}_x)^2 \sqrt{\det \bh}]^{(n+1)}
\]
As in the first variation for $n$ odd, there are many terms in this integrand which combine to give an $n+1$st term because $n+1$ is even. However, we isolate the terms which involve $u_{n+1}$, $U$, $\dot{\phi}_{n+1}$, and $\Phi$ and combine the lower order terms:
\[
\boxed{I_2 = - \int_{\gamma} \left( 4[(n+1) u_{n+1}(s) + U(s)] u_2(s) \dot{\phi}_0(s)^2 + F(\dot{\phi}_0, u_2) \right) dA_{\gamma}(s) }
\]
\subsubsection{$I_3$, odd}
Here,
\begin{align*}
I_3 &= \int_{\gamma} \frac{1}{2} (||\dot{S}||^2 (\Delta x) \sqrt{\det \bh} )^{(n)} dA_{\gamma}(s) \\
&= \int_{\gamma} [\dot{\phi}^2 \partial_{s_a} (\sqrt{ \det \bh} \; \bh^{ax})]^{(n)} + (2-n) [\dot{\phi}^2 \sqrt{ \det \bh} \; \bh^{xx} ]^{(n+1)} + [\dot{\phi}^2 \partial_x(\sqrt{\det \bh} \; \bh^{xx})]^{(n)}
\end{align*}
Write the first term as $F(\dot{\phi}_0, u_2)$ as no $(n+1)$st or $(n+1, \log)$ coefficients appear. We compute the middle and last terms as follows:
\begin{align*}
[\dot{\phi}^2 \sqrt{ \det \bh} \; \bh^{xx} ]^{(n+1)} &= 2 \dot{\phi}_0 \dot{\phi}_{n+1} \overline{q}_0 \bh^{xx}_0 + \dot{\phi}_0^2 \overline{q}_{n+1} \bh^{xx}_0 + \dot{\phi}_0^2 \overline{q}_0 \bh^{xx}_{n+1}  + F(\dot{\phi}_0, u_2)\\
[\dot{\phi}^2 \partial_x(\sqrt{\det \bh} \; \bh^{xx})]^{(n)} &= [\dot{\phi}^2 \partial_x(\sqrt{\det \bh}) \bh^{xx} + \dot{\phi}_0^2 \sqrt{\det \bh}\; \partial_x(\bh^{xx})]^{(n)} \\
&= \dot{\phi}_0^2 ((n+1) \overline{q}_{n+1} + \overline{Q} + (n+1) \overline{h}_{n+1}^{xx} + \overline{\mathfrak{h}}^{xx}) + F(\dot{\phi}_0, u_2)
\end{align*}
using
\begin{align*}
\overline{h}^{xx} &= \overline{h}_0^{xx}(s) + \overline{h}_2^{xx}(s) x^2 + \dots + \overline{h}_{n+1}^{xx}(s) x^{n+1} + \overline{\mathfrak{h}}^{xx}(s) x^{n+1} \log(x) + \dots \\
\sqrt{ \det \bh}(x,s) & = \overline{q}(s) + \dots + \overline{q}_{n+1}(s) x^{n+1} + \overline{Q}(s) x^{n+1} \log(x) + \dots
\end{align*}
Combining the two lower order polynomials and noting $\overline{q}_0 = 1 = \overline{h}_0^{xx}$, we find
\[
\boxed{ I_3 = \frac{1}{2}\int_{\gamma} \left( 2 (2-n) \dot{\phi}_0 \dot{\phi}_{n+1} + \dot{\phi}_0^2 [3(\overline{q}_{n+1} + \overline{h}_{n+1}^{xx}) + \overline{Q} + \overline{\mathfrak{h}}^{xx}] ) + F(\dot{\phi}_0, u_2(s))\right) dA_{\gamma}}
\]

\subsubsection{$I_4$, odd}
Again, $n+1$ is even so there will be many lower order terms. Thus we decompose the integrand into the principal part with $(n+1)$st order terms and the remainder
\[
I_4 = \frac{1}{2} \int_{\gamma} - (||\dot{S}||^2 ||\nabla x||^2 \sqrt{\det \bh})^{n+1} = -\frac{1}{2} \int_{\gamma} (x^{-2} \dot{\phi}^2 ||\nabla x||^2 \sqrt{\det \bh})^{n+1}
\]
Write this again as 
\[
\boxed{I_4 = -\frac{1}{2} \int_{\gamma} \left( \dot{\phi}^2_0[\overline{h}^{xx}_{n+1} + \overline{q}_{n+1}] + 2 \dot{\phi}_0 \dot{\phi}_{n+1} + F(\dot{\phi}_0, u_2(s)) \right) dA_{\gamma}}
\]

\subsubsection{$I_5$, odd}
We compute
\[
I_5 = - \frac{1}{2} \int_{\gamma} (||\dot{S}||^2 ||\nabla x||^2 \sqrt{\det \bh})^{(n+1, \log)} = - \frac{1}{2} \int_{\gamma} (\dot{\phi}^2 x^{-2} ||\nabla x||^2 \sqrt{\det \bh})^{(n+1, \log)} 
\]
Extracting the $x^{n+1} \log(x)$ terms is straightforward similar to the previous sections,
\[
\boxed{I_5 = - \frac{1}{2} \int_{\gamma} \dot{\phi}^2_0 [ \overline{\mathfrak{h}}^{xx} + \overline{Q}] + 2 \dot{\phi}_0 \dot{\Phi}  }
\]
having noted that 
\[
[x^{-2}||\nabla x||^2]^{(n +1, \log)} = [\bh^{xx}]^{(n+1, \log)} =: \overline{\mathfrak{h}}^{xx}
\]
\subsubsection{$I_6$, odd}
In parallel with the even case, we have
\begin{align*}
[dx(\ddot{S}) \sqrt{\det \bh}]^n &= [\ddot{\phi} dx(\nu) \sqrt{ \det \bh}]^n\\
&= - (n+1) \ddot{\phi}_0 u_{n+1} + F(\ddot{\phi}_0, u_2)
\end{align*}
So that 
\[
\boxed{I_6 = - (n+1) \int_{\gamma} \ddot{\phi}_0u_{n+1} + F(\ddot{\phi}_0, u_2) }
\]

\subsubsection{The Full Expression, Odd Case}
In summary, we proved that
\begin{align*}
D^2\mathcal{V}(Y) &= \int_{\gamma} -(n+1) \ddot{\phi}_0u_{n+1} + \dot{\phi}_0^2 \left[  (\overline{q}_{n+1} + \overline{h}_{n+1}^{xx}) - 4(n+2) u_2 U - 4(n+1) u_2 u_{n+1} \right] \\
& + \int_{\gamma} \dot{\phi}_0\left[(1 - n) \dot{\phi}_{n+1} - \dot{\Phi}\right] + F(\dot{\phi}_0, u_2)  dA_{\gamma}
\end{align*}
As with the even case, we compute in \S \ref{SecondVarComp}
\[
\overline{h}_{n+1}^{xx} + \overline{q}_{n+1} = - 4(n-1)u_2 u_{n+1} + \text{Tr}_{T\gamma}(k_{n+1,0}) + F(\dot{\phi}_0, u_2)
\]
for some function $F$. This time,  $\text{Tr}_{T\gamma}(k_{n+1,0}) \not \equiv 0$, since $n + 1$ is even and we include it in the formula. This yields
\begin{empheq}[box=\fbox]{align*}
D^2 \mathcal{V}(Y) &= \int_{\gamma} -(n+1) \ddot{\phi}_0 u_{n+1} + (1-n) \dot{\phi}_0 \dot{\phi}_{n+1} + \dot{\phi}_0^2 \left[-8n u_2 u_{n+1} + \text{Tr}_{T\gamma}(k_{n+1,0}) \right] \\
& - \dot{\phi}_0\left[ 4(n+2) \dot{\phi}_0 u_2 U + \dot{\Phi}\right] +  F(\ddot{\phi}_0, \dot{\phi}_0, u_2)
\end{empheq}
\noindent We've written things more suggestively to reflect the parallels with the even dimensional formula. This finishes the proof of \ref{codim1SecondVar} in the odd case. \qed \nl \nl
As an application of our second variation formula, we consider an even minimal submanifold, $Y$, flowed by an isometry to produce $\{Y_t\}$. Such a family has constant renormalized volume so $D^2 \RV(Y) = 0$.
\section{Application: Variation via Killing Vectors in $\H^{n+1}$}
\label{Killing}
In this section, we let $M = \H^{n+1}$ with 
\[
g = \frac{dx^2 + (dy_1^2 + \dots + dy_n^2)}{x^2}
\]
In particular $K = 0$ and $k_{n+1} = 0$ in \eqref{kEquation}, so we can make the corresponding simplifications to the first and second variational formulae. Consider the killing vector fields, $\{\partial_{y_i}\}$, applied to these formula - we prove
\begin{proposition} \label{EvenL2Result}
For $n$ even, $Y^n \subseteq \H^{n+1}$ minimal with closed boundary, $\partial Y = \gamma$, and graphical expansion given by $u(s,x)$ as in theorem \ref{AsymptoticExpansion}, we have
\begin{align*}
\langle u_2, u_{n+1} \rangle_{L^2(\gamma)} &= 0
\end{align*}
\end{proposition}
\noindent For $n$ odd, we have
\begin{proposition}
\label{OddL2Result}
For $n$ odd, $Y^n \subseteq \H^{n+1}$ minimal with closed boundary $\partial Y = \gamma$, and graphical expansion given by $u(s,x)$ as in theorem \ref{AsymptoticExpansion}, we have that 
\[
\langle u_2, 2(n^2 - 4n - 1) u_{n+1} - 8  U \rangle_{L^2(\gamma)} = B(u_2)
\]
where $B(u_2)$ denotes a boundary integral over $\gamma$ with integrand determined by $u_2$ and its derivatives.
\end{proposition}
\subsection{Proof: Codimension $1$, even}
We reference the expansion for the normal vector in \S \ref{Normal} and take
\begin{align*}
S_t := \dot{\phi}_k t \bnu \\
\dot{\phi}_k := \overline{\langle \partial_{y_k}, \bnu \rangle} 
\end{align*}
where $\overline{\langle \cdot, \cdot \rangle}$ denotes the inner product on $\H^{n+1}$ with respect to the compactified metric. Here, $\{y_k\}$ are the directions which are not $x$ in the metric expansion $g = \frac{dx^2 + dy_1^2 + \dots + dy_n^2}{x^2}$, and hence $\partial_{y_k}$ is a killing vector. Using, $\bg = g_{euc}$ we compute:
\begin{align*}
\bnu &= \frac{1}{\sqrt{1 - ||\overline{\nabla} u||^2 + || \overline{\nabla} u||^4}} \left[ (1 - ||\overline{\nabla} u||^2) \partial_z - \overline{\nabla} u \right] \\
&= K(x,s) \left( [1 - \bh^{ab} u_a u_b - 2 u_a u_x \bh^{ax} - u_x^2 \bh^{xx}] N - [u_a \bh^{ab} - u_x \bh^{xb}] \partial_{s_b} + [-u_a \bh^{ax} - u_x \bh^{xx}] \partial_x \right) \\
&= K(x,s) \left( \sum_{b = 0}^{n-1} d_{b}(x,s) \partial_{s_b} + D(x,s) N - (u_x + R_x) \partial_x \right)
\end{align*}
where $\partial_{s_a}$ is identified with $F_*(\partial_{s_a})$ via an abuse of notation. From this decomposition, we see that 
\begin{align*}
\F(K) &= 1, \qquad & K &= 1 + O(x^2)\\
\F(d_{b}) &= 1, \qquad & d_{b} &= O(x^2) \\
\F(D) &= 1, \qquad & D &= 1 + O(x^2) \\
\F(x^{-2} R_x) &= -1, \qquad & R_x &= 2 u_2 x + O(x^3)
\end{align*}
We now compute
\[
\dot{\phi}_k = \left[\sum_{b = 0}^{n-1} d_{b}(s,x) \overline{\langle\partial_{y_k}, \partial_{s_b}\rangle} + D(x,s) \overline{\langle \partial_{y_k}, N \rangle} \right] K(x,s)
\]
From which:
\begin{align*}
(\dot{\phi}_k)_0 &= \overline{\langle \partial_{y_k}, N \rangle} D_0 = \overline{\langle \partial_{y_k}, N \rangle} \\
(\dot{\phi}_k)_{n+1} &= \sum_{b = 0}^{n-1} [d_{b}(s,x)]_{n+1} [\overline{\langle\partial_{y_k}, \partial_{s_b}\rangle}]_0 + [D(x,s)]_{n+1} [\overline{\langle \partial_{y_k}, N \rangle}]_0 + [\overline{\langle \partial_{y_k}, N \rangle}]_0 [K(x,s)]_{n+1}
\end{align*}
We compute
\begin{align*}
[d_{b}]_{n+1} &= [u_a \bh^{ab} - u_x \bh^{xb}]_{n+1} \\
&= [u_b]_{n+1} \\
&= \partial_{s_b} u_{n+1}
\end{align*}
Similarly
\begin{align*}
[D]_{n+1} &= [1 - \bh^{ab} u_a u_b - 2 u_a u_x \bh^{ax} - u_x^2 \bh^{xx}]_{n+1} \\
&= - [u_x^2 \bh^{xx}]_{n+1} \\
&= - 4 (n+1) u_2 u_{n+1}
\end{align*}
And finally
\begin{align*}
K(x,s) &= \frac{1}{\sqrt{1 - ||\overline{\nabla} u||^2 + || \overline{\nabla} u||^4}} \\
\implies [K]_{n+1} &= -\frac{1}{2} [1 - ||\overline{\nabla} u||^2 + || \overline{\nabla} u||^4]_{n+1} \\
&= \frac{1}{2} [||\overline{\nabla} u||^2]_{n+1} \\
&=  2 (n+1) u_2 u_{n+1}
\end{align*}
This tells us 
\begin{align*}
(\dot{\phi}_k)_{n+1} &= \sum_{b = 0}^{n-1} [d_{b}(s,x)]_{n+1} [\overline{\langle\partial_{y_k}, \partial_{s_b}\rangle}]_0 + [D(x,s)]_{n+1} [\overline{\langle \partial_{y_k}, N \rangle}]_0 + [\overline{\langle \partial_{y_k}, N \rangle}]_0 [K(x,s)]_{n+1} \\
&= \sum_{b = 0}^{n-1} (\partial_{s_{b}} u_{n+1}) \overline{\langle\partial_{y_k}, \partial_{s_b}\rangle} - 2 (n+1) u_2 u_{n+1} \overline{\langle \partial_{y_k}, N \rangle}
\end{align*}
Finally, we note that $\ddot{\phi} = 0$, so we have
\begin{align*}
D^2 \mathcal{V}(\dot{S}_k, \dot{S}_k) \Big|_Y &= 0 \\
D^2 \mathcal{V}(\dot{S}_k, \dot{S}_k) \Big|_Y & = \int_{\gamma} (1 - n) (\dot{\phi}_k)_0 (\dot{\phi}_k)_{n+1} -8n (\dot{\phi}_k)_0^2 u_2 u_{n+1} \\
&= \int_{\gamma} (1 - n) \overline{\langle \partial_{y_k}, N \rangle} \left[\sum_{\ell = 0}^{n-1} \overline{\langle\partial_{y_k}, (\partial_{s_{\ell}} u_{n+1})\partial_{s_\ell}\rangle} - 2(n+1) \overline{\langle \partial_{y_k}, N \rangle}  u_{n+1} u_2  \right] \\
&  \qquad - 8n \overline{\langle \partial_{y_k}, N \rangle}^2 u_2 u_{n+1}
\end{align*}
Summing over $k = 1, \dots, n$, we note that
\begin{align*}
\sum_{k = 1}^n \overline{\langle \p_{y_k}, N \rangle} \sum_{\ell = 0}^{n} \overline{\langle \p_{y_k}, (\p_{s_{\ell}} u_{n+1}) \p_{s_{\ell}} \rangle} & = \overline{\langle N, \overline{\nabla}^{\gamma} u_{n+1} \rangle} \\
&= 0 \\
\sum_{k = 1}^n \overline{\langle \p_{y_k}, N \rangle}^2 &= 1
\end{align*}
when evaluated on $\gamma$. This is precisely because on $\gamma$, $\overline{\nabla}^{\gamma} u_{n+1} \in T\gamma$. Combining terms, we finally obtain
\begin{align*}
[-2(n+1)(1-n) - 8n] \int_{\gamma} u_2 u_{n+1} &= 0 \\
\implies \int_{\gamma} u_2 u_{n+1} &= 0
\end{align*}
having noted that $-2(n+1)(1-n) - 8n = 0$ has no integer solutions. \qed 
\subsection{Proof: Codimension $1$, odd}
The expression for $\dot{\phi}_k$ are the same in the odd dimensional case,
\[
\dot{\phi}_k = \left[\sum_{b = 0}^{n-1} d_{b}(s,x) \overline{\langle\partial_{y_k}, \partial_{s_{b}}\rangle} + D(x,s) \overline{\langle \partial_{y_k}, N \rangle} \right] K(x,s)
\]
and we compute 
\begin{align*}
(\dot{\phi}_k)_0 &= \left[\sum_{b = 0}^{n-1} (d_{b})_0 \overline{\langle \partial_{y_k}, \partial_{s_{b}} \rangle} + (D)_0 \overline{\langle \partial_{y_k}, N\rangle} \right] (K)_0 \\
&= \left[\sum_{b = 0}^{n-1} 0 \cdot \overline{\langle \partial_{y_k}, \partial_{s_{b}} \rangle} + 1 \cdot \overline{\langle \partial_{y_k}, N\rangle} \right] \cdot 1\\
&= \overline{\langle \partial_{y_k}, N \rangle}
\end{align*}
However
\begin{align*}
(\dot{\phi}_k)_{n+1} &= \left[\sum_{b = 0}^{n-1} (d_{b})_{n+1} \overline{\langle \partial_{y_k}, \partial_{s_{b}} \rangle} + (D)_{n+1} \overline{\langle \partial_{y_k}, N\rangle} \right] (K)_0 \\
& (D)_0 \overline{\langle \partial_{y_k}, N\rangle}  (K)_{n+1} + R(u_2) 
\end{align*}
where $R(u_2)$ denotes a polynomial in $\{u_2, \dots, u_{n-1}\}$ as with our convention for $F$. We compute 
\begin{align*}
(\dot{\phi}_k)_{n+1} &= \left[\sum_{b = 0}^{n-1} (\partial_{s_{b}} u_{n+1} + R(u_2)) \overline{\langle \partial_{y_k}, \partial_{s_{b}} \rangle} + D_{n+1}\overline{\langle \partial_{y_k}, N\rangle} \right] \\
& + \overline{\langle \partial_{y_k}, N\rangle} (K)_{n+1}  + R(u_2) \\
\end{align*}
recall the expansion of 
\begin{align*}
u &= u_2 x^2 + \dots + u_{n+1} x^{n+1} + U x^{n+1} \log(x) + O(x^{n+2})  \\
D(x,s) & = 1 - \bh^{ab} u_a u_b - 2u_a u_x \bh^{ax} - u_x^2 \bh^{xx}
\end{align*}
Further noting that $\F(\bh^{ab}) = \F(\bh^{xx}) = 1$ and $\bh^{ax} = O(x^3)$, we have
\begin{align*}
[D]_{n+1} &=  -4 u_2 [(n+1)u_{n+1} + U] + R(u_2) 
\end{align*}
Similarly,
\begin{align*}
K(x,s) &= \frac{1}{\sqrt{1 - ||\overline{\nabla} u||^2 + || \overline{\nabla} u||^4}} \\
\implies [K]_{n+1} & = \frac{1}{2} [||\overline{\nabla} u||^2 = 2 (n+1) u_2 u_{n+1} + R(u_2)
\end{align*}
in the same way that was deduced for the even case. With this, we have
\begin{align*}
(\dot{\phi}_k)_{n+1} &= \left[\sum_{b = 0}^{n-1} (\partial_{s_{b}} u_{n+1}) \overline{\langle \partial_{y_k}, \partial_{s_{b}} \rangle} - 2 u_2 [(n+1)u_{n+1} + U]\overline{\langle \partial_{y_k}, N\rangle} \right] + R(u_2) 
\end{align*}
for some conglomerate lower order term $R(u_2)$. According to our second variation formula for odd dimension submanifolds \S \ref{MainThm}, we also need the $x^{n+1} \log(x)$ coefficient, $\dot{\Phi}$, of $\dot{\phi}_k$:
\begin{align*}
\dot{\Phi}_k &= (\dot{\phi}_k)_{n+1, \log} \\
&= \left[\sum_{b = 0}^{n-1} (d_{b})_{n+1, \log} \overline{\langle \partial_{y_k}, \partial_{s_{b}} \rangle} + (D)_{n+1, \log} \overline{\langle \partial_{y_k}, N\rangle} \right] (K)_0  + \left[\sum_{b = 0}^{n-1} (d_{b})_0 \overline{\langle \partial_{y_k}, \partial_{s_{b}} \rangle} + (D)_0 \overline{\langle \partial_{y_k}, N\rangle} \right] (K)_{n+1, \log} \\
&= \left[\sum_{b = 0}^{n-1} (d_{b})_{n+1, \log} \overline{\langle \partial_{y_k}, \partial_{s_{b}} \rangle} + (D)_{n+1, \log} \overline{\langle \partial_{y_k}, N\rangle} \right] (K)_0  + \overline{\langle \partial_{y_k}, N\rangle} (K)_{n+1, \log}
\end{align*}
We immediately see that 
\begin{align*}
(d_{b})_{n+1, \log} &= (\partial_{s_{b}} U)  \\
(D)_{n+1, \log} &= -4 (n+1) u_2 U \\
(K)_{n+1, \log} &= 2 (n+1) u_2 U
\end{align*}
so that
\begin{align*}
\dot{\Phi}_k &= \left[\sum_{b = 0}^{n-1} U_{b} \overline{\langle \partial_{y_k}, \partial_{s_{b}} \rangle} - 2 (n+1) u_2 U\overline{\langle \partial_{y_k}, N\rangle} \right]
\end{align*}
Now using the formula for the second variation of renormalized volume (with $k_{n+1} = 0$ in $\H^{n+1}$), we have
\begin{align*}
0 &= D^2 \V(Y)(\dot{\phi}_k, \dot{\phi}_k) \\
& = \int_{\gamma} (1 - n) (\dot{\phi}_k)_0 (\dot{\phi}_k)_{n+1} + \dot{\phi}_0^2 [-8n u_2 u_{n+1}] - (\dot{\phi}_k)_0 \left[ 4(n+2) (\dot{\phi}_k)_0 u_2 U + \dot{\Phi}_k\right] + R(u_2)
\end{align*}
Summing over $k$ and combining terms, we get
\[
\begin{gathered}
0 = \int_{\gamma} 2(n^2 - 4n - 1) u_2 u_{n+1}  -8 u_2 U + R(u_2) \\
\end{gathered}
\]
which parallels the even case up to an error term $R(u_2)$ and the presence of $U$. This proves the theorem. \qed

\section{Appendix} 

\subsection{Metric on $TM$}
\label{MetricTM}
We construct a frame for all of $TM = T\Gamma \oplus N \Gamma$ using Fermi coordinates on $\Gamma$.
%
Coordinatize our space as follows: Let $p \in \gamma \subseteq \partial M$ be labeled by geodesic normal coordinates on $\gamma$ about some base point $p_0$, i.e.
\[
p = f(s) := \overline{\exp}_{p_0} (s^a E_a)
\]
for $\{E_a\}$ an ONB at $p_0$ spanning $T \gamma$. We then coordinatize $\Gamma$ as points $(s,x) \leftrightarrow (f(s), x) \in \gamma \times [0, \eps)$. Then for $U(\Gamma) \subseteq M$ sufficiently small, we define
\[
\begin{gathered}
W: B_1(0)^{m-1} \times [0, \eps) \times B_1(0)^{n+1-m} \to M \\
W(s, x, z):= \overline{\exp}_{(f(s), x)}(z^i X_i)
\end{gathered}
\]
for $\{X_i\}$ an ONB for $N(\Gamma)$. Both exponential maps are taken with respect to $\bg$ restricted to $\gamma$ and $\Gamma$ respectively. Abusing notation slightly, we define
\begin{align} \label{MetricCoeffs} 
\overline{g}_{ij} &:=  \overline{g}(W_*(\partial_{z_i}), W_*(\partial_{z_j})) =  \delta_{ij} + [\overline{\overline{\Gamma}}_{ai}^j + \overline{\overline{\Gamma}}_{aj}^i] s^a  + [\overline{\overline{\Gamma}}_{xi}^j + \overline{\overline{\Gamma}}_{xj}^i ]x + O(z_i z_j, s_a s_b, s_a x, z_i x, z_i s_a, x^2) \\ \nonumber
\overline{g}_{ab} &:= \overline{g}(W_*(\partial_{s_a}), W_*(\partial_{s_b})) = \delta_{ab} + [\overline{\overline{\Gamma}}_{ka}^b + \overline{\overline{\Gamma}}_{kb}^a] z^k + [\overline{\overline{\Gamma}}^a_{bx} + \overline{\overline{\Gamma}}^b_{ax}] x + O(z_i z_j, s_a s_b, s_a x, z_i x, z_i s_a, x^2)  \\ \nonumber
\overline{g}_{ai} &:= \overline{g}(W_*(\partial_{s_a}), W_*(\partial_{z_i})) =  [\overline{\overline{\Gamma}}_{ca}^i + \overline{\overline{\Gamma}}_{ci}^a] s^c + [\overline{\overline{\Gamma}}_{cx}^i + \overline{\overline{\Gamma}}_{ci}^x]x + O(z_i z_j, s_a s_b, s_a x, z_i x, z_i s_a, x^2) \\ \nonumber
\overline{g}_{ax} &:= \overline{g}(W_*(\partial_{s_a}), W_*(\partial_x))  = 0\\ \nonumber
\overline{g}_{ix} &:= \overline{g}(W_*(\partial_{z_i}), W_*(\partial_{x})) = 0 \\ \nonumber
\overline{g}_{xx} &:= \overline{g}(W_*(\partial_x), W_*(\partial_x)) = 1 \nonumber
\end{align}
where $\overline{\overline{\Gamma}}_{\cdot\cdot}^{\cdot}$ are the Christoffel symbols for $\gamma \times [0,\eps) \subseteq \partial M \times [0, \eps)$ equipped with $dx^2 + k_0(s,z)$, i.e. $\bg$ evaluated to lowest order in $x$. For the first $3$ expansions, see \cite{mahmoudi2004constan} among other sources. The last $3$ equations follow because the metric splits along the $x$ direction:
\[
\overline{g} = dx^2 + k = dx^2 + [k_0 + O(x^2, s^2, z^2, sx, sz, xz)]
\]
i.e. the metric is block diagonal with a $1 \times 1$. Recall the index notation
\begin{equation*}
\begin{gathered}
a,b,c,d \leftrightarrow s_a, s_b, s_c, s_d \\
i,j,k, \ell \leftrightarrow z_i, z_j, z_k, z_{\ell} \\
i,j,k,\ell \leftrightarrow w_i, w_j, w_k, w_{\ell} \\
\alpha, \beta, \gamma, \delta \leftrightarrow \{y_{\alpha}, y_{\beta}, y_{\gamma}, y_{\delta} \} \subseteq \{s_a, x\} \\
\sigma, \mu, \nu, \tau, \omega \leftrightarrow \{y_{\sigma}, y_{\mu}, y_{\nu}, y_{\tau}, y_{\omega}\} \subseteq \{s_a, x, z_i\} 
\end{gathered}
\end{equation*}
We will also often conflate $\partial_{s_a}, \partial_x, \partial_{z_i}$ with their pushforwards by $W$ as well. Given our asymptotics for $\overline{g}_{\mu \nu}$ in terms of $s$, $z$, and $x$, we can evaluate at $z = u(s,x)$ to derive asymptotics for a frame for $TY$.

\subsection{Metric on $TY$}
\label{MetricTY}
We construct a frame for $TY$ and derive an expansion for the metric, $\overline{g} \Big|_Y$, in this frame. Recall the map
\[
\begin{gathered}
G: \Gamma \to Y \embed M  \\
G(s, x) = W(s, z = u(s,x), x) = (F(s, u(s,x), x)
\end{gathered}
\]
We consider the frame for $TY$ given by 
\begin{align*}
v_a &:= G_*(\partial_{s_a}) = \partial_{s_a} + u_a^i \partial_{z_i}\\
v_x &= G_*(\partial_x) = \partial_x + u_x^i \partial_{z_i}
\end{align*}
%
Where we notationally identified $\partial_{s_a}$, $\partial_x$, and $\partial_{z_i}$ with their pushforwards by $W$. We will denote the above as 
\[
G_*(\partial_{y_{\alpha}}) = \partial_{y_{\alpha}} + u_{\alpha}^i \partial_{z_i}
\] 
The induced metric is then given by
\begin{align} \label{TYMetricAsymptotics} 
\overline{h}_{ab} & = \bg_{ab} + \bg_{aj} u_b^j  + u_a^i \bg_{bi} + u_a^i u_b^j \bg_{ij}  \\ \nonumber
\bh_{ax} &= u_x^k \bg_{ka} + u_a^i u_x^k \bg_{ik} \\ \nonumber
\overline{h}_{xx} &= 1 + u_x^k u_x^{\ell} \bg_{k \ell} 
\end{align}
using the metric notation of section \S \ref{SFFParity}. As a point of notation, we let $\va \in \{v_a, v_x\}$ so that $\{v_{\alpha}\}$ is a basis for $TY$, with $\alpha$ taking on the $x$ and $a$ subscripts. Now assume $m$ is even. Evaluating at $s = 0$ and $z = u(s,x)$ and using equation \ref{MetricCoeffs} and lemma \ref{ChristoffelParityCylinder} applied to the symbols $\{\bGamma_{\sigma \omega \tau}\}$ by converting from $g \to \bg$, we get that
\begin{align*} 
\overline{h}_{ab}\Big|_{(s = 0, z = u, x)} & = \delta_{ab} + O(x^2), \qquad &\F\left(\overline{h}_{ab} \Big|_{(s = 0, z = u, x)}\right) &= 1 \\
\overline{h}_{ax}\Big|_{(s = 0, z = u, x)} & = O(x^3), \qquad &\F\left(\overline{h}_{ax} \Big|_{(s = 0, z = u, x)}\right) &= -1 \\
\overline{h}_{xx}\Big|_{(s = 0, z = u, x)} & = \delta_{ij} + O(x^2), \qquad &\F\left(\bh_{xx} \Big|_{(s = 0, z = u, x)}\right) &= 1
\end{align*}

\subsubsection{The $\bT$ matrix}
\label{TMatrix}
We use the previous section to define a frame for $TM$ using the decomposition $TM = TY \oplus N \Gamma$, which holds at points $p \in M$ with $x < \eps$. For $U(\Gamma)$ and $U(Y)$ open neighborhoods around $\Gamma$ and $Y$ respectively, consider the map
\begin{align*}
R&: U(\Gamma) \to U(Y) \\
R(s, z, x) &= W(s, z + u(s,x), x)
\end{align*}
and define 
\[
\bT := R^*(\bg)
\]
so that 
\[
\bT = \begin{pmatrix}
\bT_{xx} & \bT_{xa} & \bT_{xj} \\
\bT_{ax} & \bT_{ab} & \bT_{aj} \\
\bT_{ix} & \bT_{ib} & \bT_{ij}
\end{pmatrix}
\]
where each entry is $\overline{g}(R_*(\partial_{\cdot}), R_*(\partial_{\cdot}))$. We recall the index notation (see equation \eqref{IndexNotation}) of 
\[
\{v_{\sigma}\} = \{v_a, v_x, v_i\} = \{R_*(\partial_{s_a}), R_*(\partial_x), R_*(\partial_{z_i})\}
\]
which is a frame for all of $TM$ for $x < \eps$ small.
We compute $T$ in these coordinates as $\bT_{\sigma\rho}:= \bg(v_{\sigma}, v_{\rho})$. Note that $\bT_{\alpha \beta}$ have been computed in the previous section \S \ref{MetricTY}. For the new entries, we have 
\begin{align*}
\bT_{ij} &= \bg_{ij} \\
\bT_{ia} &= \bg_{ai} + u_a^j \bg_{ij} \\
\bT_{ix} &= \bg_{xi} + u_x^k \bg_{ik}
\end{align*}
and immediately from \ref{MetricCoeffs} and lemma \ref{ChristoffelParityCylinder}, we get that
\begin{align*}
\bT_{ij}\Big|_{z = u} &= \delta_{ij} + O(x^2), \quad & \F(\bT_{ij}) &= 1 \\
\bT_{ia}\Big|_{z = u} &= O(x^4), \quad & \F(x^{-2}\bT_{ia}) &= 1 \\
\bT_{ix}\Big|_{z = u} &= O(x^3), \quad & \F(x^{-2}\bT_{ix}) &= -1
\end{align*}
we can also invert the metric and get the same asymptotics and $\F(\cdot)$ values.
\subsection{Projected basis for the normal bundle}
\label{NormalFrame}
We prove the following, again abusing notation by writing $\partial_{(\cdot)}$ for $W_*(\partial_{(\cdot)})$ where needed:
\begin{lemma}
For any $p \in \gamma$ and a neighborhood $U(p) \subseteq M$, there exists a frame $\{\overline{w}_1, \dots, \overline{w}_{n+1-m}\}$ for $N(U(p) \cap Y)$ which is orthonormal with respect to $\overline{g}$ on $\gamma$ such that for $m$ even (odd)
\begin{align*}
\overline{g}(w_i, \partial_x) &= O(x), \qquad &\F(\overline{g}(w_i, \partial_x)) &= -1 \\
\overline{g}(w_i, \partial_{s_a}) &= O(x^2), \qquad &\F(\overline{g}(w_i, \partial_{s_a})) &= 1 \\
\overline{g}(w_i, \partial_{z_j}) &= \delta_{ij} + O(x^2), \qquad &\F(\overline{g}(w_i, \partial_{z_j})) &= 1
\end{align*}
in fact, $\overline{g}(w_i, \partial_{s_a})$ is even up to $m+2$ ($m+3$) when $m$ is even (odd).
\end{lemma}
\noindent \Pf For notational brevity, we handle $m$ even, noting that all $m$ related indices will be shifted up by $1$ when $m$ is odd by \ref{AsymptoticExpansion}. Recall the frame for $TY = \text{span} \{v_a, v_x\}$ as given in \ref{MetricTY}. Now setting $\overline{N}^k := W_*(\partial_{z_k})$ for notation, we define
\begin{align*}
w_k & := \Pi_{N(Y)} \overline{N}^k = W_*(\partial_{z_k}) - \Pi_{TY}(W_*(\partial_{z_k})) \\
&= W_*(\partial_{z_k}) - \left[\overline{h}^{a b} \overline{g}(v_a, W_*(\partial_{z_k})) v_b + \overline{h}^{ax} (v_a, W_*(\partial_{z_k})) v_x \right] \\
& - \left[ \overline{h}^{xb} (v_x, W_*(\partial_{z_k})) v_b + \overline{h}^{xx} (v_x, W_*(\partial_{z_k})) v_x \right] 
\end{align*}
Now using \S \ref{MetricTY} and \eqref{MetricCoeffs}, we get
\begin{align} \label{NormalProjections}
\overline{g}(w_i, \partial_x) &= O(x), \qquad &\F(\overline{g}(w_i, \partial_x)) &= -1 \\ \nonumber
\overline{g}(w_i, \partial_{s_a}) &= O(x^2), \qquad &\F(\overline{g}(w_i, \partial_{s_a})) &= 1 \\ \nonumber
\overline{g}(w_i, \partial_{z_j}) &= \delta_{ij} + O(x^2), \qquad &\F(\overline{g}(w_i, \partial_{z_j})) &= 1 
\end{align} 
using the established parity of the metric coefficients. Note that the $\{w_k\}$ are not normalized but we compute using \eqref{NormalProjections} and \S \ref{TMatrix}
\begin{align*}
\overline{g}(w_k, w_k) &= c_k(x,s) \\
&= 1 + O(x^2) \\
\F(c_k) &= 1
\end{align*}
so we define 
\[
\overline{w}_k = \frac{w_k}{\sqrt{\bg(w_k, w_k)}}
\]
which still obeys \eqref{NormalProjections}.
\subsection{Simons Operator}
\label{Simons}
In this section, we compute the Simons operator for $Y$, the graph over $\Gamma$ given by $u(s,x)$ where $||u(s,x)|| = O(x^2)$ and is even in $x$ up to order $m$. Note that this includes $u = 0$, corresponding to the boundary cylinder $\Gamma = \gamma \times \R^+$. In particular, we show that the Simons operator is $O(x^2)$ in its leading coefficient. Recall our notation of $\{v_{\alpha}\} = \{v_a, v_x\}$ for a basis of $TY$. Then
\[
\tilde{A}(X) = g((\nabla_{v_{\alpha}} v_{\beta})^N, X) h^{\alpha \gamma} h^{\beta \delta} (\nabla_{v_{\gamma}} v_{\delta})^N
\]
We have
\[
X = X^j w_j
\]
where $\{w_j = w_j(x,s)\}$ is the basis for $N(Y)$ as in \ref{NormalFrame}. We compute
\begin{align*} 
\tilde{A}(X) &= X^{j}[ g((\nabla_{v_a} v_b)^N, w_j) h^{a\gamma} h^{b\delta} (\nabla_{v_{\gamma}} v_{\delta})^N \\
&+ g((\nabla_{v_{a}} v_x)^N, w_j) h^{a\gamma} h^{x\delta} (\nabla_{v_{\gamma}} v_{\delta})^N \\
&+ g((\nabla_{v_x} v_b)^N, w_j) h^{x\gamma} h^{b\delta} (\nabla_{v_{\gamma}} v_{\delta})^N \\
&+ g((\nabla_{v_x} v_x)^N, w_j) h^{x\gamma} h^{x\delta} (\nabla_{v_{\gamma}} v_{\delta})^N]
\end{align*}
and we expand
\[
(\nabla_{v_{\alpha}} v_{\beta})^N = (\tilde{\bGamma}_{\alpha \beta}^{\sigma} v_{\sigma})^N =  \tilde{\bGamma}_{\alpha \beta}^{j} w_j
\]
We reference the following Christoffel symbols (computed following \cite{graham1991einstein} and lemma \ref{YChristoffelParity})
\begin{align*}
\tilde{\bGamma}_{ab}^i &= O(1), \qquad & \F(\tilde{\bGamma}_{ab}^i) &= 1 \\
\tilde{\bGamma}_{ax}^i &= O(x) \qquad & \F(\tilde{\bGamma}_{ax}^i) &= -1 \\
\tilde{\bGamma}_{xb}^i &= O(x) \qquad & \F(\tilde{\bGamma}_{xb}^i) &= -1\\
\tilde{\bGamma}_{xx}^i &= O(x^2) \qquad & \F(\tilde{\bGamma}_{xx}^i) &= 1
\end{align*}
This immediately tells us that 
\begin{align*}
g((\nabla_{v_\alpha} v_{\beta})^N, w_j) &= O(x^{-2}) \\
\F\left( x^2 g((\nabla_{v_\alpha} v_{\beta})^N, w_j) \right) &= 1
\end{align*}
and so
\begin{align*}
\tilde{A}(X) &= X^j [g((\nabla_{v_{\alpha}} v_{\beta})^N, w_j) h^{\alpha \gamma} h^{\beta \delta} \tilde{\bGamma}_{\delta \gamma}^k) w_k \\
&= F^j(\{X^k\}) w_j = O(x^2) w_j
\end{align*}
i.e. $F^j$ is some linear function of $\{\dot{\phi}^k\}$ and 
\begin{align*}
F^j(\{X^k\}) &= f_{k}^j(s,x) X^k(s,x) \\
f_{k}^j(s,x) &= O(x^2) \\
\F(f_k^j) &= 1
\end{align*}
where the last line holds by parity of the Christoffels and metric coeffiicents. This shows that the Simons operator gives a quadratic error term that preserves parity.

\subsection{Equivalence of Hadamard and Riesz Regularization}
\label{equivOfDef}
In this section, we demonstrate that Hadamard regularization and Riesz regularization are equal

\subsubsection{$m$ even}
Under \textbf{Hadamard regularization}, the renormalized volume is the constant term in the following expansion
\[
V(Y \cap \{x > \epsilon\}) = \int_{Y \cap \{x > \epsilon\}} dA = a_0 \epsilon^{-m + 1} + a_2 \epsilon^{-m + 3} + \dots + a_{m-2} \epsilon^{-1} + a_m + O(\epsilon \log(\eps))
\]
Such an expansion follows because
\begin{align*}
\int_{Y \cap \{x > \epsilon\}} dA_{Y} &= \int_{\gamma} \int_{\epsilon < x < b} \frac{\sqrt{\det{\overline{h}}}}{x^m} d V_{\gamma} dx + \int_{x > b} dA_{Y}\\
&= \int_{\gamma} \int_{\epsilon < x < b} \left(\overline{h}_0(\vec{s}) x^{-m} + (\text{even terms}) + \overline{h}_{m-2}(\vec{s})x^{-2} + \overline{R}_m(s,x) ) dx \right) dV_{\gamma} + \int_{x > b} dA_{Y} \\
&= \sum_{k =0}^{m/2-1}\epsilon^{2k-m + 1} c_{2k} + A_m + O(\epsilon \log(\eps))
\end{align*}
where we have used that the expansion of the volume form is even up until $\overline{h}_m$ (cf corollary \ref{NoLogTerms}). In addition, $\overline{R}_m(s,x) = O(\log(x))$ is the remainder term. We set
\begin{align} \nonumber
c_{2k} &:= \frac{1}{m - 2k - 1} \int_{\gamma} \bh_{2k}(s) dA_{\gamma}(s) \qquad 2k \leq m \\ \label{HadamardExpression}
A_m &:= \int_{\gamma}\int_{0 \leq x \leq b}\overline{R}_m(s,x) dx dV_{\gamma} - \sum_{k = 0}^{m/2 - 1} c_{2k} b^{-m + 2k + 1} + \int_{x > b} dA_{Y} 
\end{align}
Note that $A_m$ is finite (recall that $\{x \geq b\}$ is a compact region) and actually independent of $b$. We denote
\[
{}^H\! \int_Y dA := A_m
\]
Under \textbf{Riesz regularization}, consider the meromorphic function 
\begin{align*}
\zeta(z) &:= \int_{Y} x^z dA_Y \\
&= \frac{a_0}{z - (-m + 1)} + \dots + \frac{a_{m-1}}{z} + D_m + O(z) 
\end{align*}
We then define
\[
{}^R\!\int_Y dA := \FPz \zeta(z)
\]
we compute this as 
\begin{align*}
\int_{Y}  x^z dA_Y &= \int_{Y \cap \{x \leq b\}} x^z dA_Y + \int_{Y \cap \{x > b\}} x^z dA_Y\\
&= \int_{\gamma} \int_{0 \leq x \leq b} x^{z-m} \sqrt{\det{\overline{h}}} d V_{\gamma} \wedge dx + \int_{Y \cap \{ x > b\}} x^z dA_Y \\ 
&= \int_{\gamma} \int_{0 \leq x \leq b} \left(\overline{h}_0(\vec{s}) x^{z-m} + (\text{even terms}) + \overline{h}_{m-2}(\vec{s})x^{z-2} + \overline{R}_m(s,x) x^z) dx \right) dA_{\gamma}  + \int_{Y \cap \{x > b\} } x^z dA_Y
\end{align*}
where $b \ll 1$ so that we can use such an expansion. For the first integral, we write 
\[
\int_{\gamma} \int_{0 \leq x \leq b} \left(\overline{h}_0(\vec{s}) x^{z-m} + (\text{even terms}) + \overline{h}_{m-2}(\vec{s})x^{z-2}) dx \right) dV_{\gamma} \wedge dx = \sum_{k = 0}^{m/2-1} a_{2k} \frac{b^{z - m + 2k + 1}}{z - m + 2k + 1}
\]
having assumed $\text{Re}(z) \gg 0$ and setting
\[
a_{2k}:= \int_{\gamma} \overline{h}_{2k}(s) dA_{\gamma}(s) = (m - 2k - 1) c_{2k}
\]
Now we define 
\begin{align*}
D_m(b,z) &:= \int_{\gamma} \int_{0 \leq x \leq b} \overline{R}_{m}(s,x)x^z  + \sum_{k = 0}^{m/2-1} \frac{b^{z-m+2k+1} - 1}{z - m + 2k + 1} a_{2k} + \int_{Y \cap \{x \geq b\} } x^z dA_Y \\
\implies \int_{Y} x^z dA_Y &= \sum_{k = 0}^{m/2-1} \frac{a_{2k}}{z - m + 2k + 1} + D_m(b,z) 
\end{align*}
As before, one can compute that $D_{m}(b,z)$ is finite, independent of $b$, and holomorphic at $z = 0$. We can compute the finite part at $z = 0$ as 
\begin{align*}
\FPz \int_Y x^z dA_Y &= \sum_{k = 0}^{m/2 - 1} \frac{a_{2k}}{-m + 2k + 1} + D_m(b,0) \\
& = \sum_{k = 0}^{m/2-1} \frac{b^{-m+2k+1}}{-m + 2k + 1} a_{2k} + \int_{\gamma} \int_{0 \leq x \leq b} \overline{R}_m(s,x) dx dA_{\gamma}(s) + \int_{Y \cap \{x \geq b\} } dA_Y \\
&= -\sum_{k = 0}^{m/2-1} b^{-m+2k+1}c_{2k} + \int_{\gamma} \int_{0 \leq x \leq b} \overline{R}_m(s,x) dx dA_{\gamma}(s) + \int_{Y \cap \{x \geq b\} } dA_Y
\end{align*}
this is the same as \eqref{HadamardExpression} so that 
\[
\boxed{{}^H\!\int_Y dA = {}^R\!\int_Y dA}
\]
\subsubsection{$m$ odd}
Under Hadamard regularization, the renormalized volume is still the constant coefficient in the following expansion
\[
V(Y \cap \{x > \epsilon\}) = \int_{Y \cap \{x > \epsilon\}} dA_Y = a_0 \epsilon^{-m + 1} + a_2 \epsilon^{-m + 3} + \dots + a_{m-2} \epsilon^{-2} + \tilde{a} \log(\epsilon^{-1}) +  a_m  + o(1)
\]
Such an expansion follows because
\begin{align*}
\int_{Y \cap \{x > \epsilon\}} dA_{Y} &= \int_{\gamma} \int_{\epsilon < x < b} \left(\overline{h}_0(\vec{s}) x^{-m} + (\text{even terms}) + \overline{h}_{m-1}(\vec{s})x^{-1} + \overline{R}_m(s,x) ) dx \right) dV_{\gamma} + \int_{x > b} dA_{Y} \\
&= \sum_{k =0}^{(m-1)/2 - 1}\epsilon^{2k-m + 1} c_{2k} + c_{m-1}\log(\epsilon^{-1}) + A_m + O(\epsilon)
\end{align*}
where we used that for $m$ odd, the expansion of the volume form is even up until $\overline{h}_{m+1}$, at which point the $x^{m+1}\log(x)$ term also appears. Denote
\begin{align*}
c_{2k} &:= \frac{1}{m - 2k - 1} \int_{\gamma} \bh_{2k}(s) dA_{\gamma}(s), \qquad \qquad 2k \leq m - 1 \\
A_m &:= \int_{\gamma}\int_{0 \leq x \leq b}\overline{R}_m(s,x) dx dV_{\gamma} - \sum_{k = 0}^{(m-1)/2 - 1} c_{2k} b^{-m + 2k + 1} + c_{m-1} \log(b) +  \int_{x > b} dA_{Y} \\
{}^H \! \int_Y dA &:= A_m
\end{align*}
analogous to the even case. For Riesz regularization, we have 
\[
{}^R \! \int_Y dA := \FPz \zeta(z) = \FPz\int_Y x^z dA
\]
and we want to show that this gives $A_m$ as above. Using the same expansion, we get 
\begin{align*}
\int_{Y }  x^z dA_{Y} &= \int_{\gamma} \int_{0 < x < b} x^z \frac{\sqrt{\det{\overline{g}}}}{x^m} d V_{\gamma} dx + \int_{x > b}  x^z dA_{Y}\\
&= \int_{\gamma} \int_{0 < x < b} \left(\overline{h}_0(\vec{s}) x^{z-m} + (\text{even terms}) + \overline{h}_{m-1}(\vec{s})x^{z-1} + x^z\overline{R}_m(s,x) ) dx \right) dV_{\gamma} + \int_{x > b} x^z dA_{Y} \\
&= \sum_{k =0}^{(m-1)/2}\frac{b^{z - m + 2k + 1}}{z - m + 2k + 1} a_{2k}  + \int_{\gamma} \int_{0 < x < b} x^z \overline{R}_{m}(s,x) dx dA_{\gamma} + \int_{x > b} x^z dA_Y
\end{align*}
for 
\[
a_{2k} = \int_{\gamma} \bh_{2k} dA_{\gamma} = (m - 2k - 1) c_{2k}
\]
as before. Then the same analysis as before gives us 
\[
\FPz \left[ \int_{\gamma} \int_{0 < x < b} x^z \overline{R}_{m}(s,x) dx dA_{\gamma} + \int_{x > b} x^z dA_Y \right] = \int_{\gamma} \int_{0 < x < b} \overline{R}_{m}(s,x) dx dA_{\gamma} + \int_{x > b} dA_Y
\]
and also 
\begin{align*}
\sum_{k =0}^{(m-1)/2}\frac{b^{z - m + 2k + 1}}{z - m + 2k + 1} a_{2k} &= \frac{b^{z-m + 1}}{z - m + 1} a_{0} + \dots + \frac{b^{z - 2}}{z - 2} a_{m - 3} + \frac{b^z}{z} a_{m-1} \\
\implies \FPz \sum_{k =0}^{(m-1)/2}\frac{b^{z - m + 2k + 1}}{z - m + 2k + 1} a_{2k} &= \frac{b^{-m + 1}}{- m + 1} a_{0} + \dots + \frac{b^{- 2}}{- 2} a_{m - 3} + \FPz\left(\frac{b^z}{z}\right) a_{m-1}
\end{align*}
We now expand
\[
\begin{gathered}
\frac{b^z}{z} = \frac{1}{z} + \log(b) + O(z) \\
\implies \FPz \left( \frac{b^z}{z} \right) = \log(b)
\end{gathered}
\]
and so 
\[
\FPz \int_{Y} x^z dA_Y = \int_{\gamma} \int_{0 \leq x \leq b} \overline{R}_{m}(s,x) dx dA_{\gamma} - \sum_{k = 0}^{(m-1)/2 - 1} c_{2k} b^{-m + 2k + 1} + c_{m-1} \log(b) + \int_{x > b} dA_Y
\]
thus 
\[
{}^R \! \int_Y dA = {}^H \! \int_Y dA
\]
in the odd case. It is interesting to note that $\zeta(z)$ \textit{does} have a pole at $z = 0$ in the odd case, in contrast with the even case. Note that for $m$ odd, \textit{renormalized volume depends on the choice of special bdf} - see \cite{albin2009renormalizing} for details.
\subsection{Degeneracy of Minimal Hypersurfaces in $M^{n+1}$}
\label{Degeneracy}
In this section, we summarize the relevant results from section $4$ of \cite{alexakis2010renormalized} on the degeneracy of a minimal submanifold $Y^m \subseteq M$. For $X \in N(Y)$, recall the Jacobi operator:
\[
J_Y(X) = \Delta^{\perp}_Y(X) + \tilde{A}(X) - \text{Tr}_{TY}[R(\cdot, X) \cdot]
\]
We can view $J_Y$ as a map between weighted H\"older spaces
\[
J_Y: x^{\mu} \Lambda^{2,\alpha}_0(N(Y)) \to x^{\mu} \Lambda^{0,\alpha}_0(N(Y))
\]
Moreover, the indicial roots of this operator are $\mu_1 = -1$ and $\mu_2 = m$. When $-1 < \mu < m$, we know from \cite{rafe1991elliptic} that $J_Y$ is Fredholm and index $0$ in the codimension $1$ case. A minimal submanifold $Y^m \subseteq M$ is said to be \textit{nondegenerate} if the kernel of this map is just $\{0\}$. Again in the codimension $1$ case, one can show that any boundary variation of a nondegenerate submanifold, i.e.
\[
\gamma_{\psi} = \{\overline{\exp}_p (\psi(p) N(p)) \}, \qquad \psi \in C^{\infty}(\gamma)
\]
can be extended to a Jacobi field on $Y$, $X = \dot{\phi}(p,x) \nu(p,x)$, so that $\dot{\phi}(p,0) = \psi(p)$. This follows by an inverse function theorem argument which is a direct adaption of the $n = 2$ case described in \cite{rafe1991elliptic}.  \nl \nl
When $Y$ is degenerate, $L_Y$ is still index $0$ but the kernel is a non-trivial finite dimensional space. We consider the kernel of $L_Y$ acting on functions $\psi \in L^2$, i.e. such that $\psi_0 = 0$,
\[
K = \{\psi \; | \; J_Y(\psi) = 0, \; \psi_0 = 0\}
\]
Elements of $K$ have the following asymptotic expansion:
\[
\psi \sim \begin{cases} \psi_{m+1}(s) x^{m+1} + \dots & m \text{ even}\\
\psi_{m+1}(s) x^{m+1} + \Psi(s) x^{m+1} \log(x) + \dots & m \text{ odd}
\end{cases}
\]
i.e. they vanish to order $m+1$. Consider 
\[
V = \{f \; | \; \exists \psi \in K \; \st \; f = \psi_{n+1} \}
\]
which is also finite dimensional. Recall equation \eqref{codim1EvenVariation}, i.e. for $Y^n$ a critical point of renormalized volume and $n$ even, we have 
\[
\int_{\gamma} u_{n+1}(s) \dot{\phi}_0(s) dA_{\gamma} = 0
\]
From the above, if $\dot{\phi}_0(s) \leftrightarrow \phi(s,x)$ a Jacobi field, and $\psi_{n+1}(s) \leftrightarrow \psi(s,x)$ an $L^2$ Jacobi field, then 
\begin{align*}
0 &= \int_Y J_Y(\phi(s,x)) \psi(s,x) dA_Y  \\
& = \int_Y [(\Delta_Y \phi) \psi + |A_Y|^2 \phi \psi - \Ric_{TY}(\nu, \nu) \phi \psi] dA_Y \\
& = \int_Y [(\Delta_Y \psi) \phi + |A_Y|^2 \psi \phi - \Ric_{TY}(\nu, \nu) \psi \phi] dA_Y - \int_{\gamma} (n+1)(\psi_{n+1}) \dot{\phi}_0 dA_{\gamma} \\
& = \int_Y J_Y(\psi) \phi - (n+1)\int_{\gamma} \psi_{n+1} \dot{\phi}_0 dA_{\gamma} \\
&= -(n+1) \int \psi_{n+1} \dot{\phi}_0 dA_{\gamma}
\end{align*}
The equality in line three is most easily seen by switching to the compactified metric, keeping track of powers of $x$, and integrating by parts twice. This tells us that $\dot{\phi}_0(s)$ must be orthogonal to all elements of $V$, i.e. $\langle \dot{\phi}_0(s), f\rangle_{L^2(\gamma)} = 0$. This then tells us that any $u_{n+1} \in V$ lies in a finite dimensional space.
\subsection{Computing $\overline{h}_{n+1}^{xx}$ and $\overline{q}_{n+1}$}
\label{SecondVarComp}
In this section, we consider $Y^n \subseteq M^{n+1}$ and compute the $(n+1)$st coefficient for the metric coefficient $\overline{h}^{xx}(s,x)$ and the volume form prefactor $\overline{q}(s,x)$. We aim to show
\begin{proposition*}
For $Y^n \subseteq M^{n+1}$ minimal, we have
\[
\overline{h}^{xx}_{n+1} + \overline{q}_{n+1} = \begin{cases}-4(n-1) u_2 u_{n+1} & n \text{ even} \\
-4(n-1) u_2 u_{n+1} + \text{Tr}_{T\gamma}(k_{n+1,0}) + R(u_2) & n \text{ odd} 
\end{cases}
\]
where 
\[
k_{n+1,0} := \frac{1}{(n+1)!} \left(\frac{d}{dx} \right)^{n+1} k(s,x,0) \Big|_{x = 0}
\]
\end{proposition*}
\noindent \rmk \; When $M = \H^{n+1} / \Gamma$ for $\Gamma$ a coconvex subgroup, $k_{n+1} = 0$.
\subsubsection{$\overline{h}_{n+1}^{xx}$, even}
Expansion by minors of the inverse yields
\[
\overline{h}^{xx} = \frac{1}{\det \overline{h}} \det (\{\overline{h}_{ab}\})
\]
where $a, b = 1, \dots, n - 1$ are the coordinates corresponding to coordinates for $s$. Note that $\F(\bh_{ab}) = 1$ via
\begin{align*}
v_{a} &= \partial_{s_a} + u_a \partial_z \\
\overline{h}_{ab} &= \bg_{ab} + u_b \bg_{az} + u_a \bg_{bz} + u_a u_b \bg_{zz}
\end{align*}
%
Similarly, we have that
\begin{align*}
\det \overline{h} &= 1 + x^2 \overline{q}_2 + \dots + x^n \overline{q}_n + x^n \log(x) \overline{Q} + x^{n+1} \overline{q}_{n+1} + O(x^{n+2} \log(x))\\
\implies \frac{1}{\det \overline{h}} &= 1 + (\text{even terms up to order $n$}) - \overline{Q} x^n \log(x) - \overline{q}_{n+1} x^{n+1} + O(x^{n+2})
\end{align*}
Because $\F(\det (\{\overline{h}_{ab}\})) = 1$, we can compute
\begin{align*}
\overline{h}^{xx}_{n + 1} &= \left([\det \overline{h}]^{-1}\right)_{n+1} \cdot \left( \det (\{\overline{h}_{ab}\}) \right)_{0} + \left([\det \overline{h}]^{-1}\right)_{0} \cdot \left( \det (\{\overline{h}_{ab}\}) \right)_{n+1} \\
&= - \overline{q}_{n+1} + (\det \{\bh_{ab}\})_{n+1}
\end{align*}
where $\{\bh_{ab}\}$ is the submatrix of the metric for $TY$ but evaluated on just the $\{v_a, v_b\}$ basis (i.e. pushforwards of vectors parallel to the boundary curve). Note that $\F(\bh_{ab}) = 1$ for each $\bh_{ab}$ and $[\bh_{ab}]^{(0)} = \delta_{ab}$, i.e. to lowest order in $x$, this submatrix is the identity matrix. We can write this submatrix as
\begin{align*}
M &:= \{\bh_{ab}\} \\
&= Id + x^2 M_2 + \dots + x^n M_n + x^{n+1} M_{n+1} + O(x^{n+2}) \\
&= Id + R 
\end{align*}
where $R$ is even in $x$ up to order $n+1$ and $O(x^2)$. In particular, we see that $[R^k]^{n+1} = 0$ for any $k \geq 2$. Using the series expansion for the determinant of a matrix, we have
\begin{align*}
\det (M) &= \det(Id) + \text{Tr}(R) + O(R^2) \\
\implies [\det(M)]^{n+1} &= [\det(Id)]^{n+1} + [\text{Tr}_{Id}(R)]^{n+1} + [O(R^2)]^{n+1} \\
&= [\text{Tr}_{Id}(R)]^{n+1}
\end{align*}
So to find $(\det \{ \bh_{ab}\})_{n+1}$, it suffices to compute the $n+1$ part of its trace. Easily, we have 
\begin{align*}
\text{Tr}_{Id}(R) &= x^2 \text{Tr}_{Id}(M_2) + \dots + x^n \text{Tr}_{Id}(M_n) + x^{n+1} \text{Tr}_{Id}(M_{n+1}) + O(x^{n+2}) \\
\implies [\text{Tr}_{Id}(R)]^{n+1} &= \text{Tr}_{Id}(M_{n+1}) \\
&= \sum_{a = 1}^{n+1} [\bh_{aa}]^{n+1} 
\end{align*}
Recalling that
\begin{align*}
\overline{h}_{ab} &= \bg_{ab} + u_b \bg_{az} + u_a \bg_{bz} + u_a u_b \bg_{zz} \\
\F(\bg_{az}) &= \F(u_a) = 1 \\
u_a &= O(x^2) \\
\implies [\bh_{aa}]^{n+1} &= [\bg_{aa}]^{n+1}
\end{align*}
here, we recall that the metric coefficients $\bg_{ab}$ are be evaluated at $z = u(s,x)$ with their definitions in section \S \ref{MetricTM}.
%
To compute $\bg_{aa}$ we first recall the expansion of $\bg$ from the expansion for odd dimensional P.E. metrics in equation \eqref{kEquation}
\begin{align*}
\bg &= dx^2 + k_0(s,z) + x^2 k_2(s,z) + \dots + k_n(s,z) x^n + K(s,z) x^n \log(x) + O(x^{n+2}) \\
k_i(s,z) &= k_{i,0}(s) + k_{i,1}(s) z + \dots + k_{i,n/2}(s) z^{n/2} + O(z^{n/2 + 1}) \\
K(s,z) &= K_{0}(s) + k_{1}(s) z + \dots + k_{n/2}(s) z^{n/2} + O(z^{n/2 + 1}) 
\end{align*}
In particular, we note from remark \ref{NoKn+1Term} that there is no $x^{n+1}$ term when $n + 1$ is odd. Evaluating at $z = u = O(x^2)$ with $\F(u) = 1$, 
\begin{align*}
[\bg_{aa}]^{n+1} &= [k_0(s,z = u)(\partial_{s_a}, \partial_{s_a})]^{n+1}
\end{align*}
Of course, from \eqref{MetricCoeffs}
\[
k_{0}(s,z)(\partial_{s_a}, \partial_{s_b}) = \delta_{ab} + [\overline{\bGamma}_{ab}^z + \overline{\bGamma}_{ba}^z] z + O(z^2, sz, s^2)
\] 
where $\overline{\bGamma}_{\cdot \cdot}^{\cdot}$ are the Christoffels as in equation \eqref{MetricCoeffs}. When evaluated at $z = u$ and $s = 0$, the $O(z^2, sz, s^2)$ terms will not contribute a $n+1$st term. Thus 
\[
[k_0(s,z = u)]^{n+1} = -2 \bGamma_{\gamma, aaz} u_{n+1}
\] 
where $\bGamma_{\gamma, \cdot}$ denotes the restriction of the christoffels in lemma \ref{ChristoffelCylinderLemma} restricted to $T \gamma$. In sum,
\begin{align*}
[\bh_{aa}]^{n+1}&= [\bg_{aa}]^{n+1} = -2 \bGamma_{\gamma, aaz} u_{n+1}
\end{align*}
With this, we get
\begin{align*}
\left( \det (\{\overline{h}_{ab}\}) \right)_{n+1} & = \sum_{a = 1}^{n-1} -2 u_{n+1} \bGamma_{\gamma, aaz} \\
& = -2 u_{n+1} H_{\gamma, k_0}
\end{align*}
Moreover, using the fact that $H_{\gamma, k_0} = 2 (n-1) u_2$ (see lemma \ref{u2Lemma}), we get
\[
\boxed{\overline{h}^{xx}_{n+1} = -\overline{q}_{n+1} - 4(n-1) u_2 u_{n+1}}
\]
\subsubsection{$\overline{h}_{n+1}^{xx} + \overline{q}_{n+1}$ for $n$ odd}
Via the same reasoning as in the even case, we have that 
\[
\overline{h}_{n+1}^{xx} + \overline{q}_{n+1} = (\det \{\overline{h}_{ab}\})_{n+1} + R(u_2)
\]
where $R$ is some function of $u_2$ and its derivatives. This time
\begin{align*}
M &:= \{\bh_{ab}\} \\
&= Id + x^2 M_2 + \dots + x^{n-1} M_{n-1} + x^{n+1} M_{n+1} + O(x^{n+2}) \\
&= Id + R 
\end{align*}	
Note that from the expansions for $\bh_{ab}$, we see that $M_{2k}$ have coefficients which are dependent on $u_2(s)$ and its derivatives for $2k \leq n-1$. Thus
\begin{align*}
[\det (\{\bh_{ab}\})]_{n+1} &= \text{Tr}(M_{n+1}) + R(u_2)
\end{align*}
Note that because $n+1$ is now even, we have
\begin{align*}
\bg &= dx^2 + k_0(s,z) + x^2 k_2(s,z) + \dots + k_n(s,z) x^n + k_{n+1} x^{n+1} + O(x^{n+2}) \\
k_i(s,z) &= k_{i,0}(s) + k_{i,1}(s) z + \dots + k_{i,n/2}(s) z^{n/2} + O(z^{n/2 + 1})
\end{align*}
Note that because $n+1$ is even, $k_{n+1}$ term factors into our computation. A similar analysis to the even case yields
\begin{align*}
\bh_{aa} &= \bg_{aa} + 2u_a \bg_{az} + u_a^2 \bg_{zz} \\
[\bg_{aa}]^{n+1} &= [k_0(s,z = u)(\partial_{s_a}, \partial_{s_a})]^{n+1} + k_{n+1,0}(\p_{s_a}, \p_{s_a}) + R(u_2) \\
&= -2 u_{n+1} \bGamma_{\gamma, aaz} + k_{n+1,0}(\p_{s_a}, \p_{s_a}) + R(u_2) \\
[u_a \bg_{az}]^{n+1} &= R(u_2) \\
[u_a^2 \bg_{zz}]^{n+1} &= R(u_2)
\end{align*}
Summing as in the even case, we have 
\begin{align*}
\left( \det (\{\overline{h}_{ab}\}) \right)_{n+1} &= -4(n-1) u_2 u_{n+1} +  \text{Tr}_{T\gamma}(k_{n+1,0}) + R(u_2) \\
\implies \overline{h}_{n+1}^{xx} + \overline{q}_{n+1} &= -4(n-1) u_2 u_{n+1} +  \text{Tr}_{T\gamma}(k_{n+1,0}) + R(u_2)
\end{align*}
%
%
\subsection{Second variation of Mean Curvature - Parity for $\ddot{S}$} \label{SecondVarDetails}
In this section, we compute the second variation of mean curvature for a family of minimal surfaces $\{Y_t\}$ and prove theorem \ref{SecondDerivParity} and proposition \ref{SecondVariationOfMC}
\begin{theorem*}
Let $\{Y_t\} \subseteq M^{n+1}$ be a family of minimal of $m$-dimensional minimal submanifolds. Let $Y = Y_0$ and $\overline{h}$ denote a compactified metric on $Y$. Then for
\[
Y_t = \{ \exp_{\overline{h}, p}(S_t(p) \bnu(p)) \; | \; p \in Y\}
\]
and $\{w_i\}$ the normal basis described in section \S \ref{NormalFrame}, we have 
\[
\ddot{S} = \frac{d^2}{dt^2} \Big|_{t = 0} S_t = \ddot{S}^i w_i
\]
and $\F(\ddot{S}^i) = 1$.
\end{theorem*}
\begin{proposition*} 
Let $\{Y_t^m\} \subseteq M^{n+1}$ be a family of minimal of $m$-dimensional minimal submanifolds. Let $Y = Y_0$ and $\overline{h}$ denote a compactified metric on $Y$. Then for
\[
Y_t = \{ \exp_{\overline{h}, p}(S_t(p)) \; | \; p \in Y\}
\]
The second variation of mean curvature is given by
\[
\frac{d^2}{dt^2} H_t = J_Y^{\perp}(\ddot{S}) + Q^{\perp}(\dot{S})
\]
where $Q^{\perp}$ is a quadratic differential functional in $\dot{S}$ and 
\begin{align*}
Q^{\perp}(\dot{S}) &= Q^i(s,x) w_i \\
\F(Q^i) &= 1
\end{align*}
\end{proposition*}
\noindent We first show that $\ddot{S} \in NY$, we then sketch the proof of how one computes $Q^{\perp}(\dot{S})$ in codimension $1$. 

\subsubsection{Normality of $\ddot{S}$}
We first show that $\ddot{S} = \nabla_{F_t} F_t \Big|_{t = 0} \in NY$ since $F_t$ is a normal variation. Recall that the image of $\sigma_p(t) = F(p,t)$ is a geodesic curve starting at $p$. We write
\[
F_{*}(\partial_t) \Big|_{q = \sigma_p(t)} = F_t\Big|_{q = \sigma_p(t)} = A(t) \tau(t) \Big|_{q = \sigma_p(t)}
\]
for $\tau(t)$ a unit normal tangent vector evaluated at $q = F_t(p)$ on the path produced by $F_t(p)$. We compute 
\begin{align*}
\ddot{S}(p) &= \nabla_{F_t} F_t \Big|_{t = 0} \\
&= (A(t) \tau(t))(A(t))\Big|_{t = 0} \; \tau(0) + A^2(0)\nabla_{\tau(t)} \tau(t) \Big|_{t = 0} \\
&= \dot{A} \tau 
\end{align*}
where the second term vanishes since $\nabla_{\tau} \tau = 0$ since $\tau$ is the tangent vector to a geodesic curve. But $\tau(0) \in NY$, so $\ddot{S} \in NY$.

\subsubsection{$\ddot{S}$ Computation in codimension $1$}
For brevity, we sketch the proof in codimension $1$. We have $\ddot{S} = (x^{-1}\ddot{\phi}) \nu$  and can let $\nu(t) = \nu(F(t,p))$ be a normal vector for $Y_t$ at the point $F(t,p)$. By abuse of notation, we will absorb the prefactor of $x^{-1}$ and denote $\ddot{\phi} = x^{-1} \ddot{\phi}$, converting to the proper normalization at the end. We now compute
\[
\nabla_{F_t} \nu(t) \Big|_{t = 0} = - \nabla^Y \dot{\phi} \in TY
\]
this follows since 
\[
g(\nu(t), \nu(t)) \equiv 1, \qquad g(\nu(t), \Fa) \equiv 0
\]
for all $t$. We show
\begin{proposition} \label{codim1SecondVarOfMC}
For $\{Y_t^n\} \subseteq M^{n+1}$ a family of minimal submanifolds and $\dot{S} = \dot{\phi}(s,x) \bnu$, $\ddot{S} = \ddot{\phi} \bnu$, the second variation of mean curvature is given by
\begin{align*}
\frac{d^2}{dt^2} H(t) \Big|_{t = 0} &= J_Y(\ddot{\phi}) + G\dot{(\phi}, \nabla \dot{\phi}, D^2 \dot{\phi}) = 0
\end{align*}
where $\F(\dot{\phi}) = \F(G(\dot{\phi}, \nabla \dot{\phi}, \dot{\phi})) = 1$.
\end{proposition}
\noindent \Pf We compute 
\begin{align*}
\dot{h}^{\alpha \beta} &= -2 \dot{\phi} A^{\alpha \beta} \\
\ddot{h}^{\alpha \beta} & = 4 \dot{\phi}^2 (A \circ A)^{\alpha \beta} + 2 \ddot{\phi} A^{\alpha \beta} - 2 \dot{\phi}^2 R_{\nu \nu}^{\alpha \beta} + 2 \dot{\phi}^{\alpha} \dot{\phi}^{\beta}
\end{align*}
A more lengthy computation shows that 
\begin{align*}
\dot{A}_{\alpha \beta} &= \nabla_{F_t} A_{\alpha \beta}(t) \\
&= [\dot{\phi} R(\nu, \va, \vb, \nu) + \dot{\phi}_{\alpha \beta} - \dot{\phi} (A \circ A)_{\alpha \beta}] \nu \\
& - A_{\alpha \beta} (\nabla^Y \dot{\phi})
\end{align*}
and also
\begin{align*}
h^{\alpha \beta} \ddot{A}_{\alpha \beta} &= [J_Y(\ddot{\phi}) + Q_1(\dot{\phi}, \dot{\phi}) + Q_2(\dot{\phi} \nu, \nabla^Y \dot{\phi}) + Q_3( \nabla^Y \dot{\phi}, \nabla^Y \dot{\phi})] \nu - 4 ||A||^2 \dot{\phi} \nabla^Y \dot{\phi} \\
Q_1(\dot{\phi}, \dot{\phi}) &= \dot{\phi}^2 T_1(\nu, \nu) - 4 \dot{\phi}^2 g(A_{(\cdot, \cdot)}, R(\nu, \cdot, \cdot, \nu)) \\
Q_2(\dot{\phi}, \nabla^Y \dot{\phi}) &= \Ric_Y(\dot{\phi} \nu, - \nabla^Y \dot{\phi}) \\
Q_3(\nabla^Y \dot{\phi}, \nabla^Y \dot{\phi}) &= 2 g (\tilde{A}(\nabla^Y \dot{\phi}), \nabla^Y \dot{\phi})
\end{align*}
when $\dot{\phi}$ is a Jacobi field. In sum, we compute 
\begin{align*}
\frac{d^2}{dt^2} H(t) \Big|_{t = 0} &= \ddot{h}^{\alpha \beta} A_{\alpha \beta} + 2 \dot{h}^{\alpha \beta} \dot{A}_{\alpha \beta} + h^{\alpha \beta} \ddot{A}_{\alpha \beta} \\
& = [4 \dot{\phi}^2 \langle A \circ A, A \rangle - 2\dot{\phi}^2 \langle R(\nu, \cdot, \nu, \cdot), A \rangle
+ 2 A(\nabla \dot{\phi}, \nabla \dot{\phi})]   \\
& + 4[\dot{\phi}^2 \langle R_{\nu, \cdot, \nu, \cdot}, A \rangle - \dot{\phi} \langle D^2 \dot{\phi}, A \rangle
+ \dot{\phi}^2 \langle A \circ A, A \rangle] \\
& + [J_Y(\ddot{\phi}) + Q_1(\dot{\phi}, \dot{\phi}) + Q_2(\dot{\phi} \nu, \nabla^Y \dot{\phi}) + Q_3( \nabla^Y \dot{\phi}, \nabla^Y \dot{\phi})]
\end{align*}
In particular, when $\dot{\phi} \nu$ is a Jacobi field, this equation is another proof that $\ddot{S} \in N(Y)$. We reframe this as
\[
\frac{d^2}{dt^2} H(t) \Big|_{t = 0} = J_Y(\ddot{\phi}) + G\dot{(\phi}, \nabla \dot{\phi}, D^2 \dot{\phi})
\]
where $\F(\dot{\phi}) = \F(G(\dot{\phi}, \nabla \dot{\phi}, \dot{\phi})) = 1$.
\bibliography{Renormalized_Volume_via_Riesz_Regularization}{}
\bibliographystyle{plain}
\end{document}